\documentclass{amsart}
\usepackage{bbm,amsmath,amsfonts,amssymb}

 \newtheorem{theorem}{Theorem}[section]
 \newtheorem{lemma}[theorem]{Lemma}
 \newtheorem{corollary}[theorem]{Corollary}
 \newtheorem{proposition}[theorem]{Proposition}
 \newtheorem{definition}[theorem]{Definition}
 
 \newtheorem{remark}[theorem]{Remark}
 
 \newcounter{figures}[section]

\newcommand{\bN}{{\mathbb N}}

\newcommand{\bC}{{\mathbb C}}
\newcommand{\bR}{{\mathbb R}}
\newcommand{\bL}{{\mathbb L}}

  \def\cC{\mathcal{C}}
  
  \def\cE{\mathcal{E}}

  \def\cH{\mathcal{H}}

  \def\cP{\mathcal{P}}

  \def\cX{\mathcal{X}}

\def\LL{{\bL}}
\def\dd{{\delta}}
\def\ddj{{\dd_j}}
\def\ddelta{{\delta}}
\def\XX{{\cX}}
\def\supp{{\rm{supp}\, }}

\def\eps{{\varepsilon}}
\def\LLam{\Gamma}
\def\TTheta{\Theta}
\def\ssigma{\sigma}
\def\sss{\sigma}
\def\TT{T}
\def\SSS{S}
\def\Id{{\rm Id}}
\def\ONE{{\mathbbm 1}}

\def\diam{{\rm diam \,}}

\def\gg{\sigma}
\def\sig{\sigma}
\def\ph{\varphi}
\def\lam{\lambda}
\def\lamm{\lambda}
\def\tt{\tau}
\def\uu{u}

\def\Css{{K(\sig)}}
\def\CPhi{{K(\sig_*)}}
\def\Cdiam{{c^\diamond}}
\def\Cstar{{c^\star}}
\def\RRe{{\rm Re}\,}
\def\ww{w}

\def\bb{{\it b}}

\def\R{\mathbb{R}}

\def\Lp{{\LL^p}}
\def\til{\tilde}
\def\Dom{D}
\def\UCB{{\rm UCB}}
\def\D{D}

\def\W{w_{\alpha, \beta}}
\def\lip{{\rm Lip}\,}
\def\dst{{\delta_\star}}

\begin{document}

\title[Heat kernel generated frames]
{Heat kernel generated frames in the setting of Dirichlet spaces}

\author{T. Coulhon}
\address{Equipe AGM, CNRS-UMR 8088,
Universit\'{e} de Cergy-Pontoise, Cergy-Pontoise, France}
\email{thierry.coulhon@u-cergy.fr}

\author{G. Kerkyacharian}
\address{
Laboratoire de Probabilit\'{e}s et Mod\`{e}les Al\'{e}atoires, CNRS-UMR 7599,
Universit\'{e} Paris VI et Universit\'{e} Paris VII, rue de Clisson, F-75013 Paris}
\email{kerk@math.jussieu.fr}

\author{P. Petrushev}
\address{Department of Mathematics\\University of South Carolina\\
Columbia, SC 29208}
\email{pencho@math.sc.edu}



\thanks{Corresponding author: Pencho Petrushev,
Email: pencho@math.sc.edu,
Phone:+803-777-6686,
Fax:+803-777-6527}

\date{May 10, 2012}

\begin{abstract}
Wavelet bases and frames consisting of band limited functions of nearly exponential localization on $\R^d$
are a powerful tool in harmonic analysis by making various spaces of functions
and distributions more accessible for study and utilization,
and providing sparse representation of natural function spaces (e.g. Besov spaces) on $\R^d.$
Such frames are also available on the sphere and in more general homogeneous spaces,
on the interval and ball.
The purpose of this article is to develop band limited well-localized frames
in the general setting of Dirichlet spaces with doubling measure and
a local scale-invariant Poincar\'e inequality
which lead to heat kernels with  small time Gaussian bounds and H\"{o}lder continuity.
As an application of this construction, band limited frames are developed in the context of
Lie groups or homogeneous spaces with polynomial volume growth,
complete Riemannian manifolds with Ricci curvature bounded from below and
satisfying the volume doubling property, and other settings.
The new frames are used for decomposition of Besov spaces in this general setting.
\end{abstract}

\maketitle

\tableofcontents


\section{Introduction}\label{introduction}
\setcounter{equation}{0}

Decomposition systems (bases or frames) consisting of band limited functions of nearly
exponential space localization have had significant impact in
theoretical and computational harmonic analysis,
PDEs, statistics, approximation theory and their applications.
Meyer's wavelets \cite{M} and the frames (the $\varphi$-transform)
of Frazier and Jawerth \cite{F-J1, F-J2, F-J-W} are the most
striking examples of such decomposition systems playing a pivotal role in the solution of
numerous theoretical and computational problems.
The key to the success of wavelet type bases and frames is rooted in their ability to capture a great deal
of smoothness and other norms in terms of respective coefficient sequence norms and
provide sparse representation of natural function spaces (e.g. Besov spaces) on $\bR^d.$
Frames of a similar nature have been recently developed in non-standard settings such as
on the sphere \cite{NPW1,NPW2} and more general homogeneous spaces
\cite{PesG}, on the interval \cite{PX1,KPX1} and ball \cite{PX2, KPX2} with weignts,
and extensively used in statistical applications (see e.g. \cite{KKPPW, KKPPP}).

The primary goal of this paper is to extend and refine the construction of band limited frames
with elements of nearly exponential space localization to the general setting
of strictly local regular Dirichlet spaces with doubling measure and  local scale-invariant Poincar\'e inequality
which lead to a markovian heat kernel with small time Gaussian bounds and H\"{o}lder continuity.
The key point of our approach is to be able to deal with
(a) different geometries,
(b) compact and noncompact spaces,
and
(c) spaces with nontrivial weights,
and at the same time to allow for the development and frame decomposition of
Besov and Triebel-Lizorkin spaces with complete range of indices.
This will enable us to cover and shed new light on the existing frames and space decompositions
and develop band limited localized frames in the context of
Lie groups or homogeneous spaces with polynomial volume growth,
complete Riemannian manifold with Ricci curvature bounded from below and satisfying the volume doubling condition,
and other new settings.
To this end we will make advances on several fronts:
development of functional calculus of positive self-adjoint operators with associated heat kernel
(in particular, localization of the kernels of related integral operators),
development of lower bounds on kernel operators,
development of a Shannon sampling theory,
Littlewood-Paley analysis, and
development of dual frames.

As a first application of our frames we shall develop rapidly and characterize
the classical Besov spaces $B^s_{pq}$ with positive smoothness and $p\ge 1$.
Classical and nonclassical Besov and Triebel-Lizorkin spaces in the general framework of this paper
with full range of indices
and their frame and heat kernel decompositions are developed in the follow-up paper \cite{KP}.

In this preamble we outline the main components and points of this undertaking, including
the underlying setting, a general scenario for realization of the setting and examples,
and a description of the main results.


\subsection{The setting}\label{sec:setting}
We now describe precisely all the ingredients we need to develop our theory.

I. We assume that $(M,\rho,\mu)$ is a metric measure space, which satisfies the conditions:

(a) $(M, \rho)$ is a locally compact  metric space with distance $\rho(\cdot, \cdot)$
and $\mu$ is a~positive Radon measure
such that the following {\em volume doubling condition} is valid
\begin{equation}\label{doubling}
0 < \mu(B(x,2r) ) \leq 2^d \mu(B(x,r))<\infty
\quad\hbox{for $x \in M$ and $r>0$.}
\end{equation}
Here $B(x,r)$ is the open ball centered at $x$ of radius $r$ and
$d >0$ is a constant that plays the role of a~dimension.
Note that $(M,\rho,\mu)$ is also a homogeneous space in the sense of Coifman and Weiss \cite{CW1, CW2}.

(b) The {\em reverse doubling condition} is assumed to be valid, that is, there exists a constant
$\beta >0$ such that
\begin{equation}\label{reverse-doubling}
\mu(B(x,2r) ) \ge 2^\beta \mu(B(x,r))
\quad\hbox{for $x \in M$ and $0< r \le \frac{\diam M}{3}$.}
\end{equation}
It will be shown in \S2 that this condition is a consequence of the above doubling condition
if $M$ is connected.

(c) The following {\em non-collapsing condition} will also be stipulated:
There exists a constant $c>0$ such that
\begin{equation}\label{non-collapsing}
\inf_{x\in M}\mu(B(x,1) )\ge c, \quad x\in M.
\end{equation}
As will be shown in \S2 in the case $\mu (M)<\infty$
the above inequality follows by (\ref{doubling}).
Therefore, it is an additional assumption only when $\mu(M)=\infty$.

Since we consider in this paper inhomogeneous function spaces only,
it would be natural to make only purely local assumptions, and in particular to assume doubling only for balls
with radius bounded by some constant, which would enlarge considerably our range of examples.
This would require however more work and more space. On the other hand, our next assumptions
on the heat kernel are local, in the sense that they are required for small time only.
Clearly, by assuming global doubling and global heat kernel bounds,
one can treat homogeneous  spaces as well.


II. Our main assumption is that the local geometry of the space $(M,\rho,\mu)$ is related to
an essentially self-adjoint positive operator $L$ on $L^2(M, d\mu)$
such that the associated semigroup $P_t=e^{-tL}$ consists of integral operators with
(heat) kernel $p_t(x,y)$ obeying the conditions:

(d) {\it Small time  Gaussian upper bound}:
\begin{equation}\label{Gauss-local}
p_t(x,y) \le \frac{ C\exp \{-\frac{c\rho^2(x,y)}t\}}{\sqrt{\mu(B(x,\sqrt t))\mu(B(y, \sqrt t))}}
\quad\hbox{for} \;\;x,y\in M,\,0<t\le 1.
\end{equation}
One can see that by combining the results in  \cite{CCO, Ouhabaz} and \cite{CS}, this estimate
and the doubling condition (\ref{doubling}) coupled with the fact that
$e^{-tL}$ is actually a holomorphic semigroup on $L^2(M,d\mu)$, i.e. $e^{-zL}$ exists for $ z\in \bC$, $\RRe z \ge  0$,
imply that $e^{-zL}$ is an integral operator with kernel $p_z(x,y)$ satisfying the following estimate:
For any $z=t+iu$, $0<t\le 1$, $u\in\bR$, $x, y\in M$,
\begin{equation}\label{hol1}
|p_z(x,y)| \le \frac{C\exp\big\{-c \RRe \frac{\rho^2(x,y)}{z}\big\}}
{\sqrt{\mu(B(x,\sqrt t))\mu(B(y, \sqrt t))}}.
\end{equation}

(e) {\em H\"{o}lder continuity}: There exists a constant $\alpha>0$ such that
\begin{equation}\label{lip}
\big|  p_t(x,y) - p_t(x,y')  \big|
\le C\Big(\frac{\rho(y,y')}{\sqrt t}\Big)^\alpha
\frac{\exp\{-\frac{c\rho^2(x,y)}t \}}{\sqrt{\mu(B(x,\sqrt t))\mu(B(y, \sqrt t))}}
\end{equation}
for $x, y, y'\in M$ and $0<t\le 1$, whenever $\rho(y,y')\le \sqrt{t}$.

\smallskip

(f) {\it Markov property}:
\begin{equation}\label{hol3}
\int_M p_t(x,y) d\mu(y) \equiv 1
\quad\hbox{for $t >0$,}
\end{equation}
which readily implies, by analytic continuation,
\begin{equation}\label{Markov}
\int_M p_z(x,y) d\mu(y) \equiv 1
\quad\hbox{for $z=t+iu$,\; $t >0$.}
\end{equation}
Above $C,c>0$ are structural constants that will affect almost all constants in what follows.

The main results in this article will be derived from the above conditions.
However, it is perhaps suitable to exhibit a more tangible general scenario
that guarantees the validity of these conditions.

\subsection{Realization of the setting in the framework of Dirichlet spaces}\label{sec:Dirichle-spaces}

We would like to point out that in a general framework of Dirichlet spaces
the needed Gaussian bound, H\"{o}lder continuity, and Markov property of the heat kernel follow from
the local scale-invariant Poincar\'e inequality and the doubling condition on the measure,
which in turn are equivalent to the parabolic Harnack inequality.
The point is that situations where our theory is applicable are quite common and it just amounts
to verifying the local scale-invariant Poincar\'e inequality and the doubling condition on the measure.
We shall further illustrate this point on several examples and, in particular,
on the ``simple" example of $[-1, 1]$ with the heat kernel induced by the Jacobi operator,
seemingly not covered in the literature.

We shall operate in the framework of strictly local regular Dirichlet spaces
(see \cite{FUKU,Ouhabaz,ALB,Sturm0, Sturm1,Sturm, BM, BM1, Davies}).
To be specific, we assume that $M$ is a locally compact separable metric space equipped with
a positive Radon measure $\mu$ such that every open and nonempty set has positive measure.
Also, we assume that $L$ is a~positive symmetric operator on (the real) $L^2(M, \mu)$ with domain $D(L)$,
dense in  $L^2(M, \mu)$. We shall denote briefly $\LL^p:= L^p(M, \mu)$ in what follows.
One can associate with $L$ a symmetric non-negative  form
$$
\cE (f,g) = \langle Lf, g \rangle = \cE(g,f), \;
\cE (f,f) = \langle Lf, f \rangle \geq 0,
$$
with domain $D(\cE)= D(L)$.
We consider on $D(\cE)$ the prehilbertian structure induced by
$$\|f\|_{\cE}^2 = \|f \|_2^2 + \cE(f,f)$$
which in general is not complete (not closed),
but closable (\cite{Davies}) in $\bL^2$.
Denote by $\overline{\cE}$ and $D(\overline{\cE})$ the closure of $\cE$ and its domain.
It gives rise to a self-adjoint extension  $\overline{L}$ (the Friedrichs extension) of $L$ with domain
$D(\bar L)$ consisting of all $f \in D(\overline{\cE})$ for which there exists
$u \in \bL^2$ such that $\overline{\cE}(f,g)= \langle u, g \rangle$ for all $g \in D(\overline{\cE})$
and $\bar Lf =u$.
Then $\bar L$ is positive and self-adjoint, and
$$
D(\overline \cE)= D((\overline L)^{1/2}), \quad
\overline{\cE} (f,g) = \langle (\overline L)^{1/2} f,  (\overline L)^{1/2}g \rangle.
$$

Using the classical spectral theory of positive self-adjoint operators, we can associate with $\overline L$
a self-adjoint strongly continuous contraction semigroup $ P_t =e^{-t \bar L}$ on $\bL^2(M,\mu)$.
Then
$$e^{-t \bar L}= \int_0^\infty e^{-\lambda t} dE_\lambda,$$
where $E_\lambda$ is the spectral resolution associated with $\bar L$.
Moreover this semigroup has a holomorphic extension to the complex half-plane $\RRe z  >0.$

\smallskip

Our next assumption is that $P_t$ is a submarkovian semigroup:
$0\le f\le 1$ and $f \in \bL^2$ imply $0\le P_t f\le 1$.
Then $P_t$ can be extended as a contraction operator on
$ \bL^p$, $1\leq p\leq \infty$, preserving positivity, satisfying $ P_t 1 \leq 1$,
and hence yielding a strongly continuous contraction  semigroup on $\bL^p, 1\le p<\infty$.
A sufficient condition for this \cite{ALB,FUKU}, which can be verified on $D(L)$, is that
for every $\varepsilon >0$ there exists
$\Phi_\varepsilon: \bR \mapsto [-\varepsilon, 1+ \varepsilon]$
such that $\Phi_\varepsilon$ is non-decreasing, $\Phi_\varepsilon\in {\rm Lip}\, 1$,
$ \Phi_\varepsilon (t)=t$ for $t \in [0,1]$ and
$$
\Phi_\varepsilon(f) \in D(\overline{\cE})\;\;
\hbox{and}\;\;
\overline{\cE}( \Phi_\varepsilon(f),  \Phi_\varepsilon(f)) \leq \cE(f,f),
\quad \forall f \in D(L).
$$
(in fact, this can be done easily only if $ \Phi_\varepsilon (f) \in D(L)$).

Under the above conditions, $(D(\overline{\cE}), \overline{\cE})$ is called a Dirichlet space and
$D(\overline{\cE})\cap \bL^\infty$ is an algebra.

We assume that the form $\overline{\cE}$ is strongly local,
i.e. $\overline{\cE}(f,g)=0$ for $f,g\in D(\overline{\cE})$
whenever $f$ is with compact support and $g$ is constant on a neighbourhood of the support of $f$.
We also assume that $\overline{\cE}$ is regular, meaning that
the space $\cC_c(M)$ of continuous functions on $M$ with compact support has the property that
the algebra $\cC_c(M)\cap D(\overline{\cE})$ is dense in $\cC_c(M)$ with respect to the sup norm,
and dense in $D(\overline{\cE})$ in the norm $\sqrt{\overline{\cE}(f,f)+\|f\|_2^2}$.

We next give a sufficient condition for strong locality and regularity (\cite{FUKU}, Chapter 3)
which can be verified for $D(L)$: $\overline{\cE}$ is strongly local and regular if
(i) $ D(L)$ is a subalgebra of $\cC_c(M)$ verifying the strong local condition:
$0= \cE(f,g) = \langle Lf,g \rangle$ if $ f, g \in D(L)$,
$f$ is with compact support, and $g$ is constant on a neighbourhood of the support of $f$,
and
(ii) for any compact $K$ and open set $U$ such that $K \subset U$ there exists
$u \in D(L)$, $u\ge 0$, $\supp u \subset U$, and $u \equiv 1$ on $K$
(thus $D(L)$ is a dense subalgebra of $\cC_c(M)$  and dense in $ D(\overline{\cE})$).

\smallskip

Under the above assumptions, there exists a bilinear symmetric form $d\Gamma$  defined on
$D (\overline{\cE}) \times D(\overline{\cE})$ with values in the signed Radon measures on $M$
such that
$$
\cE( \phi f,g) +\cE(f, \phi g) -\cE( \phi, fg)= 2 \int_M \phi d\Gamma(f,g)
\quad\hbox{for}\quad f, g, \phi \in \cC_c(M)\cap D(\overline{\cE}),
$$
which obviously verifies
$\overline{\cE}(f,g)=\int_Md\Gamma(f,g)$ and  $d\Gamma(f,f) \geq 0$.

In fact, if $D(L)$ is a subalgebra of $\cC_c(M)$, then $d\Gamma $ is absolutely continuous
with respect to $\mu,$ and
$$
d\Gamma(f,g)(u) = \Gamma(f,g)(u) d\mu(u), \quad
\Gamma(f,g) = \frac 12 (L(fg ) -f Lg -g Lf) \quad
\forall f,g \in D(L).
$$
In other words, $\overline{\cE}$ admits a "carr\'e du champ" (\cite{BH}, Chapter 1, \S 4):
There exists a bilinear function
$
D (\overline{\cE} ) \times D (\overline{\cE})\ni f, g \mapsto \Gamma(f,g) \in \LL^1
$
such that $\Gamma(f,f)(u) \geq 0$,
$$
\overline{\cE}( \phi f,g) +\overline{\cE}(f, \phi g) - \overline{\cE}( \phi, fg)
= 2 \int_M \phi (u)\Gamma(f,g)(u) d\mu(u)
\quad \forall f,g, \phi \in  D(\overline{\cE}) \cap \bL^\infty,
$$
and
$
\overline{\cE}(f,g) = 2 \int_M \Gamma(f,g)(u) d\mu(u).
$

\smallskip

One can define an intrinsic distance on $M$ by
$$
\rho(x,y)=\sup\{u(x)-u(y): u\in D(\overline{\cE})\cap\cC_c(M),
d\Gamma(u,u)= \gamma(u)(x) d\mu (x), \gamma(u)(x) \leq 1\}.
$$
We assume that $\rho: M\times M\to [0,\infty]$
is actually a true metric that generates the original topology on $M$
and that $(M,\rho)$ is a complete metric space.

As a consequence of this assumption, the space $M$ is connected,
the closure of an open ball $B(x,r)$ is the closed ball
$\overline{B(x,r)} := \{y \in M, \rho(x,y) \leq r\}$,
and the closed balls are compact (see \cite{Sturm,Sturm0,Sturm1}).

\smallskip

We are now in a position to describe an {\em  optimal scenario} when the needed
Gaussian bound (\ref{Gauss-local}), H\"{o}lder continuity (\ref{hol1}), and Markov property (\ref{lip})
on the heat kernel can be effectively realized.
In the framework of strictly local regular Dirichlet spaces with a complete intrinsic metric,
the following two properties are equivalent \cite{Sturm, HS}:

(i) The heat kernel satisfies
\begin{equation}\label{Gauss-local-bis}
\frac{ c_1'\exp \{-\frac{c_1\rho^2(x,y)}t\}}{\sqrt{\mu(B(x,\sqrt t))\mu(B(y, \sqrt t))}}
\le  p_t(x,y)
\le \frac{ c_2'\exp \{-\frac{c_2\rho^2(x,y)}t\}}{\sqrt{\mu(B(x,\sqrt t))\mu(B(y, \sqrt t))}}
\end{equation}
for $x,y\in M$ and $0<t\le 1$.

(ii)
(a)
$(M,\rho,\mu)$ is a local doubling measure space:
There exists $d>0$ such that $\mu(B(x,2r)) \le 2^d \mu(B(x,r)$ for $x\in M$ and $0<r<1$.

(b) Local scale-invariant Poincar\'e inequality holds:
There exists a constant $C>0$ such that for any ball $B=B(x,r)$ with $0<r\le 1$, $x\in M$,
and any function $f\in D(\overline{\cE})$,
$$
\int_B|f-f_B|^2\le Cr^2 \int_Bd\Gamma(f,f).
$$
Moreover, it is also well-known that the above property is equivalent to a local parabolic Harnack inequality,
and, furthermore, any of these equivalent properties implies the validity of
(\ref{lip}) and (\ref{hol3}) (see \cite{Sturm}, \cite{HS}, \cite{S}, \cite{GS}, and the references therein).

Consequently, given a situation which fits into the framework of strictly local regular
Dirichlet spaces with a complete intrinsic metric it suffices to only verify
the local Poincar\'e inequality and the global doubling condition on the measure and
then our theory applies in full.

\smallskip

In a future work we will further develop this theory under the more general assumption of the small time
sub-Gaussian estimate:
\begin{equation}\label{subGauss-local}
p_t(x,y) \le \frac{ C\exp \left\{-c\big(\frac{\rho^{m}(x,y)}t\big)^{\frac{1}{m-1}}\right\}}{\sqrt{\mu(B(x,t^{1/m}))\mu(B(y, t^{1/m}))}}
\quad\hbox{for} \;\;x,y\in M,\,0<t\le 1,
\end{equation}
where $m\ge 2$.

\subsection{Examples}

There is a great deal of set-ups which fit in the general framework of this article.
We next briefly describe several benchmark examples which are indicative for the versatility and depth of our methods.

\subsubsection{Uniformly elliptic divergence form operators on $\R^d$}

Given a uniformly elliptic symmetric matrix-valued function $\{a_{i,j}(x)\}$ depending on $x\in\R^d$,
one can define an operator
$$L = -\sum_{i,j=1}^d\frac{\partial}{\partial x_i}\left(a_{i,j}\frac{\partial}{\partial x_j}\right)$$
on $L^2(\R^d,dx)$  via the associated quadratic form.
Thanks to the uniform ellipticity condition, the intrinsic metric associated with this operator
is equivalent to the Euclidean distance.
The Gaussian upper and lower estimates of the heat kernel  in this setting hold for all time
and are due to Aronson, the H\"older regularity of the solutions is due to Nash \cite{N},
the Harnack inequality was obtained by Moser \cite{MO, MO1}.

\subsubsection{Domains in $\R^d$}

One can define  uniformly elliptic divergence form operators on $\R^d$ by choosing boundary conditions.
In this case the upper bounds of the heat kernels are well understood (see for instance \cite{Ouhabaz}).
The problem for establishing Gaussian lower bounds is much more complicated. One has to choose Neumann conditions
and impose regularity assumptions on the domain. For the state of the art, we refer the reader to \cite{GS}.

\subsubsection{Riemannian manifolds and Lie groups}

The conditions from \S\ref{sec:Dirichle-spaces} are verified for
the Laplace-Beltrami operator of a Riemannian manifold with non-negative Ricci curvature \cite{LY},
also for manifolds with Ricci curvature bounded from below if one assumes in addition that they satisfy
the volume doubling property,
also for manifolds that are quasi-isometric to such a manifold \cite{G1, SD, Parma},
also for co-compact covering manifolds whose deck transformation group has polynomial growth \cite{SD, Parma},
for sublaplacians on polynomial growth Lie groups \cite{VSC, Robinson} and their homogeneous spaces \cite{M}.
We would like to point out that the case of the sphere endowed with the natural Laplace-Beltrami operator
treated in \cite{NPW1,NPW2} and the case of more general compact homogeneous spaces endowed with
the Casimir operator considered in \cite{PesG} fall into the above category.
One can also consider variable coefficients operators on Lie groups, see  \cite{SStr}.

We refer the reader to \cite[Section 2.1]{GS} for further details on the above examples.
For more references on the heat kernel in various settings, see \cite{VSC, Davies, Grig, S}.

\subsubsection{Heat kernel on $[-1, 1]$ generated by the Jacobi operator}\label{subsec:Jacobi}

To show the flexibility of our general approach to 
frames and spaces through heat kernels we consider in \S\ref{sec:Jacobi} the ``simple" example of
$M =[-1, 1]$ with $d\mu(x) = \W(x)dx$, where $\W(x)$ is the classical Jacobi weight:
$$
\hbox{
$\W(x)=w(x)= (1-x)^\alpha (1+x)^\beta$, \quad $\alpha , \beta >-1$.
}
$$
The Jacobi operator is defined by
$$
Lf(x) = - \frac{\big[w(x) a(x)  f'(x)\big]'}{w(x)}
\quad\hbox{with} \;\; a(x):=1-x^2
$$
and
$D(L) = C^2[-1,1]$.
As is well-known (\cite{Sz}),
$LP_k =\lambda_k P_k$, where $P_k$ ($k\ge 0$) is the $k$th degree (normalized) Jacobi polynomial and
$\lambda_k = k(k+\alpha + \beta +1)$.
Integration by parts gives
$$
\cE(f,g) :=\langle Lf, g \rangle = \int_{-1}^1 a(x) f'(x) g'(x) \W(x) dx.
$$
In \S\ref{sec:Jacobi} it will be shown that in this case the general theory applies,
resulting in a complete strictly local Dirichlet space with an intrinsic metric defined by
$$
\rho(x, y)= |\arccos x - \arccos y|,\quad x,y \in [-1, 1],
$$
which is apparently compatible with the usual topology on $[-1, 1]$.
It will be also shown that in this setting the measure $\mu$ verifies the doubling condition
and the respective scale-invariant Poincar\'e inequality is valid.
Therefore, the example under consideration fits in the general setting described in \S\ref{sec:Dirichle-spaces}
and our theory applies.
In particular, the associated heat kernel with a representation
$$
p_t(x, y) = \sum_{k\ge 0} e^{-\lambda_k t}P_k(x)P_k(y)
$$
has Gaussian bounds (see \S\ref{sec:Jacobi}),
which to the best of our knowledge appears first in the present article.
Another consequence of this is that our theory covers completely the construction of frames and the development of
Besov and Triebel-Lizorkin spaces on $[-1, 1]$ with Jacobi weights from \cite{KPX1, PX1}.

Finally, we would like to point out that there are other examples,
e.g. the development of frames and weighted Besov and Triebel-Lizorkin spaces on the unit ball $B$ in $\R^d$
in \cite{KPX2, PX2},
which perfectly fit in our general setting, but we will not pursue in this article.

\subsection{Outline of the paper}

This paper is organized as follows:

In \S2 we give some auxiliary results which are instrumental in proving our main results.
In particular, we collect all needed facts about
doubling measures and related kernels,
construction of maximal $\delta$-nets, and
integral operators.

In \S3 we develop some components of a non-holomorphic functional calculus related to a positive self-adjoint
operator $L$ in the general set-up of the paper. In particular, we establish the nearly exponential localization
of the kernels of operators of the form $f(\sqrt L)$ under suitable conditions on $f$.
These localization results are crucial for the development of the Littlewood-Paley theory in our setting.
They also enable us to explore the main properties of the spectral spaces and develop the linear
approximation theory from spectral spaces through the machinery of Jackson-Bernstein inequalities and interpolation.
In this section we also give the main properties of finite dimentional spectral spaces.

In \S4 we establish a sampling theorem in the spirit of the Shannon theory
and develop a cubature rule/formula in the compact and non-compact case,
which is exact for spectral functions of a given order.
This cubature rule is a critical component in the development of our frames.

Our main results are placed in \S5, where we construct pairs of dual frames of the form:
$\{\psi_{j\xi}: \xi\in\XX_j, j\ge 0\}$, $\{\tilde\psi_{j\xi}: \xi\in\XX_j, j\ge 0\}$,
where each $\XX_j$ is a $\delta_j$-net on $M$ for an appropriate $\delta_j$.
The frame elements $\psi_{j\xi}$, $\tilde{\psi}_{j\xi}$ are band limited and well-localized functions,
which allow for decomposition of functions and distributions from various spaces
(in particular, Besov and Triebel-Lizorkin spaces)
of the form
$$
f= \sum_{j\ge 0} \sum_{\xi \in \XX_j} \langle f, \tilde\psi_{j\xi}\rangle \psi_{j\xi}.
$$
The most critical point in this paper is the construction of the dual frame $\{\tilde\psi_{j\xi}\}$.
We develop it in two settings:
(i) in the general case, and
(ii) in the case when the spectral spaces have the polynomial property under multiplication (see \S5.3).
In the second case the construction is simple and elegant, however, the setting is somewhat restrictive,
while in the first case the construction is much more involved, but the localization of $\tilde\psi_{j\xi}$
is inverse polynomial of an arbitrarily fixed order.

In \S6 we develop the classical and most commonly used Besov spaces $B^s_{pq}$
with indices $s>0$, $1\le p\le \infty$, and $0<q\le \infty$ in the setting of this paper.
These spaces are defined through Littlewood-Paley decomposition and characterized as approximation spaces
of linear approximation from spectral spaces. A frame decomposition of $B^s_{pq}$ is also established.
In full generality,  classical and non-classical Besov and Triebel-Lizorkin spaces and their frame
decomposition in the general setting of the paper are developed in \cite{KP}.

Section 7 is an appendix, where we place the proofs of the Poincar\'e inequality for the Jacobi operator
and the doubling property of the respective measure. Gaussian bounds of the associated heat
kernel are also established.

\smallskip

\noindent
{\bf Notation.}
Throughout this article we shall use the notation
$|E|:= \mu(E)$ for $E\subset M$,
$\LL^p:=L^p(M, \mu)$,
$\|\cdot\|_p:=\|\cdot\|_{\LL^p}$,
and $\|T\|_{p\to q}$ will denote the norm of a bounded operator $T: \LL^p \to \LL^q$.
$\UCB$ will stand for the space of all uniformly continuous and bounded functions on $M$
and $\LL^\infty$ will be in most cases identified with $\UCB$.
$D(T)$ will stand for the domain of a given operator $T$.
We shall denote by $C^\infty_0(\bR_+)$ the set of all compactly supported $C^\infty$ functions
on $\bR_+:=[0, \infty)$.
In most cases ``$\sup$" will mean ``${\rm ess \,sup}$".
Positive constants will be denoted by $c$, $C$, $c_1$, $c'$, $\dots$ and they may vary
at every occurrence, $a\sim b$ will stand for $c_1\le a/b\le c_2$.

\section{Doubling metric measure spaces: Basic facts}\label{conseq-doubling}
\setcounter{equation}{0}

In this section we put together some simple facts
related to metric measure spaces $(M, \rho, \mu)$ obeying
the doubling, inverse doubling and non-collapsing conditions
(\ref{doubling})-(\ref{non-collapsing})
and integral operators acting on functions defined on such spaces.

\subsection{Consequences of doubling and clarifications} 

The~doubling condition (\ref{doubling}) readily implies
\begin{equation}\label{D1}
 |B(x,\lambda r)| \le (2\lambda)^d  |B(x, r)|,
\quad x\in M, \;\lambda > 1, \; r>0,
\end{equation}
and, therefore, due to $B(x,r) \subset B(y, \rho(y,x) +r)$
\begin{equation}\label{D2}
|B(x, r)| \le 2^d \Big(1+ \frac{\rho(x,y)}{r}\Big)^d  |B(y, r)|,
\quad x, y\in M, \; r>0.
\end{equation}
In turn, the reverse doubling condition yields
\begin{equation}\label{D3}
|B(x,\lambda r)| \ge (\lambda/2)^\beta |B(x, r)|,
\quad \lambda > 1, \; r>0,\; \hbox{$0< \lambda r<\frac{\diam M}{3}$.}
\end{equation}
Also, the non-collapsing condition (\ref{non-collapsing}) coupled with (\ref{D1}) implies
\begin{equation}\label{est-muB-unify1}
\inf_{x\in M}|B(x, r)| \ge \hat{c}r^d,
\quad 0< r \le 1,
\end{equation}
where $\hat{c}=c2^{-d}$ with $c>0$ the constant from (\ref{non-collapsing}).

Note that $|B(x, r)|$ can be much larger than $cr^d$ as is evidenced by the case of the Jacobi operator
on $[-1, 1]$, considered in \S\ref{subsec:Jacobi} and \S\ref{sec:Jacobi},  see (\ref{measure-ball}).

Several clarifying statements are in order.
We begin with a claim which, in particular, shows that the non-collapsing condition is automatically
obeyed when $\mu(M) <\infty$.


\begin{proposition}\label{prop:prop-doubl-space}
Let $(M, \rho, \mu)$ be a metric measure space which obeys
the doubling condition $(\ref{doubling})$. Then

$(a)$ $\mu(M) <\infty$ if and only if $\diam M <\infty$.
Moreover, if $\diam M =D <\infty$, then
\begin{equation}\label{est-muB}
\inf_{x\in M}|B(x, r)| \ge r^d|M|(2D)^{-d},
\quad  0< r \leq D.
\end{equation}

$(b)$
$\mu (\{x\}) >0$ for some $x\in M$ if and only if $\{x\}=B(x, r)$ for some $r>0$.

\end{proposition}

\noindent
{\bf Proof.}
We first prove (a).
Note that if $\diam M =D <\infty$, then
$M= B(x,D)$ for any $x\in M$ and hence $|M| =| B(x,D)| <\infty$.

In the other direction, let $|M|< \infty$. Assume on the contrary that $\diam M= \infty$.
Then inductively one can construct a sequence of points $\{x_0, x_1, \dots\} \subset M$ such that
if $d_j:= \rho(x_0, x_j)$, then $1\le d_1< d_2< \cdots$ and
$d_{j+1}>3d_j$, $j\ge 0$.
One checks easily that $B(x_j, \frac{d_j}{2}) \cap B(x_k, \frac{d_k}{2})=\emptyset$
if $j\ne k$.
On the other hand, using (\ref{doubling}),
$$
0<|B(x_0, 1)| \le |B(x_j, 2d_j)| \le 4^d |B(x_j, d_j/2)|.
$$
Therefore, we have a sequence of disjoint balls $\{B(x_j, \frac{d_j}{2})\}_{j\ge 1}$ in $M$
such that $|B(x_j, \frac{d_j}{2})| \ge 4^{-d}|B(x_0, 1)| >0$ and hence $|M|=\infty$.
This is a contradiction that proves the claim.

Estimate (\ref{est-muB}) is immediate from (\ref{D1}).


To prove (b), we first note that if $\{x\}=B(x, r)$ for some $r>0$, then (\ref{doubling})
implies $\mu(\{x\}) >0$.
For the other implication, let $ \mu(\{x\}) >0$ and assume that $\{x\} \not= B(x,r)$ for all $r>0$.
Then we use this to construct inductively a sequence $\{x_1, x_2, \dots\} \subset M$
such that if $d_j:= \rho(x, x_j)$, then $d_1> d_2> \cdots >0$ and
$d_{j+1}<\frac{d_j}{3}$, $j\ge 1$. Clearly, the latter inequality yields
$B(x_j, \frac{d_j}{2}) \cap B(x_k, \frac{d_k}{2})=\emptyset$
if $j\ne k$.
On the other hand by our assumption, (\ref{doubling}), and the fact that
$x\in B(x_j, 2d_j)$ we infer
$$
0<\mu(\{x\}) \le |B(x_j, 2d_j)| \le 4^d|B(x_j, d_j/2)|.
$$
Now, as above we conclude that $|M|=\infty$ which is a contradiction.
$\qed$

\smallskip

We next show that the reverse doubling condition (\ref{reverse-doubling}) is not quite restrictive.


\begin{proposition}\label{prop:rev}
If $M$ is connected, then the reverse doubling condition holds,
i.e. there exists $\beta >0$ such that
$$
|B(x,2r) | \ge 2^\beta |B(x,r)|
\quad\hbox{for $x\in M$ and $0<r<\frac{\diam M}{3}$.}
$$
\end{proposition}

\noindent
{\bf Proof.}
Suppose $0<r < \frac{\diam M}{3}$. 
Then there exists $y\in M$ such that $d(x,y)=3r/2$, for otherwise
$
B(x,3r/2)= \overline{B(x,3r/2)}\ne M
$
is simultaneously open and close, which contradicts the connectedness of $M$.
Evidently,
$B(x,r) \cap  B(y, r/2)=\emptyset$ and $B(y,r/2) \subset B(x,2r)$, which yields
$|B(x,2r)| \ge |B(y,r/2)| + |B(x,r)|$.
On the other hand
$B(x,r) \subset B(y, 5r/2)$
which along with (\ref{D1}) implies
$|B(x,r)| \le 10^d B(y,r/2)$
and hence
$| B(x,2r) | \ge (10^{-d} +1) | B(x,r) |
=  2^{\beta}  | B(x,r) |. $
$\qed$

\subsection{Useful notation and estimates}

The localization of various operator kernels in what follows will be governed by
symmetric functions of the form
\begin{equation}\label{def-Dxy}
D_{\delta, \sig}(x,y)
:= \big(|B(x, \delta)||B(y, \delta)|\big)^{-1/2}
\Big(1+ \frac{\rho(x,y)}{\delta}\Big)^{-\sig},
\quad x,y \in M.
\end{equation}
Here $\delta, \sig>0$ are parameters that will be specified in every particular case.

We next give several simple properties of $D_{\delta, \sig}(x,y)$
which will be instrumental in various proofs in the sequel.
Note first that (\ref{D1})-(\ref{D2}) readily yield
\begin{equation}\label{E1}
 D_{\delta, \sig}(x,y)
 \le 2^{d/2}|B(x, \delta)|^{-1}\Big(1+ \frac{\rho(x,y)}{\delta}\Big)^{\sig-d/2},
\end{equation}
\begin{equation}\label{E11}
D_{\lambda \delta, \sig}(x,y) \le (2/\lambda)^d D_{\delta, \sig}(x,y),
\quad 0<\lambda < 1,
\end{equation}
\begin{equation}\label{E12}
D_{\lambda \delta, \sig}(x,y) \le \lambda^\sigma D_{\delta, \sig}(x,y),
\quad  \lambda>1.
\end{equation}
Furthermore, for $0<p<\infty$ and $\sig>d(1/2+1/p)$
\begin{equation}\label{INT}
\|D_{\delta,\sigma}(x,\cdot)\|_p
=\Big(\int_M \big[D_{\delta,\sig}(x,y)\big]^p d\mu(y)\Big)^{1/p}
\le c(p) |B(x,\delta)|^{1/p-1},
\end{equation}
where $c(p)=\big(\frac{2^{dp/2}}{2^{-d}-2^{-(\sigma -d/2)p}}\big)^{1/p}$ is decreasing as a function of $p$,
and
\begin{equation}\label{Comp}
\int_M D_{\delta,\sig}(x, u) D_{\delta,\sig}(u,y) d\mu(u)
\le cD_{\delta,\sig}( x, y)\quad \hbox{if } \; \sig >  2d,
\end{equation}
with $c = \frac{2^{\sigma+d+1}}{2^{-d}- 2^{d-\sigma}}$.

The above two estimates follow readily by the following lemma
which will be needed as well.


\begin{lemma}\label{lem:tech-est}
$(a)$ If $\gg > d$, then for $\delta >0$
\begin{equation}\label{tech1}
\int_M(1+\delta^{-1}\rho(x, y))^{-\gg}d\mu(y) \le
c_1|B(x, \delta)|,
\quad x\in M,
\quad \big(c_1 = (2^{-d}- 2^{-\sigma})^{-1}\big).
\end{equation}

$(b)$ If $\gg>d$, then for $x, y\in M$ and $\delta>0$
\begin{align}\label{tech2}
\int_M\frac{1}{(1+\delta^{-1}\rho(x, u))^{\gg}(1+\delta^{-1}\rho(y, u))^{\gg}} d\mu(u)
&\le
 2^\sigma c_1 \frac{|B(x, \delta)|+|B(y, \delta)|}{(1+\delta^{-1}\rho(x, y))^{\gg}}\notag\\
&\le
2^\sigma (2^d+1) c_1\frac{|B(x, \delta)|}{(1+\delta^{-1}\rho(x, y))^{\gg-d}}.
\end{align}

$(c)$ If $\gg  > 2d $, then for $x, y\in M$ and $\delta>0$
\begin{align}\label{tech3}
\int_M\frac{1}{|B(u, \delta)|(1+\delta^{-1}\rho(x, u))^{\gg}(1+\delta^{-1}\rho(y, u))^{\gg}} d\mu(y)
\le
\frac{ c_2}{(1+\delta^{-1}\rho(x, y))^{\gg}},
\end{align}
with $c_2 = \frac{2^{\sigma+d+1}}{2^{-d}- 2^{d-\sigma}}$.

\end{lemma}

\noindent
{\bf Proof.}
Denote briefly
$E_0:=\{y\in M: \rho(x, y) < \delta\} = B(x, \delta)$ and
$$E_j:=\{y\in M: 2^{j-1}\delta \leq \rho(x, y) < 2^j\delta\}
= B(x, 2^j\delta) \setminus  B(x, 2^{j-1}\delta), \quad j\ge 1.$$
Then using (\ref{doubling}) we get
\begin{align*}
&\int_M(1+\delta^{-1}\rho(x, y))^{-\gg}d\mu(y)
= \sum_{j\ge 0} \int_{E_j}(1+\delta^{-1}\rho(x, y))^{-\gg}d\mu(y)\\
&\qquad\le |B(x,\delta)|+ (2^d-1)\sum_{j\ge 0}\frac{|B(x, 2^j\delta)|}{(1+2^{j})^\gg}\\
&\qquad\le |B(x, \delta)|\Big(1+(2^d-1)\sum_{j\ge 0}\frac{2^{jd}}{(1+2^{j})^\gg}\Big)
\le \frac{|B(x, \delta)|}{2^{-d} -2^{-\sigma}},
\end{align*}
which gives (\ref{tech1}).


For the proof of (\ref{tech2}), we note that the triangle inequality implies
$$
\frac{1+\delta^{-1}\rho(x, y)}{(1+\delta^{-1}\rho(x, u))(1+\delta^{-1}\rho(y, u))}
\le \frac{1}{1+\delta^{-1}\rho(x, u)}
+\frac{1}{1+\delta^{-1}\rho(y, u)}
$$
and hence
\begin{equation}\label{tech6}
\frac{(1+\delta^{-1}\rho(x, y))^\gg}{(1+\delta^{-1}\rho(x, u))^\gg(1+\delta^{-1}\rho(y, u))^\gg}
\le \frac{2^\gg}{(1+\delta^{-1}\rho(x, u))^\gg}
+\frac{2^\gg}{(1+\delta^{-1}\rho(y, u))^\gg}.
\end{equation}
We now integrate and use (\ref{tech1}) to obtain (\ref{tech2}).

For the proof of (\ref{tech3}), we use the above inequality and (\ref{D2})
to obtain
\begin{align}\label{tech4}
&\frac{(1+\delta^{-1}\rho(x, y))^\gg}
{|B(u, \delta)|(1+\delta^{-1}\rho(x, u))^\gg(1+\delta^{-1}\rho(y, u))^\gg}\\
&\qquad\qquad
\le \frac{2^{\gg+d}}{|B(x, \delta)|(1+\delta^{-1}\rho(x, u))^{\gg-d}}
+\frac{2^{\gg+d}}{|B(y, \delta)|(1+\delta^{-1}\rho(y, u))^{\gg-d}}\notag
\end{align}
and integrating and applying again (\ref{tech1}) we arrive at (\ref{tech3}).
$\qed$

\subsection{Maximal \boldmath $\delta$-nets}\label{sec:dd-nets}

For the construction of decomposition systems (frames) we shall need
maximal $\delta$-nets on $M$.


\begin{definition}\label{def:net}
We say that $\cX\subset M$ is a $\delta$-net on $M$ $(\delta>0)$ if $\rho(\xi, \eta) \ge \delta$
for all $\xi, \eta\in\cX$,
and $\cX\subset M$ is a maximal $\delta$-net on $M$
if $\cX$ is a $\delta$-net on $M$ that cannot be enlarged,
i.e. there does not exist $x\in M$ such that $\rho(x, \xi) \ge \delta$ for all $\xi\in\cX$
and $x\not\in \cX$.
\end{definition}

We collect some simple properties of maximal $\delta$-nets in the following proposition.


\begin{proposition}\label{prop:delta-net}
Suppose $(M, \rho, \mu)$ is a metric measure space obeying the doubling condition $(\ref{doubling})$
and let $\delta >0$.

$(a)$ A maximal $\delta$-net on $M$ always exists.

$(b)$ If $\cX$ is a maximal $\delta$-net on $M$, then
\begin{equation}\label{net-prop}
M=\cup_{\xi\in\cX} B(\xi, \delta)
\quad\hbox{and}\quad
B(\xi, \delta/2)\cap B(\eta, \delta/2)=\emptyset
\quad\hbox{if}\;\; \xi\ne\eta, \; \xi, \eta\in\cX.
\end{equation}

$(c)$ Let $\cX$ be a maximal $\delta$-net on $M$. Then $\cX$ is countable or finite
and there exists a~disjoint partition $\{A_\xi\}_{\xi\in\cX}$ of $M$ consisting
of measurable sets such that
\begin{equation}\label{B-Axi-B}
B(\xi, \delta/2) \subset A_\xi \subset B(\xi, \delta), \quad \xi\in\cX.
\end{equation}
\end{proposition}

\noindent
{\bf Proof.}
For (a)
observe that a maximal $\delta$-net is a maximal set in the collection of all $\delta$-net on $M$
with respect to the natural ordering of sets (by inclusion) and hence by Zorn's lemma
a maximal $\delta$-net on $M$ exists.

Part (b) is immediate from the definition of maximal $\delta$-nets.

To prove (c) we first fix $y\in M$ and observe that for any $n>\delta$ , $n \in \bN$,
by (\ref{D1})-(\ref{D2}) it follows that
$|B(y, n)| \le c(n, \delta) |B(\xi, \delta/2)|$
for $\xi \in \cX \cap B(y, n)$,
where $c(n, \delta)$ is a constant depending on $n$ and $\delta$.
On the other hand, by (\ref{net-prop})
$$
\sum_{\xi \in \cX \cap B(y, n)}|B(\xi, \delta/2)| \le |B(y, 2n)|\le 2^d|B(y, n)|.
$$
Therefore, $\# (\cX \cap B(y, n)) \le 2^dc(n, \delta)<\infty$,
which readily implies that $\cX$ is countable or finite.

Let us order the elements of $\cX$ in a sequence: $\cX=\{\xi_1, \xi_2, \dots\}$.
We now define the sets $A_\xi$ of the claimed cover of $M$ inductively.
We set
$$
A_{\xi_1}:= B(\xi_1, \delta)\setminus \cup_{\eta\in\cX, \eta\ne\xi_1} B(\eta, \delta/2)
$$
and if $A_{\xi_1}, A_{\xi_2}, \dots, A_{\xi_{j-1}}$ have already been defined, we set
$$
A_{\xi_{j}}:= B(\xi_{j}, \delta)\setminus
\big[\cup_{\nu\le j-1}A_{\xi_\nu}\cup_{\eta\in\cX, \eta\ne\xi_{j}} B(\eta, \delta/2)\big].
$$
It is easy to see that the sets $A_{\xi_1}, A_{\xi_2}, \dots$
have the claimed properties.
$\qed$

\smallskip

Discrete versions of estimates (\ref{Comp}) and (\ref{tech1}) will be needed.
Suppose $\cX$ is a maximal $\delta$-net on $M$ and
$\{A_\xi\}_{\xi\in\cX}$ is a~companion disjoint partition of $M$ as in Proposition~\ref{prop:delta-net}.
Then
\begin{equation}\label{discr-tech1}
\sum_{\xi\in\cX} |A_\xi|\big(1+\delta^{-1}\rho(x, \xi)\big)^{-d-1} \le 2^{2d+2}|B(x, \delta)|
\end{equation}
and
\begin{equation}\label{discr-tech11}
\sum_{\xi\in\cX} \big(1+\delta^{-1}\rho(x, \xi)\big)^{-2d-1} \le 2^{3d+2}.
\end{equation}
Furthermore, for any $\dst \ge \delta$
\begin{equation}\label{basic-est}
\sum_{\xi\in \cX} \frac{|A_\xi|}{|B(\xi, \dst)|}\big(1+\dst^{-1}\rho(x, \xi)\big)^{-2d-1} \le 2^{3d+2},
\end{equation}
and if $\sig \ge 2d+1$
\begin{equation}\label{discr-comp}
\sum_{\xi\in\cX} |A_\xi| D_{\dst,\sig}(x, \xi) D_{\dst,\sig}(y,\xi)
\le 2^{\sigma+3d+3 }D_{\dst,\sig}( x, y).
\end{equation}
Also, for $\sig \ge 2d+1$
\begin{equation}\label{discr-comp-2}
\sum_{\xi\in\cX}  \big(1+\delta^{-1}\rho(x, \xi)\big)^{-\sig}
\big(1+\delta^{-1}\rho(y, \xi)\big)^{-\sig}
\le 2^{\sigma+2d+3}\big(1+\delta^{-1}\rho(x, y)\big)^{-\sig}.
\end{equation}

We next prove (\ref{basic-est}). The proofs of (\ref{discr-tech1}) and (\ref{discr-tech11}) are similar.
Observe first that by (\ref{D2})
$|B(x, \dst)| \le 2^d (1+\dst^{-1}\rho(x, \xi))^d|B(\xi, \dst)|$.
On the other hand,
for $u\in A_\xi \subset B(\xi, \delta)$
$$
1+\dst^{-1}\rho(x, u) \le 1+\dst^{-1}\rho(x, \xi)+\dst^{-1}\rho(\xi, u)
\le 2\big(1+\dst^{-1}\rho(x, \xi)\big).
$$
Therefore,
\begin{align*}
\frac{|A_\xi|}{|B(\xi, \dst)|}\big(1+\dst^{-1}\rho(x, \xi)\big)^{-2d-1}
&\le \frac{2^d|A_\xi|}{|B(x, \dst)|}\big(1+\dst^{-1}\rho(x, \xi)\big)^{-d-1}\\
&\le \frac{2^{2d+1}}{|B(x, \dst)|}
\int_{A_\xi}\big(1+\dst^{-1}\rho(x, u)\big)^{-d-1} d\mu(u).
\end{align*}
This leads to
\begin{align*}
&\sum_{\xi\in\cX_j}\frac{|A_\xi|}{|B(\xi, \dst)|}\big(1+\dst^{-1}\rho(x, \xi)\big)^{-2d-1}\\
&\qquad\qquad\qquad\qquad\le \frac{2^{2d+1}}{|B(x, \dst)|}
\int_M\big(1+\dst^{-1}\rho(x, u)\big)^{-d-1} d\mu(u)
\le 2^{3d+2},
\end{align*}
where for the last inequality we used (\ref{tech1}).
Thus (\ref{basic-est}) is established.

For the proof of (\ref{discr-comp}), we observe that using (\ref{tech6})
\begin{align*}
&|A_\xi| D_{\dst,\sig}(x, \xi) D_{\dst,\sig}(y,\xi)
= D_{\dst,\sig}(x,y)
\frac{|A_\xi|(1+\dst^{-1}\rho(x, y))^\gg}
{|B(\xi, \delta)|(1+\dst^{-1}\rho(x, \xi))^\gg(1+\dst^{-1}\rho(y, \xi))^\gg}\\
&\qquad\quad
\le D_{\dst,\sig}(x,y)
\Big[\frac{2^{\gg}|A_\xi|}{|B(\xi, \dst)|(1+\dst^{-1}\rho(x, \xi))^{\gg}}
+\frac{2^{\gg}|A_\xi|}{|B(\xi, \dst)|(1+\dst^{-1}\rho(y, \xi))^{\gg}}
\Big].
\end{align*}
Now, summing up and applying (\ref{basic-est}) we arrive at (\ref{discr-comp}).

Estimate (\ref{discr-comp-2}) follows in a similar manner from (\ref{tech6}) and (\ref{discr-tech11}).

\subsection{Integral operators}\label{kernel-oper}

We shall mainly deal with integral (kernel) operators.

The kernels of many operators will be controlled by the quantities $D_{\delta,\sig}(x,y)$,
introduced in (\ref{def-Dxy}).
Our first order of business is to establish a Young-type inequality for such operators.


\begin{proposition}\label{prop:young}
Let $H$ be an integral operator with kernel $H(x,y)$, i.e.
$$
Hf(x)=\int_M H(x,y) f(y) d\mu(y),
\;\;\hbox{and let}\;\;\;
|H(x,y)|\le c'D_{\delta,\sig}(x,y)
$$
for some $0<\delta\le 1$ and $\sig \ge 2d+1$.
If $1\leq p \leq q \leq \infty$, then
\begin{equation}\label{young-1}
\|Hf\|_q \le  c\delta^{d(\frac 1q - \frac 1p)}\|f\|_p,
\quad f\in \LL^p,
\end{equation}
where $c=c'\hat{c}^{d(1/r-1)}2^{2d+1}$ with $\hat{c}$ being the constant from $(\ref{est-muB-unify1})$ .
\end{proposition}

This result is immediate from the following well-known lemma.


\begin{lemma}\label{lem:young}
Suppose $\frac 1p-\frac 1q =1-\frac 1r$, $1\le  p,q,r \le \infty$, and
let $H(x,y)$ be a~measurable kernel, verifying the conditions
\begin{equation}\label{young1}
\|H(\cdot,y)\|_r \le K
\quad \hbox{and} \quad
\|H(x, \cdot)\|_r \le K.
\end{equation}
If
$Hf(x) = \int_M H(x,y) f(y) d\mu(y)$, then
$$
\|Hf\|_q \le K\|f\|_p
\quad\hbox{for}\quad f\in \LL^p.
$$
\end{lemma}

For the proof, see e.g. \cite[Theorem 6.36]{F}.

\smallskip


\noindent
{\bf Proof of Proposition~\ref{prop:young}.}
Pick $1\le r\le\infty$ so that $1/p-1/q=1-1/r$.
By (\ref{INT}) and (\ref{est-muB-unify1}) we obtain
$$\|H(\cdot, y)\|_r \le c' c(r) |B(y, \delta)|^{1/r-1} \le  c' c(1) (\hat{c}\delta)^{d(1/r-1)}$$
and a similar estimate holds for $\|H(x,\cdot)\|_r$.
These estimates and the~above lemma imply (\ref{young-1}).
$\qed$

\smallskip

We shall frequently use the following well-known result (\cite{DS}, Theorem~6, p. 503).


\begin{proposition}\label{prop:kernel-oper}
An operator $T: \LL^1 \to \LL^{\infty}$ is bounded if and only if
$T$ is an integral operator with kernel $K\in L^\infty(M\times M)$, i.e.
$$
\hbox{
$Tf(x)=\int_MK(x, y)f(y)d\mu(y)$ a.e. on $M$,
}
$$
and if this is the case
$\|T\|_{1\to \infty} = \|K\|_{L^\infty}$.
Moreover, the boundedness of $T$ can be expressed in the bilinear form
$|\langle Tf, g \rangle| \le c\|f\|_{\LL^1}\|g\|_{\LL^1}$, $\forall f, g\in \LL^1$.
\end{proposition}

We next use this to derive a useful result for products of integral and non-integral operators.


\begin{proposition}\label{prop:prod-oper}
In the general setting of a doubling metric measure space $(M, \rho, \mu)$,
let $U, V: \LL^2 \to \LL^2$ be integral operators
and suppose that for some $ 0<  \delta  \leq 1$ and $\sig \ge d+1$ we have
\begin{equation}\label{local-UV}
|U(x,y)| \leq c_1D_{\delta, \sig}(x,y)
\quad\hbox{and}\quad  |V(x,y)|
\le c_2D_{\delta,\sig}(x,y).
\end{equation}
Let $R: \LL^2\to \LL^2$ be a bounded operator, not necessarily an integral operator.
Then $U R  V $ is an integral operator
with the following upper bound on its kernel
\begin{equation}\label{local-URV}
|U  R  V (x,y)|
\le \|U(x,\cdot)\|_2  \| R \|_{2 \to 2}\|V(\cdot,y)\|_2
\le  \frac{c\|R \|_{2 \to 2}}{\sqrt{|B(x, \delta)||B(y, \delta)|}}
\end{equation}
with $c:= c_1 c_22^{2d+1}$.
\end{proposition}

\noindent
{\bf Proof.}
By Proposition~\ref{prop:young} we get
\begin{align*}
\|U  R  V\|_{1\rightarrow \infty}
&\leq  \|U\|_{2\rightarrow \infty}   \|R\|_{2\rightarrow 2}  \|V\|_{1\rightarrow 2}
\le c\delta^{-d}\|R\|_{2\rightarrow 2}
\end{align*}
and, therefore, $U  R  V$ is a kernel operator.
Formally, we have
\begin{align}\label{kernel-URV}
(U R  V) f
&= \int_M U(x,u)  (R  V) f (u)d\mu(u)\notag\\
&= \int_M U(x,u) \int_M R[V(\cdot,y)](u) f(y) d\mu(y)d\mu(u)\\
&= \int_M \Big( \int_M U(x,u)R[ V(\cdot,y)](u)d\mu(u)\Big)f(y)d\mu(y)\notag
\end{align}
and hence the kernel of $U R  V$ is given by
\begin{equation}\label{split}
H(x,y) = \int_M U(x,u) R[V(\cdot,y)](u)d\mu(u)
= \big\langle U(x,\cdot), \overline{ R[V(\cdot,y)]}\big\rangle.
\end{equation}
This along with (\ref{local-UV}) and (\ref{INT}) leads to
$$
|H(x,y)| \le  \|U(x,\cdot)\|_2 \|R[V(\cdot,y)]\|_2
\le \frac{c_1c_2  [c(2)]^2\|R\|_{2\to 2}}{|B(x,\delta)|^{1/2}|B(y,\delta)|^{1/2}},
$$
which confirms (\ref{local-URV}),
taking into account that $[c(2)]^2 \leq 2^{2d+1}$ by (\ref{INT}) if $\sigma \geq d+1$.

It remains to justify the manipulations in (\ref{kernel-URV}).
Observe first that in order to prove (\ref{split}) it suffices to establish identities (\ref{kernel-URV})
for all $f\in \LL^2$ such that $\supp f\subset B(a, R)$ an arbitrary ball on $M$.
To this end we shall need Bochner's integral.
In particular, we shall use the following results (e.g. \cite{Yosida}, pp. 131--133):
Suppose $B$ is a separable Banach space
and $F: (M, \mu , \Sigma) \mapsto B$ is measurable in the following sense:
$\forall \ell \in B^*$, $x \mapsto \ell(F(x)) $ is measurable.
Then Bochner's integral $\int_M^{(B)} F(x) d\mu(x)$ is well defined and
takes its value in $B$ if and only if
$$
\int_M  \| F(x) \|_B d\mu(x) <\infty.
$$
Furthermore, if $\int_M^{(B)} F(x) d\mu(x)$ exists, then
$\ell\Big(\int_M^{(B)} F(x) d\mu(x)\Big)= \int_M \ell(F(x)) d\mu(x)$
for any $\ell\in B^*$.
Also, if
$T: B \to B$ is a bounded linear operator, then
\begin{equation}\label{operator-Bochner}
T \Big(\int^{(B)}_M F(x) d\mu(x)\Big) = \int^{(B)}_M T(F(x)) d\mu(x).
\end{equation}
We shall utilize Bochner's integral in our setting with $B=\LL^2$.

Suppose $f\in \LL^2$ and $\supp f\subset B(a, R)$, $a\in M$, $R>0$.
Then using (\ref{local-UV}), (\ref{INT}), and (\ref{D2}) we obtain
\begin{align}\label{int-norm}
\int_M \|V(\cdot, y)f(y)\|_2 d\mu(y)
&\le  c\int_{B(a,R)} |f(y)||B(y,\delta)|^{-1/2} d\mu(y)\notag\\
&\le  c\|f\|_2 \Big(\int_{B(a,R)} |B(y,\delta)|^{-1}d\mu(y)\Big)^{1/2}\\
&\le \frac{c\| f\|_2}{\sqrt{|B(a,\delta)|}}
\Big(\int_{B(a,R)} \big(1+ \delta^{-1}\rho(y,a)\big)^d d\mu(y) \Big)^{1/2} <\infty.\notag
\end{align}
Therefore,
$\int_M^{(B)}V(\cdot, y)f(y) d\mu(y)$ exists and for any $g\in\LL^2$
\begin{align*}
\big\langle  \int^{(B)}_M V(\cdot,  y)  f(y)d\mu(y) , g \big\rangle
&= \int_M \Big(\int_M  \overline{g(x)} V(x,  y)  d\mu(x) \Big)  f(y) d\mu(y)\\
&= \int_M  \overline{g(x)} \Big(\int_M  V(x, y)  f(y) d\mu(y)\Big) d\mu(x))
= \langle  Vf, g \rangle.
\end{align*}
Here the shift of the order of integration is justified by Fubini's theorem and the fact that
\begin{align*}
\int_M \int_M |V(x,y)||f(y)| |g(x)| d\mu(x) d\mu(y)
&\le \|g\|_2\Big\|\int_M |V(\cdot,y)||f(y)|d\mu(y)\Big\|_2\\
&\le \|g\|_2\int_M \|V(\cdot,y)f(y)\|_2d\mu(y)<\infty,
\end{align*}
where we used (\ref{int-norm}).
Therefore,
$
Vf=\int^{(B)}_M V(\cdot,  y)  f(y)d\mu(y).
$
We now use (\ref{operator-Bochner}) to obtain
$$
RVf = R\Big[\int^{(B)}_M V(\cdot,  y) f(y)d\mu(y)\Big]
= \int^{(B)}_M  R[V(\cdot,  y) ]f(y)d\mu(y),
$$
which implies
\begin{align*}
(URV)f(x)
&= \int_M U(x,u) (RV)f(u) d\mu(u)
= \big\langle \int^{(B)}_M   R[V(\cdot,  y) ] f(y)d\mu(y),  \overline{U(x, \cdot)} \big\rangle\\
&= \int_M  \Big(\int_M U(x, u)  R[V(\cdot,  y) ] (u) d\mu(u) \Big)  f(y)d\mu(y).
\end{align*}
Consequently,
$H(x, y)$ is given by (\ref{split})
and the proof is complete.
$\qed$

\section{Functional calculus}\setcounter{equation}{0}

The aim of this section is to develop the functional calculus of operators of the form $f(\sqrt L)$
associated with smooth and non-smooth functions $f$.
The calculus of smooth operators is in the spirit of \cite{DXO, Ouhabaz} and
will be needed in most part of this article, including the construction of frames and the Littlewood-Paley theory,
while the non-smooth calculus will be needed for estimation of the kernels of the spectral projectors and
lower bound estimates.

\subsection{Smooth functional calculus}\label{local-kernels}

We shall be operating in the setting described in \S\ref{sec:setting}.
More precisely, we assume that $(M, \rho, \mu)$ is a metric measure space obeying conditions
(\ref{doubling})-(\ref{non-collapsing}) and
$L$ is an essentially self-adjoint positive operator on $\LL^2$ such that the semi-group
$e^{-tL}$, $t>0$, has a kernel $p_t(x, y)$ verifying (\ref{Gauss-local})-(\ref{Markov}).


\begin{theorem}\label{thm:local-kernels}
Let $g: \bR \to \bC$ be a measurable function such that for some $\sigma >2d$
\begin{equation}\label{loc-ker0}
\|g\|_*:=\int_\bR |\hat g(\xi)|(1+|\xi|)^{\sigma}d\xi <\infty,
\quad\hbox{where}\quad
\hat g(\xi):=\int_\bR g(x)e^{-ix\xi} dx
\end{equation}
is the Fourier transform of $g$.
Then $g(\delta^2L) e^{- \delta^2L}$, $ 0 <\delta  \leq 1$,
is an integral operator with kernel $g(\delta^2L) e^{- \delta^2L}(x,y)$ satisfying
\begin{equation}\label{loc-ker1}
\big|g(\delta^2L) e^{- \delta^2L}(x,y)\big| \le c_\sigma \| g\|_\star D_{\delta,\sig}(x,y),
\quad\forall\, x, y\in M,
\end{equation}
and
\begin{equation}\label{loc-ker2}
\big|g(\delta^2L) e^{- \delta^2L}(x,y)  -  g(\delta^2L) e^{- \delta^2L}(x,y')\big|
\le    c_\sigma \| g\|_\star \Big(\frac{\rho(y,y')}{\delta}\Big)^\alpha D_{\delta, \sig}(x,y),
\end{equation}
for all $x,y,y'\in M$, if $\rho (y, y') \le \delta$.
Here $\alpha>0$ is the constant from $(\ref{lip})$,
$D_{\delta,\sig}(x,y)$ is defined in $(\ref{def-Dxy})$,
and $c_\sigma>0$ is a constant depending only on $\sigma$ and
the structural constants from $(\ref{hol1})-(\ref{lip})$.
Moreover,
\begin{equation}\label{loc-ker3}
\int_M  g(\delta^2L) e^{- \delta^2L}(x,y) d\mu(y) = g(0).
\end{equation}
\end{theorem}

\noindent
{\bf Proof.}
To prove (\ref{loc-ker1}) we first
show that $g(\delta^2 L) e^{-\delta^2L}$ is a kernel operator.
From (\ref{loc-ker0}) it follows that
$\|\hat{g} \|_1 <  \infty$ which implies $g(x) = \frac 1{2\pi}\int_\bR \hat{g}(\xi ) e^{ix\xi} dx$
and hence
$\|g\|_\infty\le \frac 1{2\pi}\|\hat{g} \|_1.$
Then by the spectral theorem
$$
\big\|g(\delta^2 L) e^{-\delta^2L}\big\|_{2\to 2}
= \big\|g(\delta^2 \cdot) e^{-\delta^2\cdot}\big\|_\infty
\le (2\pi)^{-1} \| \hat{g} \|_1.
$$
Therefore, invoking Proposition~\ref{prop:kernel-oper}, in order to show that
$g(\delta^2 L) e^{-\delta^2L}$ is a kernel operator it suffices to prove that
$$
\big|\langle g(\delta^2 L) e^{-\delta^2L}\ph, \psi \rangle\big| \le c\|\ph\|_1\|\psi\|_1,
\quad \forall \ph, \psi\in \LL^1\cap\LL^2.
$$
Let $E_\lambda$, $\lambda\ge 0$, be the spectral resolution associated with the operator $L$,
then $L =\int_0^\infty \lambda dE_\lambda$.
Writing the spectral decomposition of $g(\delta^2 L) e^{-\delta^2L}$
and using the Fourier inversion identity,
we obtain for $\ph, \psi \in \LL^1\cap\LL^2$
\begin{align*}
\big\langle g(\delta^2 L) e^{-\delta^2L} \ph, \psi \big\rangle
&= \int_0^\infty g(\delta^2 \lambda ) e^{-\delta^2\lambda} d\langle E_\lambda\ph, \psi \rangle\\
&= \int_0^\infty\frac 1{2\pi}\left(\int_\bR \hat{g}(\xi ) e^{i  \delta^2 \lambda\xi} d\xi\right)
e^{-\delta^2\lambda} d\langle E_\lambda \ph, \psi \rangle\\
&= \frac 1{2\pi}\int_\bR \hat{g}(\xi )
\left(\int_0^\infty  e^{- \delta^2 \lambda (1-i\xi)} d\big\langle E_\lambda \ph, \psi \big\rangle\right) d\xi \\
&= \frac 1{2\pi}\int_\bR \hat{g}(\xi )
\langle  e^{- \delta^2  (1-i\xi) L} \ph, \psi \rangle  d\xi.
\end{align*}
The above shift of the order of integration is justified by Fubini's theorem and
the fact that for any $h\in \LL^2$
\begin{align*}
\int_\bR \int_0^\infty
|\hat{g}(\xi)|
\big|e^{- \delta^2 \lambda (1-i\xi)}\big| d\|E_\lambda h\|_2^2d\xi
=\int_\bR
|\hat{g}(\xi)|d\xi
\int_0^\infty
e^{- \delta^2 \lambda} d\|E_\lambda h\|_2^2
\le \|\hat g\|_1\|h\|_2^2.
\end{align*}
To go further, we use that
$e^{- \delta^2(1-i\xi)L}$ is an integral operator with kernel
$p_z(x,y)$, $z= \delta^2(1-i\xi)$, and $\| p_z \|_\infty \le c$
to obtain for
$\ph, \psi \in \LL^1\cap\LL^2$
\begin{align}\label{kernel5}
\big\langle g(\delta^2 L) e^{-\delta^2L}\ph, \bar{\psi}\big\rangle
&= \frac 1{2\pi}   \int_\bR \hat{g}(\xi)
\Big(\int_M \int_M p_{ \delta^2  (1-i\xi) }(x,y)\phi(x) \psi(y) d\mu(x) d\mu(y) \Big)  d\xi\notag\\
&= \int_M \int_M
\Big[\frac 1{2\pi}\int_\bR \hat{g}(\xi )p_{ \delta^2  (1-i\xi) }(x,y) d\xi\Big]
\phi(x) \psi(y) d\mu(x) d\mu(y).
\end{align}
To justify the above shift of order of integration  we again use Fubini's theorem
and the fact that
\begin{align*}
\int_\bR
\int_M \int_M 
|\hat{g}(\xi ) || p_{ \delta^2  (1-i\xi) }(x,y) | | \phi(x) ||\psi(y)| d\mu(x) d\mu(y)d\xi
\le c\|\hat g\|_1\|\ph\|_1\|\psi\|_1 <\infty.
\end{align*}
This also implies
$
|\big\langle g(\delta^2 L) e^{-\delta^2L}\ph, \bar{\psi}\big\rangle|
\le c\|\hat g\|_1\|\ph\|_1\|\psi\|_1
$
for all
$\ph, \psi \in \LL^1\cap\LL^2$.
Therefore, $g(\delta^2L) e^{- \delta^2L}$ is a kernel operator and by (\ref{kernel5})
\begin{equation}\label{rep-kernel-1}
g(\delta^2L) e^{- \delta^2L}(x, y)
= \frac 1{2\pi}\int_\bR \hat{g}(u) p_{\delta^2(1-iu)}(x, y)du.
\end{equation}
From this and (\ref{hol1}) we infer
\begin{equation}\label{rep-kernel-2}
|g(\delta^2L) e^{- \delta^2L}(x,y)| \le  c'\big(|B(x, \delta)||B(y, \delta)|\big)^{-1/2}
\int_\bR|\hat{g}(u)|\exp\Big\{-\frac{c\rho^2(x,y)}{\delta^2(1+u^2)}\Big\}du.
\end{equation}
Assume $\rho(x, y)/\delta\ge 1$.
Clearly,
$\sup_{x\geq 0} x^\beta e^{-x}=( \frac \beta e)^\beta$ for $\beta >0$.
Using this with $\beta=\sig/2$ we obtain
\begin{align*}
&\exp\Big\{-\frac{c\rho^2(x,y)}{\delta^2(1+u^2)}\Big\}
\le \exp\Big\{-\Big(1+\frac{\rho^2(x,y)}{\delta^2}\Big)\frac{c}{2(1+u^2)}\Big\}\\
&\qquad\qquad \le c'\Big(1+\frac{\rho^2(x,y)}{\delta^2}\Big)^{-\sigma/2}(1+u^2)^{\sigma/2}
\le c''\Big(1+\frac{\rho(x,y)}{\delta}\Big)^{-\sigma}(1+|u|)^{\sigma}.
\end{align*}
Therefore,
\begin{align*}
|g(\delta^2L) e^{- \delta^2L}(x,y)|
&\le \frac{c(1+ \frac{\rho(x,y)}{\delta} )^{-\sig}}{(|B(x, \delta) ||B(y,\delta|)^{1/2}}
\int_\bR |\hat{g}(u)|(1+|u|)^{\sig}du\\
& = c \int_\bR |\hat{g}(u)|(1+|u|)^{\sig}du \; D_{\sigma,\delta}(x,y),
\end{align*}
which confirms (\ref{loc-ker1}).

If $\rho(x, y)/\delta < 1$, then by (\ref{rep-kernel-2})
\begin{align*}
|g(\delta^2L) e^{- \delta^2L}(x,y)|
&\le  c'\big(|B(x, \delta)||B(y, \delta)|\big)^{-1/2}\int_\bR|\hat{g}(u)|du\\
& \le c \int_\bR |\hat{g}(u)|(1+|u|)^{\sig}du \; D_{\sigma,\delta}(x,y).
\end{align*}
This completes the proof of (\ref{loc-ker1}).


We now take on (\ref{loc-ker2}).
As
$g(\delta^2L) e^{- \delta^2L} = g(\delta^2L)e^{- \frac12\delta^2L}e^{- \frac12\delta^2L}$,
the kernels of these operators are related by
$$
g(\delta^2L) e^{- \delta^2L}(x, y)
= \int_M g(\delta^2L)e^{- \frac12\delta^2L}(x, u)e^{- \frac12\delta^2L}(u, y) d\mu(u),
$$
which implies
\begin{align*}
&\big|g(\delta^2L) e^{- \delta^2L}(x,y)  -  g(\delta^2L) e^{- \delta^2L}(x,y')\big|\\
&\qquad\qquad\qquad \le \int_M \big|g(\delta^2L)e^{- \frac12\delta^2L}(x, u)\big|
\big|p_{ \delta^2/2}(u, y)- p_{\delta^2/2}(u, y')\big| d\mu(u).
\end{align*}
We use (\ref{loc-ker1}) with $\delta$ replaced by $\delta/\sqrt{2}$
and $g(\lambda)$ by $g(2\lambda)$
to estimate the first term under the integral and
(\ref{lip}) for the second term,
taking into account that
$
\exp\big\{-\frac{c\rho^2(x, y)}{\delta^2}\big\}
\le c_\sig\big(1+\frac{\rho(x, y)}{\delta}\big)^{-\sig}.
$
Thus we get
\begin{align*}
&\big|g(\delta^2L) e^{- \delta^2L}(x,y)  -  g(\delta^2L) e^{- \delta^2L}(x,y')\big|\\
&\quad
\le c\|g\|_\star\Big(\frac{\rho(y,y')}{\delta}\Big)^\alpha
\int_M D_{\delta, \sig}(x,u)D_{\delta, \sig}(u,y)d\mu(u)
\le c\|g\|_\star\Big(\frac{\rho(y,y')}{\delta}\Big)^\alpha D_{\delta, \sig}(x,y).
\end{align*}
Here for the latter estimate we used (\ref{Comp}) and that $\sig>2d$.


It remains to prove (\ref{loc-ker3}). By (\ref{Markov}), i.e.
$\int_M p_{\delta^2-iu} (x,y) dy  \equiv 1$,
and (\ref{rep-kernel-1}) we get
\begin{align*}
\int_M  g(\delta^2L) e^{- \delta^2L}(x,y) dy
&= \frac 1{2\pi}\int_\bR \hat{g}(u)\int_M p_{\delta^2-iu} (x,y)d\mu(y) du\\
&= \frac 1{2\pi}\int_\bR \hat{g}(u)du
=g(0).
\end{align*}
Here the justification of the shift of order of integration is by straightforward
application of Fubini's theorem.
$\qed$

\smallskip

Some remarks are in order.
Condition (\ref{loc-ker0}) is apparently a smoothness condition on $g$.
By Cauchy-Schwartz it follows that
$$
\int_\bR |\hat g(\xi)|(1+|\xi|^2)^{\sigma/2}d\xi
\le c \Big(\int_\bR |\hat g(\xi)|^2(1+|\xi|^2)^{\sigma+1}d\xi\Big)^{1/2}
=c\|g\|_{H^{\sig+1}}
$$
and hence (\ref{loc-ker0}) holds if $\|g\|_{H^{\sig+1}}<\infty$.
However, it will be more convenient to us to replace (\ref{loc-ker0}) by
a condition in terms of derivatives of $g$ that is easier to verify.
From
$
\xi^k \hat g(\xi) = (-i)^k\widehat{g^{(k)}}(\xi)
$
we get
$|\xi|^k|\hat g(\xi)| \le \|g^{(k)}\|_{L^1}$.
Also, $|\hat g(\xi)| \le \|g\|_{L^1}$.
Pick $k\ge \sig >2d$. Then using the above we obtain
$$
(1+|\xi|)^{k+2}|\hat g(\xi)|
\le 2^{k+1}\big(|\hat g(\xi)|+|\xi|^{k+2}|\hat g(\xi)|\big)
\le 2^{k+1}\big(\|g\|_{L^1}+\|g^{(k+2)}\|_{L^1}\big)
$$
that implies
\begin{align*}
\|g\|_*:=\int_\bR |\hat g(\xi)|(1+|\xi|)^{k}d\xi
& = \int_\bR |\hat g(\xi)|(1+|\xi|)^{k+2}(1+|\xi|)^{-2}d\xi\\
& \le c\big(\|g\|_{L^1}+\|g^{(k+2)}\|_{L^1}\big).
\end{align*}
Thus we arrive at the following


\begin{remark}\label{rem:smooth-cond}
For the norm $\|g\|_*$ from condition $(\ref{loc-ker0})$ we have
$\|g\|_* \le c\|g\|_{H^{\sig+1}}$
and
$\|g\|_* \le c\big(\|g\|_{L^1}+\|g^{(k+2)}\|_{L^1}\big)$
if $k\ge \sig >2d$.
\end{remark}


\begin{corollary}\label{cor:local-Lm-exp}
For any $m\in \bN$ and $\sig>0$ there exists a constant $c_{\sig,m}>0$ such that
the kernel of the operator
$L^m e^{-\delta^2L}$, $ 0 <\delta  \leq 1$, satisfies
\begin{equation}\label{Lm-exp-1}
\big|L^m e^{-\delta^2L}(x, y)\big| \le c_{\sig,m} \delta^{-2m}D_{\delta, \sig}(x,y)
\quad \hbox{and}
\end{equation}
\begin{equation}\label{Lm-exp-2}
\big|L^m e^{-\delta^2L}(x, y)  -  L^m e^{-\delta^2L}(x,y')\big|
\le c_{\sig,m}\delta^{-2m}\Big(\frac{\rho(y,y')}{\delta}\Big)^\alpha D_{\delta, \sigma}(x,y),
\end{equation}
if $\rho (y, y') \le \delta$.
\end{corollary}

\noindent
{\bf Proof.}
Set $g(\lam) := \lam^m \theta(\lam)e^{-\lam} $ for $\lambda \geq 0$,
where $\theta\in C^\infty(\R)$, $\supp \theta \subset [-1, \infty)$, and
$\theta(\lam)=1$ for $\lambda \ge 0$.
Since $L \geq 0$,  we can write
$$
L^m e^{-\delta^2 L}
= 2^{m}\delta^{-2m}g\big(\delta_*^2L\big)
e^{-\delta_*^2 L}
\quad\hbox{with}\quad \delta_*:=2^{-1/2}\delta
$$
and the corollary follows by Theorem~\ref{thm:local-kernels} and (\ref{E11}).
$\qed$

\smallskip

We next use Theorem~\ref{thm:local-kernels} and Remark~\ref{rem:smooth-cond} to obtain some
important kernel localization results.
Our main interest is in operators of the form $f(\delta \sqrt L)$.


\begin{theorem}\label{thm:main-local-kernels}
Let $f\in C^{2k+4}(\bR_+)$, $k > 2d$, $\supp f\subset [0, R]$ for some $R\ge 1$, and
$f^{(2\nu+1)}(0)=0$ for $\nu=0, \dots, k+1$.
Then $f(\delta \sqrt L)$, $0<\delta\le 1$, is an integral operator with kernel
$f(\delta \sqrt L)(x, y)$ satisfying
\begin{equation}\label{main-loc-ker1}
\big|f(\delta \sqrt L)(x, y)\big| \le c_k D_{\delta, k}(x,y)
\quad \hbox{and}
\end{equation}
\begin{equation}\label{main-loc-ker2}
\big|f(\delta \sqrt L)(x, y)  -  f(\delta \sqrt L)(x,y')\big|
\le c_k'\Big(\frac{\rho(y,y')}{\delta}\Big)^\alpha D_{\delta, k}(x,y)
\;\; \hbox{if}\;\; \rho (y, y') \le \delta,
\end{equation}
where
$
c_k = c_k(f)= \tilde c_kR^{2k+d+4} \big(\|f\|_{L^\infty} + \|f^{(2k+4)}\|_{L^\infty}
+ \max_{\nu \le 2k+4 } |f^{(\nu)}(0)|\big)
$
with
$\tilde c_k>0$ a~constant depending only on $k, d$,
and the constants in $(\ref{hol1})-(\ref{lip})$, and
$c_k'=c_kR^\alpha$;
as before $\alpha>0$ is the constant from $(\ref{lip})$.
Furthermore,
\begin{equation}\label{main-loc-ker3}
\int_M  f(\delta \sqrt L)(x, y) d\mu(y) = f(0).
\end{equation}
\end{theorem}

\noindent
{\bf Proof.} We first observe that it suffices to only prove the theorem when $R=1$,
then in the general case it follows by rescaling.
Indeed, assume that $f$ satisfies the hypotheses of the theorem and
set $h(\lambda):= f(R\lambda)$, $\lambda\in\R_+$.
Then $h$ verifies the assumptions with $R=1$ and
if the theorem holds for $R=1$ we obtain, using (\ref{E11}),
\begin{equation}\label{loc-ker4}
| f(\delta \sqrt L)(x,y)| =|h(\delta R^{-1} \sqrt L)(x,y)|
\le  c_k(h) D_{\delta/R, k}(x,y)\le (2R)^d c_k(h) D_{\delta, k}(x,y)
\end{equation}
and similarly
$$
\big|f(\delta \sqrt L)(x, y)  -  f(\delta \sqrt L)(x,y')\big|
\le (2R)^{d+\alpha}c_k'(h)\Big(\frac{\rho(y,y')}{\delta}\Big)^\alpha D_{\delta, k}(x,y)
\;\; \hbox{if\; $\rho (y, y') \le \frac{\delta}{R}$}.
$$
For $\frac{\delta}{R} <\rho (y, y') \le \delta$, the last estimate follows by (\ref{loc-ker4}).
It remains to observe that
\begin{align*}
c_k(h) 
&= \tilde c_k\big(\|f\|_{L^\infty} + R^{2k+4} \|f^{(2k+4)}\|_{L^\infty}  +\max_{\nu \le 2k+4 } R^\nu|f^{(\nu)}(0)|\big)\\
&\le \tilde c R^{2k+4}\big(\|f\|_{L^\infty} +  \|f^{(2k+4)}\|_{L^\infty}  +\max_{\nu \le 2k+4 } |f^{(\nu)}(0)| \big)
\end{align*}
and hence the theorem holds in general.


We now prove the theorem in the case when $R=1$.
Choose $\theta\in C^{\infty}(\bR)$ so that
$\theta(\sqrt{\cdot})\in C^{\infty}(\bR)$,
$\theta$ is even,
$\supp \theta\subset [-1, 1]$,
$\theta(\lambda)=1$ for $\lam \in [-1/2, 1/2]$,
and $0\le \theta \le 1$.
Denote
$P_k(\lambda):=\sum_{j=0}^{k+2} \frac{f^{(2j)}(0)}{(2j)!} \lambda^{2j}$
and let $f_1(\lambda)$, $g_0(\lambda)$, and $g_1(\lambda)$ be defined for $\lambda\in\R_+$ from
$$
f(\lambda)= \theta(\lambda)P_k(\lambda) + f_1(\lambda),
\quad
\theta(\lambda)P_k(\lambda) = g_0(\lambda^2) e^{-\lambda^2},
\quad
f_1(\lambda)= g_1(\lambda^2) e^{-\lambda^2}.
$$
Thus
$g_0(\lambda)= P_k(\sqrt{|\lambda|})\theta(\sqrt{|\lambda|}) e^\lambda$ for $\lambda\in\R_+$,
and we use this to define $g_0(\lambda)$ for $\lambda<0$.
Clearly,
$g_0\in C^{\infty}(\bR)$, $\supp g_0\subset [-1, 1]$ and
$$
\|g_0\|_{L^1}+\|g_0^{(k+2)}\|_{L^1} \le  c(k)\sup_{\nu \le 2k+4} |f^{(\nu)}(0)|.
$$
Therefore, by Theorem~\ref{thm:local-kernels} the kernel of the operator
$\theta(\delta\sqrt L)P_k(\delta \sqrt L)$ satisfies the desired inequalities
(\ref{main-loc-ker1})-(\ref{main-loc-ker2}) with $R=1$.

On the other hand,
$g_1(\lamm)= f_1\big(\sqrt{|\lamm|}\big)e^\lamm$ for $\lamm\in \bR_+$
and we use this to define $g_1(\lamm)$ for $\lamm<0$.
Observe that
$f_1(\delta \sqrt L)=g_1(\delta^2L)e^{-\delta^2L}$
and $\supp g_1\subset [-1, 1]$.
Furthermore,
$f_1\in C^{2k+4}(\R_+)$,
$f_1^{(\nu)}(0)=0$, $\nu=0, \dots, 2k+4$, and
\begin{equation}\label{der-estim}
\|f_1^{(j)}\|_{L^\infty} \le \|f^{(j)}\|_{L^\infty} + c\max_{\nu\le 2k+4} |f^{(\nu)}(0)|,
\quad 0\le j\le 2k+4.
\end{equation}
We next show that
$g_1\in C^{k+2}(\bR)$
and estimate the derivatives of $g_1$.
We have for $1\le m\le k+2$ and $\lamm>0$
$$
g_1^{(m)}(\lamm)=\sum_{\nu=0}^m \binom{m}{\nu}e^\lamm
\Big(\frac{d}{d\lamm}\Big)^\nu\big[f_1(\sqrt \lamm)\big]
$$
and a little calculus shows that for $\nu\ge 1$ and $\lamm>0$
$$
\Big(\frac{d}{d\lamm}\Big)^\nu\big[f_1(\sqrt \lamm)\big]
=\sum_{j=1}^\nu c_j\lamm^{-\nu+j/2}f_1^{(j)}(\sqrt \lamm),
\quad \hbox{where}\quad |c_j|\le \nu!.
$$
On the other hand, by Taylor's theorem
$|f_1^{(j)}(\sqrt \lam)| \le |\lam|^{(2m-j)/2}\|f_1^{(2m)}\|_{L^\infty}$
and hence
$$
\Big|\Big(\frac{d}{d\lamm}\Big)^\nu\big[f_1(\sqrt{|\lamm|})\big]\Big|
\le c|\lamm|^{m-\nu}\|f_1^{(2m)}\|_{\infty},
\quad 1\le \nu\le m.
$$
Exactly in the same way we obtain the same estimate for $\lamm<0$.
Denote briefly $h(\lamm):=f_1(\sqrt{|\lamm|})$.
Observe that since $f_1\in C^{2k+4}(\bR_+)$ we have $h^{(k+2)}(\lamm) =o(1)$ as $\lamm \to 0$.
This and the above inequalities yield $h^{(\nu)}(0)=0$, $\nu=0, \dots, k+2$,
and hence $h\in C^{k+2}(\bR)$, which implies $g_1\in C^{k+2}(\bR)$.
From the above we also obtain
$$
|g_1^{(m)}(\lamm)|
\le c\sum_{\nu=0}^m  e^\lamm|\lamm|^{m-\nu}\|f_1^{(2m)}\|_{L^\infty}
\le c(m+1)\|f_1^{(2m)}\|_{L^\infty},
\quad \lamm\in \bR.
$$
This in turn (with $m=k+2$) implies
$\|g_1^{(k+2)}\|_{L^1} \le c(k+3)\|f_1^{(2k+4)}\|_{L^\infty}$
and, evidently, 
$\|g_1\|_{L^1} \le e\|f_1\|_{L^\infty}$.
We now apply Theorem~\ref{thm:local-kernels} to conclude that
$f_1(\delta \sqrt L)$ is an integral operator with kernel
$f_1(\delta \sqrt L)(x, y)$ satisfying (\ref{main-loc-ker1})-(\ref{main-loc-ker2}),
where, in view of Remark~\ref{rem:smooth-cond} and (\ref{der-estim}),
the constants $c_k$, $c_k'$ are of the claimed form.

Putting the above together we conclude that
$f(\delta \sqrt L)$ is an integral operator with kernel
$f(\delta \sqrt L)(x, y)$ satisfying (\ref{main-loc-ker1})-(\ref{main-loc-ker2})
with $R=1$.

Identity (\ref{main-loc-ker3}) follows by (\ref{loc-ker3}).
$\qed$


\begin{corollary}\label{cor:LS-local-kernels}
Let $f:\R_+\to \bC$ be as in the hypothesis of Theorem~\ref{thm:main-local-kernels}.
Then for any $m\in \bN$ and $0< \delta\le 1$ the operator
$L^m f(\delta \sqrt L)$ is an integral operator with kernel $L^m f(\delta \sqrt L)(x, y)$
such that
\begin{equation}\label{L-loc-ker1}
\big|L^mf(\delta \sqrt L)(x, y)\big| \le c_{k, m}\delta^{-2m} D_{\delta, k}(x,y)
\quad \hbox{and}
\end{equation}
\begin{equation}\label{L-loc-ker2}
\big|L^m f(\delta \sqrt L)(x, y)  -  L^m f(\delta \sqrt L)(x,y')\big|
\le c_{k, m}'\delta^{-2m}\Big(\frac{\rho(y,y')}{\delta}\Big)^\alpha D_{\delta, k}(x,y)
\end{equation}
whenever $\rho (y, y') \le \delta$.
Here the constants $c_{k, m}$, $c_{k, m}'$ are as the constants $c_{k}$, $c_{k}'$
in Theorem~\ref{thm:main-local-kernels} with $R^{2k+d+4}$ replaced by $R^{2k+d+4+2m}$ and
$\tilde c_k$ depending on $m$ as well.
\end{corollary}

\noindent
{\bf Proof.}
Let $h(\lambda):= \lambda^{2m} f(\lambda)$.
Then $h(\delta\sqrt L)= \delta^{2m}L^m f(\delta\sqrt L)$
and observe that
$h^{(2\nu+1)}(0)=0$ for $\nu=0, \dots, k+1$.
Consequently, the corollary follows by Theorem~\ref{thm:main-local-kernels} applied to $h$.
$\qed$


\begin{corollary}\label{lip0}
Let $f:\R_+\to \bC$ be as in the hypothesis of Theorem~\ref{thm:main-local-kernels}.
Then there exists a constant $c>0$ such that for any $0<\delta \leq 1$
$$
\| f(\delta \sqrt L)\phi \|_q \le c\delta^{1/p-1/q}\| \phi \|_p,
\quad \forall \phi \in \bL^p,
\quad 1\leq p\leq q \leq \infty,
$$
and
$$
|f(\delta \sqrt L)\phi(x) -  f(\delta \sqrt L)\phi(y)| \le c\|\phi \|_\infty
\Big(\frac{\rho(x,y)}{\delta}\Big)^\alpha, \quad x,y \in M,
\quad\forall \phi \in \bL^\infty.
$$
\end{corollary}
This corollary is an immediate consequence of Theorem~\ref{thm:main-local-kernels}
and Proposition \ref{prop:young}.

\subsection{Non-smooth functional calculus.}

We need to establish some properties of operators of the form $f(\sqrt L)$ and their kernels
in the case of non-smooth compactly supported functions $f$. These are kernel operators with
not necessarily well localized kernels.


\begin{theorem}\label{thm:rough-kernels}
Let $f$ be a bounded measurable function on $\bR_+$ with
$\supp f \subset [0,\tt]$ for some $\tt \ge 1$.
Then $f(\sqrt L)$ is an integral operator with kernel $f(\sqrt L)(x, y)$ satisfying
\begin{equation}\label{rough1}
| f(\sqrt L) (x,y) |
\le  \frac{c\| f \|_\infty}
{\sqrt{|B(x, \tt^{-1})|| B(y, \tt^{-1})|}},\quad x,y\in M,
\end{equation}
and for $x,y,y'\in M$
\begin{equation}\label{rough2}
|f(\sqrt L)(x,y)- f(\sqrt L)(x,y')|
\le  \frac{c[\tt\rho(y,y')]^\alpha \| f \|_\infty}
{\sqrt{|B(x, \tt^{-1})|| B(y, \tt^{-1})|}}
\quad \hbox{if} \;\;\;  \rho(y, y')\le \tt^{-1}.
\end{equation}
Furthermore, if $1\leq p \leq 2 \leq q \leq \infty$,
\begin{equation}\label{rough3}
\| f(\sqrt L)  \|_{p \rightarrow q} \le  c\tt^{d(1/p - 1/q)}\|f\|_\infty,
\end{equation}

\begin{equation}\label{rough4}
\|f(\sqrt L)(. ,x)\|_2^2 = |f|^2(\sqrt L)(x,x)
\le c|B(x,\tt^{-1})|^{-1}\|f\|_\infty^2,
\quad \hbox{and}
\end{equation}

\begin{equation}\label{rough5}
\||f|^2(\sqrt L)\|_{1 \rightarrow \infty}
= \sup_{x\in M}|f|^2(\sqrt L)(x,x).
\end{equation}
Above the constants depend only on $d$ and the constants in $(\ref{hol1})$ and $(\ref{lip})$;
the~constant in $(\ref{rough3})$ depends in addition on $p, q$.
\end{theorem}

\noindent
{\bf Proof.}
Pick a function $\theta  \in C^\infty(\bR_+)$ so that
$\supp \theta \subset [0, 2]$,
$\theta(x)=1$ for $x \in [0, 1]$, and $0\le \theta\le 1$.
Then by Theorem~\ref{thm:main-local-kernels}
\begin{equation}\label{rough-6}
|\theta (\tt^{-1} \sqrt L)(x,y) | \le c_\sig D_{\tt^{-1},\sig}(x,y)
\quad\hbox{for any} \; \sig>0.
\end{equation}
Choose $\sig>3d/2$.
We have
\begin{align}\label{rough-7}
f(\sqrt L)
&= \int_0^\infty f(\sqrt \lambda) dE_\lambda
=\int_0^\infty \theta (\tt^{-1}\sqrt\lambda)f(\sqrt\lambda)\theta(\tt^{-1}\sqrt\lambda)dE_\lambda\notag\\
&=\theta(\tt^{-1} \sqrt L)   f(\sqrt L)  \theta (\tt^{-1} \sqrt L).
\end{align}
Now, (\ref{rough1}) follows by Proposition~\ref{prop:prod-oper},
using the above, (\ref{rough-6}),
and the fact that $\|f(\sqrt L)\|_{2\rightarrow 2} \leq \|f\|_\infty$.

From (\ref{rough-6})-(\ref{rough-7}) and Proposition~\ref{prop:prod-oper} we also obtain
for $1\leq p \leq 2 \leq q \leq \infty$
\begin{align*}
\|f(\sqrt L)\|_{p\rightarrow q}
&\leq  \| \theta (\tt^{-1} \sqrt L)\|_{p\rightarrow 2} \| f(\sqrt L)\|_{2\rightarrow 2}
\| \theta (\tt^{-1} \sqrt L)  \|_{2\rightarrow q}\\
&\le c\| f \|_\infty\tt^{-d(1/q - 1/p)},
\end{align*}
which confirms (\ref{rough3}).

For the proof of (\ref{rough2}), we first observe that
\begin{align*}
f(\sqrt L)
&= \int_0^{\infty} f( \sqrt \lambda) e^{\tt^{-2}(\sqrt \lambda )^2}
e^{-\tt^{-2}\lambda} dE_\lambda
= \int_0^{\infty} g(\sqrt\lambda)e^{- \tt^{-2}\lambda} dE_\lambda
= g(\sqrt L) e^{-\tt^{-2}L},
\end{align*}
where $g(u) := f(u) e^{\tt^{-2}u^2}$, $\|g\|_\infty \le e \|f\|_\infty$,
and hence
$$
f(\sqrt L) (x,y)- f(\sqrt L) (x,y')
=\int_M  g(\sqrt L)(x,u) \big[e^{-\tt^{-2}L}(u,y) -e^{-\tt^{-2}L}(u,y')\big]d\mu(u).
$$
We now use (\ref{rough1}), applied to $g(\sqrt L)$, and (\ref{lip}) to obtain
\begin{align*}
&|f(\sqrt L)(x,y)-f(\sqrt L)(x,y')|\\
&\quad \leq c(\tt\rho(y,y'))^\alpha \|g\|_\infty
\int_M \frac 1{\sqrt{|B(x,\tt^{-1})||B(u,\tt^{-1})|}}
\frac{e^{-(\tt\rho(u,y))^2}}{\sqrt{|B(u, \tt^{-1})||B(y, \tt^{-1})|}}d\mu(u)\\
& \quad \le \frac{c(\tt\rho(y,y'))^\alpha \| f \|_\infty}{\sqrt{|B(x,\tt^{-1})||B(y,\tt^{-1})|}}
\int_M \frac{e^{-(\tt\rho(u,y))^2}}{| B(u,\tt^{-1})|}d\mu(u).
\end{align*}
However, using (\ref{D2}) we have
$$
\int_M \frac{e^{- (\tt \rho(u,y))^2}}{|B(u, \tt^{-1}) | } d\mu(u)
\le \frac{2^d}{|B(y, \tt^{-1})|}
\int_M(1+ \tt\rho(u,y))^d e^{- (\tt \rho(u,y))^2} d\mu(u) \le c <\infty,
$$
where for the latter inequality we used (\ref{tech1}).
This completes the proof of (\ref{rough2}).

We now turn to the proof of (\ref{rough4}). We have
\begin{align*}
\| f(\sqrt L) (. ,y)\|_2^2
&= \int_M |f(\sqrt L) (x,y) |^2 dy
= \int_M f(\sqrt L)(x,y) \overline{f(\sqrt L) (x,y)} d\mu(y)\\
& = \int_M f(\sqrt L)(x,y)\overline{f}(\sqrt L) (y,x) d\mu(y)
=|f|^2(\sqrt L)(x,x)\\
& \le c|B(x,\tt^{-1})|^{-1}\|f\|_\infty^2,
\end{align*}
which proves (\ref{rough4}). Here for the latter estimate we used (\ref{rough1}).

Finally, using the above we have
\begin{align*}
|f|^2(\sqrt L)(x,y)
&= \int_M f(\sqrt L) (x,u)  \overline{ f(\sqrt L) (y,u)} d\mu(u)\\
&\le \Big(\int_M|f(\sqrt L) (x,u) |^2d\mu(u)\Big)^{1/2}
\Big(\int_M | f(\sqrt L) (y,u) |^2d\mu(u)\Big)^{1/2}\\
& = \big(|f |^2(\sqrt L)(x,x)\big)^{1/2}\big(|f |^2(\sqrt L)(y,y)\big)^{1/2}
\end{align*}
and hence
$
\||f|^2(\sqrt L) \|_{1 \rightarrow \infty}
= \sup_{x,y } | |f |^2(\sqrt L)(x,y)| = \sup_{x} |f |^2(\sqrt L)(x,x),
$
which confirms (\ref{rough5}).
$\qed$

\subsection{Approximation of the identity and Littlewood-Paley decomposition}\label{app-identity}

We first give a convenient approximation of the identity in $\LL^p$ statement.


\begin{proposition}\label{prop:app-identity}
Let $\varphi \in C^\infty(\bR_+)$, $\supp \varphi \subset [0, R]$, $R>0$,
$\varphi(0)=1$, and $\varphi^{(2\nu+1)}(0)=0$ for $\nu=0, 1, \dots$.
Then for any $f\in \LL^p$, $1\le p\le \infty$, $(L^\infty:=\UCB)$ one has
$$
f= \lim_{\delta\to 0}\varphi(\delta\sqrt L) f
\quad\hbox{in} \;\; \LL^p.
$$
\end{proposition}

\noindent
{\bf Proof.}
By Theorem~\ref{thm:main-local-kernels} it follows that $\varphi(\delta\sqrt L)$
is an integral operator
with kernel $\varphi(\delta\sqrt L)(x, y)$ satisfying for any $k > 2d$
\begin{equation}\label{identity1}
|\varphi(\delta\sqrt L)(x, y)| \le c_k D_{\delta, k}(x, y)
\le c|B(x, \delta)|^{-1}\big(1+\delta^{-1}\rho(x, y)\big)^{-k+d/2},
\end{equation}
where for the last inequality we used (\ref{D2}).
Now, just as in the proof of (\ref{tech1}) we obtain for $k>3d/2$ and $r>0$
$$
\int_{M\setminus B(x, r)} |\varphi(\delta\sqrt L)(x, y)|d\mu(y)
\le c(\delta/r)^{k-3d/2} \to 0
\quad\hbox{as}\;\; \delta \to 0.
$$
Indeed, suppose $2^{\ell-1}\delta \le r < 2^\ell\delta$ and denote
$E_j:= B(x, 2^j\delta) \setminus  B(x, 2^{j-1}\delta)$.
Then using (\ref{identity1}) and (\ref{D1}) we get
\begin{align*}
&\int_{M\setminus B(x, r)} |\varphi(\delta\sqrt L)(x, y)|d\mu(y)
\le c|B(x, \delta)|^{-1}\sum_{j\ge \ell} \int_{E_j}(1+\delta^{-1}\rho(x, y))^{-k+d/2}d\mu(y)\\
&\qquad\qquad\qquad\le c|B(x,\delta)|^{-1}\sum_{j\ge \ell}\frac{|B(x, 2^j\delta)|}{(1+2^{j})^{k-d/2}}
%
\le c2^{-\ell(k-3d/2)} \le c(\delta/r)^{k-3d/2}.
\end{align*}
On the other hand, from (\ref{main-loc-ker3}) and $\varphi(0)=1$ we have
$\int_M \varphi(\delta\sqrt L)(x, y) d\mu(y) =1$.
Using the above and the fact that the vector lattice set of all boundedly supported
uniformly continuous functions on $M$ is dense in $\LL^p$ (by the Stone-Daniell theorem)
one proves as usual the claimed convergence.
$\qed$

\smallskip

We next give precise meaning to what we call {\em Littlewood-Paley decomposition} of
$\LL^p$-functions in this article.


\begin{corollary}\label{cor:Littlewood-Paley}
Let $\varphi_0, \varphi \in C^\infty(\R_+)$,
$\supp \varphi_0 \subset [0, b]$ and $\supp \varphi \subset [b^{-1}, b]$ for some $b>1$,
$\varphi(0)=1$, $\varphi^{(2\nu+1)}(0)=0$ for $\nu \ge 0$, and
$
\varphi_0(\lambda)+\sum_{j\ge 1}\varphi(b^{-j}\lam) =1
$
for $\lambda\in \R_+$.
Then for any $f\in \LL^p$, $1\le p \le\infty$, $(\LL^\infty:=\UCB)$
\begin{equation}\label{Littlewood-1}
f=\varphi_0(\sqrt L)+\sum_{j\ge 1}\varphi(b^{-j}\sqrt L)f
\;\;\mbox{ in }\; \LL^p.
\end{equation}
\end{corollary}

\noindent
{\bf Proof.}
Let $\theta(\lam):= \varphi_0(\lam) + \varphi(b^{-1}\lam)$
and observe that
$\sum_{k=0}^j \varphi_k(\lam)= \theta(b^{-j}\lam)$ for $j\ge 1$.
Then the result follows by Proposition~\ref{prop:app-identity}.
$\qed$

\subsection{Spectral spaces}\label{spectral-spaces}

We adhere to the setting of this article, described in the introduction.
As before $E_\lambda$, $\lambda\ge 0$, is the spectral resolution associated with
the self-adjoint positive operator $L$ on $\LL^2:= L^2(M, \mu)$.
As elsewhere we shall be dealing with operators of the form $f(\sqrt L)$.
We denote by $F_\lambda$, $\lambda\ge 0$, the spectral resolution associated with $\sqrt L$,
that is, $F_\lam=E_{\lam^2}$.
Then $f(\sqrt L)= \int_0^\infty f(\lambda) d F_\lambda$
and the spectral projectors are defined by
$E_\lambda = \ONE_{[0, \lambda]}(L) := \int_0^\infty \ONE_{[0,\lambda]}(u) dE_u$
and
\begin{equation}\label{spect-projector}
F_\lambda = \ONE_{[0,  \lambda]}(\sqrt L)
:=\int_0^\infty  \ONE_{[0,\lambda]}(u)dF_u
=\int_0^\infty  \ONE_{[0,\lambda]}(\sqrt u)dE_u.
\end{equation}

We next list some properties of $F_\lambda$ which follow readily from Theorem~\ref{thm:rough-kernels}:
The~operator $F_\lambda$ is a kernel operator whose kernel $F_\lambda (x,y)$ is
a real symmetric nonnegative function on $M\times M$.
Also,
\begin{equation}\label{loc-ker-F}
F_\lambda (x,y)\le c|B(x, \lambda^{-1})|^{-1/2}|B(y, \lambda^{-1})|^{-1/2}
\end{equation}
and $F_\lambda (x,y)$ is in Lip $\alpha$ for some $\alpha>0$, see (\ref{rough2}).
The mapping property of $F_\lam$ on $\LL^p$ spaces is given by
\begin{equation}\label{ker-F-norm}
\| F_\lambda f \|_q \le c\lambda^{d(1/p-1/q)} \| f \|_p,
\quad 1\leq p \leq 2 \leq q \leq \infty.
\end{equation}

We define the \textit{spectral spaces} $\Sigma^p_\lam$ for $1\le p \le 2$ by
$$
\Sigma^p_\lambda = \{ f \in \LL^p: \; F_\lambda f = f\}.
$$
Notice that $F_\lambda$ is not necessarily a continuous operator on $\LL^p$ if $p>2$
and, therefore, $\Sigma^p_\lambda$ cannot be defined as above for $2<p \le \infty$.
Instead, we shall use the following characterization of $\Sigma^p_\lambda$:
{\em
A function $f \in \Sigma^p_\lambda$ for $1\le p \le 2$ if and only if
$\theta(\sqrt L) f =f$ for all $\theta\in C^\infty_0(\bR_+)$
such that $\theta \equiv 1$ on $[0,\lam]$.
}
This characterization follows by the fact that
$\Sigma^p_\lam \subset \Sigma^2_\lam$ for $1\le p \le 2$
and the boundedness of the operator $\theta(\sqrt L)$ with $\theta$ as above.


\begin{definition}\label{def:Sigma-p}
For $1\le p \le \infty$ we define
$$
\Sigma^p_\lambda
:= \{f\in \LL^p: \theta (\sqrt L)f = f \hbox{ for all }
\theta \in C^\infty_0(\bR_+), \;  \theta \equiv 1 \hbox{ on } \;  [0, \lambda]\}.
$$
Furthermore, for any compact $K \subset [0, \infty)$ we define
$$
\Sigma^p_K
:= \{f\in \LL^p: \theta (\sqrt L)f = f \hbox{ for all }
\theta \in C^\infty_0(\bR_+), \;  \theta \equiv 1 \hbox{ on } \;  K\}.
$$
\end{definition}


\begin{proposition}\label{prop:spectral-intersect}
For any $\lambda \ge 1$ and $1\le p \le \infty$
\begin{equation}\label{spectral-intersect}
\Sigma_\lambda^p = \cap_{\varepsilon >0} \Sigma_{\lambda+\varepsilon}^p.
\end{equation}
\end{proposition}

\noindent
{\bf Proof.}
Suppose  $f \in \cap_{\epsilon >0} \Sigma_{\lambda+\epsilon}^p$
and let
$\theta \in C^\infty_0(\bR_+)$, $\supp\theta \subset [0,R]$, and
$\theta \equiv 1$ on $[0, \lambda]$.
By Definition~\ref{def:Sigma-p}
$f= \theta(r^{-1}\sqrt L)f$ for each $r >1$ and hence
\begin{equation}\label{theta-theta-1}
\|f- \theta(\sqrt L)f\|_p = \|\theta(r^{-1}\sqrt L)f - \theta(\sqrt L)f \|_p,
\quad r>1.
\end{equation}
Assuming that $1<r\le 2$, Theorem~\ref{thm:main-local-kernels} implies
$$ 
|\theta(r^{-1}\sqrt L)(x, y) - \theta(\sqrt L)(x, y)| \le C_r D_{1,k}(x,y),
$$ 
where
$C_r = c_k R^{2k+d+4}\big(\|\theta(r^{-1}\cdot)- \theta(\cdot)\|_\infty
+ \|(d/d\lambda)^{2k+4}[\theta(r^{-1}\cdot)- \theta(\cdot)]\|_\infty\big)$.
We now choose $k\ge 2d+1$ and apply Proposition~\ref{prop:young} to obtain
\begin{equation}\label{theta-theta-2}
\|\theta(r^{-1}\sqrt L)- \theta(\sqrt L)) \|_{p\rightarrow p} \le cC_r.
\end{equation}
Clearly, for any $\nu\ge 0$ we have
$\lim_{r\to 1}\|(d/d\lambda)^\nu [\theta(r^{-1}\cdot)- \theta(\cdot)]\|_\infty = 0$
and, therefore,
$\lim_{r\to 1} C_r =0$.
This along with (\ref{theta-theta-1})-(\ref{theta-theta-2}) yields
$\|f- \theta(\sqrt L)f\|_p=0$, which completes the proof.
$\qed$

With the next claim we establish a Nikolski's type inequality that relates different
$\LL^p$-norms on spectral spaces.


\begin{proposition}\label{prop:Nikolski}
If $1\le p \le q \le \infty$, then
$\Sigma_\lam^p \subset \Sigma_\lambda^q$, $\Sigma_\lambda^q  \cap \bL^p =\Sigma_\lambda^p$, and
there exists a constant $c>0$ such that
\begin{equation}\label{norm-relation1}
\|g\|_q \le c\lam^{d(1/p-1/q)}\|g\|_p,
\quad g \in \Sigma_\lam^p, \; \lam  \ge 1.
\end{equation}
Furthermore, for any $g \in \Sigma_\lam^\infty$, $\lam  \ge 1$,
\begin{equation}\label{lip2}
|g(x)- g(y)| \le c \big(\lambda \rho(x,y)\big)^\alpha \|g\|_\infty,
\quad x,y \in M,
\end{equation}
with $\alpha>0$ the constant from $(\ref{lip})$.
\end{proposition}

\noindent
{\bf Proof.}
Let $g \in \Sigma_\lam^p$, $\lam\ge 1$, and set $\delta:= \lam^{-1}$.
Choose $\theta \in C^\infty_0(\bR_+)$ so that $\theta \equiv 1$ on $[0,1]$.
Then $g=\theta(\delta\sqrt L)g$
and (\ref{norm-relation1})-(\ref{lip2}) follow readily by Corollary~\ref{lip0}.
$\qed$

\subsection{Linear approximation from spectral spaces}\label{linear-app}

The purpose of this subsection is to give a~short account of linear approximation
from $\Sigma_t^p$ in $\LL^p$, $1\le p\le \infty$.
Let $\cE_t(f)_p$ denote the best approximation of $f \in \LL^p$ ($\LL^\infty:=\UCB$) from $\Sigma_t^p$,
that is,
\begin{equation}\label{def:best-app}
\cE_t(f)_p:=\inf_{g\in\Sigma_t^p}\|f-g\|_p.
\end{equation}

Our goal is to characterize the approximation space $A_{pq}^s$, $s>0$, $0<q\le\infty$,
defined as the set of all functions $f\in \LL^p$ such that
\begin{equation}\label{def:app-space-p}
\|f\|_{A_{pq}^s}
:= \|f\|_p +
\Big(\sum_{j\ge 0}\big(2^{s j}\cE_{2^j}(f)_p \big)^q\Big)^{1/q}<\infty
\quad\hbox{if}\;\; q<\infty,
\quad\hbox{and}
\end{equation}
\begin{equation}\label{def:app-space-infty}
\|f\|_{A^s_{p \infty}}
:= \|f \|_p + \sup_{j\ge 0}  2^{sj} \cE_{2^j}(f)_p <\infty
\quad\hbox{if}\;\; q=\infty.
\end{equation}
Due to the monotonicity of $\cE_t(f)_p$ we have
$\|f\|_{A_{pq}^s}\sim\|f\|_p+\Big(\int_1^\infty(t^{s}\cE_t(f)_p)^qdt/t\Big)^{1/q}$,
when $q<\infty$, and
$\|f\|_{A^s_{p \infty}}
:= \|f \|_p + \sup_{t\ge 1} t^s \cE_t(f)_p <\infty
\quad\hbox{if}\;\; q=\infty.$

To characterize $A_{pq}^s$ we shall use the well-known machinery of Bernstein and Jackson
estimates and interpolation.
In \S\ref{sec:char-Besov} it will be shown that $A_{pq}^s$ can be identified as a certain Besov space.

\subsubsection{Bernstein and Jackson estimates. Characterization of spectral spaces}
\label{sec:bernstein-char}

We begin by proving a Bernstein estimate.


\begin{theorem}\label{thm:Bernstain}
Let $1\le p \le \infty$ and $m \in \bN$.
Then there exists a constant $c^\star=c^\star(m)>0$, independent of $p$,
such that for any $g\in \Sigma_\lam^p$, $\lam\ge 1$,
\begin{equation}\label{bernstein}
\|L^m g\|_p \le c^\star\lam^{2m}\|g\|_p.
\end{equation}
\end{theorem}

\noindent
{\bf Proof.}
As in the proof of Proposition~\ref{prop:Nikolski},
pick $\theta \in C^\infty_0(\bR_+)$ so that $\theta \equiv 1$ on $[0,1]$.
Then for any $g \in \Sigma_\lam^p$ we have
$g=\theta(\delta\sqrt L)g$ with $\delta:= \lam^{-1}$ and, therefore,
$L^m g = L^m\theta(\delta\sqrt L)g$.
Then (\ref{bernstein}) follows by applying Corollary~\ref{cor:LS-local-kernels}
and Proposition~\ref{prop:young}.
$\qed$

\smallskip

Observe that from  spectral theory it readily follows that when $p=2$
the Bernstein estimate (\ref{bernstein}) holds with constant $c^\star=1$.

Our next aim is to show that the spectral spaces $\Sigma_\lambda^p$ can be characterized
by means of Bernstein estimates, in the spirit of the previous theorem,
but with a constant ($c_\nu$ below) independent of $m$.


\begin{theorem}\label{thm:BernChar}
Let $1\leq p \leq \infty$ and $\lambda > 0$.
Then the following assertions are equivalent:

$(a)$ $f \in \Sigma_\lambda^p$.

$(b)$ $f \in \cap_{m\in \bN} \Dom(L^m)$ and for any $\nu >\lambda$
there exists a constant $c_\nu>0$ such that
$$
\|L^m f\|_p \le c_\nu \nu^{2m}\|f\|_p,
\quad\forall m\ge 1.
$$

$(c)$ $$ z\in \bC \mapsto e^{-zL}f = \sum_{k \ge 0} \frac{(-z)^k}{k!} L^k f$$
is an entire function of exponential  type $\lambda^2$.
\end{theorem}

\noindent
{\bf Proof.}
Clearly, (b) $\Longleftrightarrow$ (c) using the Paley-Wiener theorem.

To prove that $(a) \Longrightarrow (b)$ we shall show that 
the constant $c^\star$ in (\ref{bernstein}) can be specified as follows:
For any $0<\varepsilon < 1$ there exists a constant $c(\eps, d) >0$ such that
\begin{equation}\label{bernstein-const}
c^\star=c(\eps, d) m^{4d+8}(1+\eps)^{2m}.
\end{equation}
Indeed, let $\theta \in C^\infty_0(\bR)$ be so that $\theta \equiv 1$ on $[-1,1]$,
$\supp \theta\subset [-1-\eps, 1+\eps]$, and also $0\le \theta\le 1$.
With $\delta:=\lambda^{-1}$ we have $f= \theta(\delta \sqrt L)f$ for any $f \in \Sigma_\lambda^p$
and we shall estimate $\|L^m\theta(\delta \sqrt L)f\|_p$.
Denote briefly $h(u):= u^{2m}\theta(u)$.
Then $h(\delta \sqrt L)=\delta^{2m}L^m\theta(\delta \sqrt L)$.
To go further,
set $k:= \lfloor 2d\rfloor+2$, hence $2d+1<k\le 2d+2$.
It is readily seen that
$$
\|h\|_\infty \le (1+\eps)^{2m}
\quad \hbox{and}\quad
\|h^{(2k+4)}\|_\infty \le c_1(\eps, d)m^{4d+8}(1+\eps)^{2m}.
$$
Now, by Theorem~\ref{thm:main-local-kernels} we infer
$$
|L^m\theta(\delta \sqrt L)(x, y)| = \delta^{-2m}|h(\delta \sqrt L)(x, y)|
\le c_2(\eps, d) m^{4d+8}(1+\eps)^{2m} \lam^{2m}D_{\delta, k}(x, y)
$$
and applying Proposition~\ref{prop:young} $(k> 2d+1)$ we arrive at
$$
\|L^m f\|_p =
\|L^m \theta(\delta \sqrt L)f\|_p \le c(\eps, d) m^{4d+8}(1+\eps)^{2m} \lam^{2m}\|f\|_p
\quad \hbox{for} \;\; f\in \Sigma_\lam^p,
$$
which confirms (\ref{bernstein-const}).

Given $\nu>\lam$, choose $0<\eps<1$ so that $(1+\eps)^2\lambda \le \nu$.
Then from above and the obvious fact that
$\sup_{m\ge 1}m^{4d+8}(1+\eps)^{-2m} \le c'(\eps, d)$ we get
$$
\|L^m f\|_p \le c(\eps, d) m^{4d+8}(1+\eps)^{-2m} \nu^{2m}\|f\|_p
\le c''(\eps, d)\nu^{2m}\|f\|_p \quad\forall f\in \Sigma_\lam^p.
$$
Thus $(a) \Longrightarrow (b)$.


Now, to prove that  $(b) \Longrightarrow (a)$,
suppose (b) holds for some function $f \in \LL^p$ and let
$\theta \in C^\infty_0(\bR_+)$, $\theta \equiv 1$ on $[0, \lambda]$,
as in Definition~\ref{def:Sigma-p}.
Assume $\supp \theta \subset [0, R]$.
Let $\eps>0$. We shall show that
$\|f-\theta(\sqrt L)f\|_p <\eps$, which implies $f\in \Sigma_\lam^p$.
Indeed, for $0<\delta < r <1$ we have
$$
\|f-\theta(\sqrt L)f\|_p
\le \|f-\theta(\delta\sqrt L)f\|_p
+ \|\theta(\delta\sqrt L)f-\theta(r\sqrt L)f\|_p
+\|\theta(r\sqrt L)f-\theta(\sqrt L)f\|_p.
$$
By Proposition~\ref{prop:app-identity}, $\|f-\theta(\delta\sqrt L)f\|_p \to 0$ as $\delta\to 0$
and hence there exists $\delta >0$ such that
$\|f-\theta(\delta\sqrt L)f\|_p < \eps/2$.
Clearly,
$\|\theta(r\sqrt L)f-\theta(\sqrt L)f\|_p \to 0$ as $r\to 1$
and hence there exists $r<1$ such that
$\|\theta(r\sqrt L)f-\theta(\sqrt L)f\|_p <\eps/2$.

It remains to show that
$\|\theta(\delta\sqrt L)f-\theta(r\sqrt L)f\|_p =0$.
Let $\lam<\nu< \lam/r$
and denote briefly
$h(u):= \big[\theta(\delta u)-\theta(r u)\big]u^{-2m}$.
Note that $\supp h \subset [\lam/r, R/\delta]$.
Then using our assumption we have
\begin{align*}
\|\theta(\delta\sqrt L)f-\theta(r\sqrt L)f\|_p
= \|h(\sqrt L)L^m f\|_p
&\le  \|h(\sqrt L)\|_{p\to p}\|L^m f\|_p\\
&\le c_\nu \|h(\sqrt L)\|_{p\to p}\nu^{2m}\|f\|_p,
\quad \forall m\ge 1.
\end{align*}
As above, set $k:= \lfloor 2d\rfloor+2$, then $2d+1<k\le 2d+2$.
Now, applying Theorem~\ref{thm:main-local-kernels} and Proposition~\ref{prop:young} it follows that
\begin{align*}
\|h(\sqrt L)\|_{p\to p}
&\le c(R/\delta)^{2k+d+4}\big[\|h\|_\infty + \|h^{(2k+4)}\|_\infty\big]\\
&\le c' m^{2k+4}(\lam/r)^{-2m}
\end{align*}
and hence
$$
\|\theta(\delta\sqrt L)f-\theta(r\sqrt L)f\|_p \le c m^{4d+8}(r\nu/\lam)^{2m}\|f\|_p.
$$
Here the constant $c$ depends on $\delta, r, R, d, \lam, \nu$, but is independent of $m$.
Since $0< r\nu/\lam <1$ by letting $m\to\infty$ we obtain
$\|\theta(\delta\sqrt L)f-\theta(r\sqrt L)f\|_p=0$.
Therefore, $(b) \Longrightarrow (a)$.
$\qed$

We now establish a Jackson estimate for approximation from $\Sigma_t^p$.


\begin{theorem}\label{thm:Jackson}
Let $1\le p \le \infty$.
Then for any $m\in \bN$ there exists a constant $c_m>0$ such that for any $t\ge 1$
\begin{equation}\label{linear-Jackson}
\cE_t(f)_p \le c_m  t^{-2m}\|L^m f\|_p
\quad\hbox{for}\;\; f \in \Dom (L^m)\cap \LL^p.
\end{equation}
\end{theorem}

\noindent
{\bf Proof.}
Let $\theta\in C^\infty(\bR)$, $\theta(u)=1$ for $u\in [0, 1]$, $0\le \theta\le 1$,
and $\supp \theta\subset [0, 2]$.
Set $\varphi(u):=\theta(u/2)-\theta(u)$.
Then
$1-\theta(u)=\sum_{j\ge 0} \varphi(2^{-j}u)$, $u\in \bR_+$.
Given $t>0$, set $\delta:=2/t$.
Assume $f \in \Dom (L^m)\cap \LL^p$.
Clearly, $\theta(\delta \sqrt L)f \in \Sigma_t^p$ and hence
$$
\cE_t(f)_p
\le \|f-\theta(\delta \sqrt L)f\|_p
\le \sum_{j\ge 0}\|\varphi(2^{-j}\delta \sqrt L)f\|_p.
$$
Denote briefly $h(u):=\varphi(u)u^{-2m}$.
Then
$\varphi(2^{-j}\delta \sqrt L)L^{-m}= (2^{-j}\delta)^{2m}h(2^{-j}\delta \sqrt L)$
and, therefore,
$$
\|\varphi(2^{-j}\delta \sqrt L)f\|_p
\le \|\varphi(2^{-j}\delta \sqrt L)L^{-m}L^m f\|_p
\le (2^{-j}\delta)^{2m}\|h(2^{-j}\delta \sqrt L)\|_{p\to p}\|L^m f\|_p.
$$
By Theorem~\ref{thm:main-local-kernels} and Proposition~\ref{prop:young} it follows that
$\|h(2^{-j}\delta \sqrt L)\|_{p\to p} \le c(d,m)$
and hence
$$
\cE_t(f)_p \le ct^{-2m}\|L^m f\|_p \sum_{j\ge 0}2^{-2mj}
\le c't^{-2m}\|L^m f\|_p,
$$
which gives (\ref{linear-Jackson}).
$\qed$

\subsubsection{Characterization of approximation spaces}

Once the Bernstein and Jackson estimates are established,
the approximation spaces $A^s_{pq}$, defined in (\ref{def:app-space-p})-(\ref{def:app-space-infty}),
can be characterized by interpolation.
In the following we shall denote by $\big(X_0, X_1\big)_{\theta, q}$ the real interpolation space between
the normed spaces $X_0$, $X_1$, see e.g. \cite{BB, BL}.


\begin{theorem}\label{thm:interpol}
Let $s>0$, $1\leq p \leq \infty$ and $0<q\le \infty$.
Then for any $r>s$
\begin{equation}\label{interpol-1}
A^s_{pq} = \big(\bL^p, \D\big(\sqrt{L_{}}\big)^r\big)_{\theta, q}, \quad s= \theta r.
\end{equation}
\end{theorem}

\noindent
{\bf Proof.}
A classical argument (e.g. \cite{DL})  using the Jackson and Bernstein estimates from
(\ref{linear-Jackson}) and (\ref{bernstein}) implies the following characterization of
the spaces $A^s_{pq}$: If~$2m >s$, then
\begin{equation}\label{interpol-2}
A^s_{pq}
= \big(\bL^p, \D(L_{}^m )\big)_{\theta, q}
= \big(\bL^p, \D(\sqrt{L_{}})^{2m}\big)_{\theta, q},
\quad s= 2\theta m.
\end{equation}
Thus (\ref{interpol-1}) holds for $r=2m$.
On the other hand, $-\sqrt{L_{}}$ is the infinitesimal generator of the subordinate semigroup
$
Q_tf = \int_0^\infty \frac{t e^{-t^2/4s}}{2s \sqrt{\pi s}} e^{-sL}f d\mu(s)
$
on $ \bL^p$,
and by a well-known result (e.g. \cite{BB})
if $1\leq r < k$, 
then
$$
\big(\LL^p, \D\big(\sqrt{L_{}}\big)^k\big)_{\theta, 1}
\subset   \D\big(\sqrt{L_{}}\big)^r
\subset  \big(\bL^p, \D\big(\sqrt{L_{}}\big)^k \big)_{\theta, \infty}, \quad \theta = r/k.
$$
Therefore, if $1\le r <2m$ and $\theta_0 =\frac r{2m}$, then
$$
A^{r}_{p1} =\big(\bL^p, \D\big(\sqrt{L_{}}\big)^{2m}\big)_{\theta_0, 1}
\subset   \D\big(\sqrt{L_{}}\big)^r
\subset  \big(\bL^p, \D\big(\sqrt{L_{}}\big)^{2m}\big)_{\theta_0, \infty}=A^{r}_{p\infty}
$$
This along with (\ref{interpol-2}) implies
$$
\big(\LL^p, \big(\bL^p, \D\big(\sqrt{L_{}}\big)^{2m}\big)_{\theta_0, 1}\big)_{\theta, q}
\subset
\big(\LL^p, \D\big(\sqrt{L_{}}\big)^{r}\big)_{\theta, q}
\subset
\big(\LL^p, \big(\LL^p, \D\big(\sqrt{L_{}}\big)^{2m}\big)_{\theta_0, \infty}\big)_{\theta, q}
$$
and by the reiteration theorem (e.g. \cite{BL}) this leads to
$$
\big(\LL^p, \D\big(\sqrt{L_{}}\big)^{r}\big)_{\theta, q}
= \big(\bL^p, \D\big(\sqrt{L_{}}\big)^{2m}\big)_{\theta \theta_0, q}= A^{s}_{pq},
\quad
s=2\theta \theta_0 m= \theta r.
$$
The proof is complete.
$\qed$


\begin{remark}\label{semig}
From the above, $A^s_{pq}=\big(\bL^p, D(L^m )\big)_{\theta, q}$,
$s= 2\theta m$, $ 0< s <2m$, but then as is well-known $($e.g. \cite{BB}$)$
$$
\| f \|_{A^s_{pq}} \sim \| f \|_p +
\Big(\int_0^1  \big(t^{-s/2} \| (e^{-tL} - \Id)^m f \|_p \big)^q \frac{dt}t\Big)^{1/q}
$$
with the usual modification for $q=\infty$.
Moreover, since $e^{-tL} $ is a holomorphic semigroup, we also have
$$
\|f\|_{A^s_{pq}} \sim \| f \|_p +
\Big(\int_0^1 \big(t^{-s/2} \| (tL)^m e^{-tL}f\|_p \big)^q \frac{dt}t\Big)^{1/q}
$$
with the usual modification for $q=\infty$.
\end{remark}

\subsection{Kernel norms}

Here we derive bounds on the $\LL^p$-norms of the kernels of operators of the form
$\theta (\delta \sqrt L)$, which will be important for the development of frames.


\begin{theorem}\label{thm:norms}
Let $\theta\in C^\infty(\bR_+)$, $\theta \ge 0$,
$\supp \theta\subset [0, R]$ for some $R>1$,
and $\theta^{(2\nu+1)}(0)=0$, $\nu=0, 1, \dots$.
Suppose that either

$(i)$\;\;
$\theta(\uu)\ge 1$ for $\uu\in [0, 1]$, or

$(ii)$
$\theta(\uu)\ge 1$ for $\uu\in [1, b]$, where $b>1$ is a sufficiently large constant.

\noindent
Then for $0 < p \le \infty$, $0<\delta \le \min\{1, \frac{\diam M}{3}\}$, and $x \in M$ we have
\begin{equation}\label{est-norm-1}
c_1|B(\xi,\delta)|^{1/p-1}
\le \| \theta (\delta \sqrt L) (x,.)\|_p
\le c_2|B(x,\delta)|^{1/p-1},
\end{equation}
where $c_1>0$ depends only on $p$ and the parameters of the space,
and $c_2>0$ depends on $p$ and the smoothness and the support of $\theta$ similarly as
in Theorem~\ref{thm:main-local-kernels}.
\end{theorem}

\noindent
{\bf Proof.}
By Theorem~\ref{thm:main-local-kernels} we have
$|\theta (\delta \sqrt L) (x,y)| \le c_\sig D_{\delta,\sig}(x,y)$
for any $\sig >0$.
Pick $\sig > d(1/2+1/p)$. Then the upper bound estimate in (\ref{est-norm-1})
follows readily by estimate (\ref{INT}).

It is not hard to see that to prove the lower bound estimate in (\ref{INT}) it suffices
to have it for $p=2$ and $p=\infty$ and use the already established upper bound.
However, clearly
$$
\|\theta (\delta \sqrt L) (x,.)\|_2^2 = \theta^2(\delta \sqrt L) (x,x)
\quad\hbox{and}\quad
\|\theta (\delta \sqrt L) (x,.)\|_\infty \ge \theta(\delta \sqrt L) (x,x),
$$
and it boils down to establishing lower bounds on
$\theta^2(\delta \sqrt L) (x,x)$ and $\theta(\delta \sqrt L) (x,x)$.

Further, let $f, g\in \bL^\infty(\bR_+)$ be bounded, $\supp f, g\subset [0, R]$, and $0\le g\le f$.
Then $f=g+h$ for some $h \geq 0$, and hence
$ f(\sqrt L)(x,x)= g(\sqrt L)(x,x)  +  h(\sqrt L)(x,x)$.
On the other hand, by (\ref{rough4}) 
$f(\sqrt L)(x,x)  = \int_M  |\sqrt f (\sqrt L)(x,y)|^2 d\mu(y)\ge 0$,
and we have similar representations of $g(\sqrt L)(x,x)$ and $h(\sqrt L)(x,x)$.
Therefore,
\begin{equation}\label{lem-norm-3}
0\le g\le f
\quad\Longrightarrow\quad
0 \le g(\sqrt L)(x,x) \le f(\sqrt L)(x,x).
\end{equation}
This allows to compare the kernels of different operators
and we naturally come to the next lemma which is interesting in their own right.


\begin{lemma}\label{lem:norms}
$(a)$
There exist constants $c_3, c_4>0$ such that for any $\tt\ge 1$
\begin{equation}\label{lem-norm-1}
c_3|B(x,\tt^{-1})|^{-1}
\le \ONE_{[0, \tt]}\big(\sqrt{L}\big)(x,x))
\le c_4|B(x,\tt^{-1})|^{-1}.
\end{equation}

$(b)$
There exists $b >1$ such that if $\tau\ge 1$ and $\tau^{-1}\le \frac{\diam M}{3}$, then
\begin{equation}\label{lem-norm-2}
 c_5|B(x,\tt^{-1})|^{-1} \le \ONE_{[\tt, b\tt]}(\sqrt{L})(x,x) \le c_6|B(x,\tt^{-1})|^{-1},
\end{equation}
where $c_5, c_6 >0$ depend only on the parameters of the space.
\end{lemma}

\noindent
{\bf Proof.}
We first show that
\begin{equation}\label{lem-norm-4}
p_t(x,y) = \lim_{\tt \to \infty} \ONE_{[0,\tt]}(\sqrt L) p_t(x,y),
\quad t>0.
\end{equation}
Indeed, we have
$$
\ONE_{[0,\tt]}(\sqrt L) e^{-t L} + \ONE_{(\tt , \infty)}(\sqrt L)e^{-t L} = e^{-t L},
$$
and since
$\ONE_{[0, \tt]}(\sqrt L)e^{-t L}$  is a kernel operator (Theorem~\ref{thm:rough-kernels}), then
$\ONE_{(\tt , \infty)}(\sqrt L)e^{-t L}$ is also a~kernel operator and
\begin{equation}\label{lem-norm-5}
\ONE_{[0,\tt]}(\sqrt L) p_t(x,y)+ \ONE_{(\tt , \infty)}(\sqrt L)p_t(x,y)
= e^{-t L}(x,y).
\end{equation}
On the other hand,
$$\ONE_{(\tt , \infty)}(\sqrt L)e^{-t L}
= e^{-\frac t4 L} \big[\ONE_{(\tt , \infty)}(\sqrt L)e^{-\frac t2 L}\big] e^{-\frac t4 L},$$
and by  spectral theory
$
\|\ONE_{(\tt , \infty)}(\sqrt L)e^{-\frac t2 L}\|_{2\rightarrow 2}  = e^{-\frac t2 \tt^2}.
$
Therefore, applying Proposition~\ref{prop:prod-oper} we arrive at
$$
\ONE_{(\tt , \infty)}(\sqrt L)p_t(x,y)
\le c\frac{e^{-\frac t2 \tt^2}}{\sqrt{|B(x,\sqrt{t/2})|  |B(y,\sqrt{t/2})|}}
\to 0 \quad\hbox{as} \;\; \tt\to\infty.
$$
This and (\ref{lem-norm-5}) imply (\ref{lem-norm-4}).


We also need these bounds on the heat kernel:
\begin{equation}\label{lem-norm-6}
c'|B(x,\sqrt t)|^{-1}
\le p_t(x, x)
\le c|B(x,\sqrt t)|^{-1}, \quad 0<t\le 1.
\end{equation}
The upper bound is immediate from (\ref{Gauss-local}).
For the lower bound we have for $\ell >1$, using (\ref{hol3}),
\begin{align*}
p_t(x,x)
&= \int_M  [p_{t/2}(x,y) ]^2 d\mu(y)
\ge \int_{B(x, 2^\ell\sqrt t)} [p_{t/2}(x, y) ]^2 d\mu(y)\\
&\ge \frac 1{|B(x, 2^\ell\sqrt t)|} \Big[\int_{B(x, 2^\ell\sqrt t)} p_{t/2}(x,y) d\mu(y)\Big]^2\\
&\ge \frac{2^{-\ell d}}{|B(x, \sqrt t)|}
\Big[1- \int_{M\setminus B(x, 2^\ell\sqrt t)} p_{t/2}(x,y) d\mu(y)\Big]^2.
\end{align*}
However, by (\ref{Gauss-local})
$p_{t/2}(x,y) \le c_\sigma D_{\sqrt t, \sigma}(x, y)$ for any $\sigma >0$
hence, just as in the proof of Proposition~\ref{prop:app-identity},
$$
\int_{M\setminus B(x, 2^\ell\sqrt t)} e^{-t/2 L}(x,y) d\mu(y)
\le c2^{-\ell}\le \frac 12
$$
for a sufficiently large $\ell$ (the constant $c$ is independent of $\ell$).
This completes the proof of the lower bound estimate in (\ref{lem-norm-6}).


We now turn to the proof of (\ref{lem-norm-1}).
Since
$\ONE_{[0, \tt]}(u) \le e e^{- \tt^{-2}u^2}$
we obtain, using (\ref{lem-norm-3}) and (\ref{lem-norm-6})
$$
\ONE_{[0, \tt]} (\sqrt{L})(x,x)
\leq e e^{- \tt^{-2}(\sqrt L)^2} (x,x)
\le c|B(x,\tt^{-1})|^{-1},
$$
which gives the right-hand side estimate in (\ref{lem-norm-1}).

For the proof of the left-hand side estimate in (\ref{lem-norm-1}),
we first note that for any $t>0$
\begin{align*}
e^{- tu^2}
&= \ONE_{[0, \tt]} (u)e^{- tu^2} + \sum_{k\geq 0} \ONE_{(2^k\tt,2^{k+1} \tt]} (u)e^{- tu^2}\\
&\leq  \ONE_{[0, \tt]} (u) + \sum_{k\geq 0}  \ONE_{[0,2^{k+1} \tt]} (u) e^{- t2^{2k}\tt^2}.
\end{align*}
From this, (\ref{lem-norm-3}), (\ref{lem-norm-4}), (\ref{lem-norm-6}),
and the right-hand side estimate in (\ref{lem-norm-1}) we obtain
\begin{align*}
c'|B(x,\sqrt t)|^{-1}
&\le p_t(x,x)\\
&\le \ONE_{[0, \tt]} (\sqrt L)(x,x) + \sum_{k\geq 1}
\ONE_{[0,2^{k+1} \tt]} (\sqrt L)(x,x) e^{- t2^{2k}\tt^2}\\
&\le \ONE_{[0, \tt]} (\sqrt L)(x,x) +
c_4\sum_{k\geq 1}e^{- t2^{2k}\tt^2}|B(x,2^{-k-1}\tt^{-1})|^{-1} \\
&\leq \ONE_{[0, \tt]} (\sqrt L)(x,x) +
c_4|B(x, \tt^{-1})|^{-1}\sum_{k\geq 1}e^{- t2^{2k}\tt^2} 2^{(k+1)d}.
\end{align*}
Here for the latter inequality we used (\ref{D1}).
Given $\tt\ge 1$ and $r\in \bN$ we choose $t$ so that $\tt\sqrt t = 2^r$.
Then from above
\begin{align*}
\frac{c'2^{-rd}}{|B(x,\tt^{-1} )|}
\leq  \frac{c'}{|B(x,\sqrt t)|}
&\le \ONE_{[0, \tt]} (\sqrt L)(x,x) +  \frac{ c_42^d2^{-rd}}{|B(x,\tt^{-1})|}
\sum_{k\ge 0}   e^{- 2^{2k} 2^{2r}} 2^{(k+r)d}\\
&\leq \ONE_{[0, \tt]} (\sqrt L)(x,x) +  \frac{c_42^{d}2^{-rd}}{|B(x,\tt^{-1})|}
\sum_{k\ge r}   e^{- 2^{2k} } 2^{kd} .
\end{align*}
Hence,
$$
\frac{2^{-rd}}{|B(x,\tt^{-1})|}
\Big(c'-  c_42^{d}\sum_{k\geq r} e^{- 2^{2k}}2^{kd}\Big)
\le \ONE_{[0, \tt]}(\sqrt L)(x,x).
$$
Taking $r\in \bN$ sufficiently large, this implies the left-hand side estimate in (\ref{lem-norm-1}).


We now take on (\ref{lem-norm-2}).
The right-hand side estimate follows by from the right-hand side estimate in (\ref{lem-norm-1}).
Using (\ref{lem-norm-1}) and the reverse doubling condition (\ref{reverse-doubling})
with $\tau^{-1} \le  \frac{\diam M}{3}$,
we obtain for $l\in \bN$
\begin{align*}
\ONE_{[\tt, 2^l\tt]}(\sqrt{L})(x,x)
&= \ONE_{[0, 2^l\tt]}(\sqrt{L})(x,x)- \ONE_{[0, \tt]} (\sqrt{L})(x,x)\\
&\ge \frac{c_3}{|B(x,2^{-l}\tt^{-1})|}- \frac{c_4}{|B(x,\tt^{-1})|}
\ge \frac{c_3 2^{l\beta}-c_4}{|B(x,\tt^{-1})|},
\end{align*}
which leads to (\ref{lem-norm-2}) with $b=2^l$ for sufficiently large $l$.
$\qed$

\smallskip


\noindent
{\em Completion of the proof of Theorem~\ref{thm:norms}.}
We now focus on the left-hand side estimate in (\ref{est-norm-1}).
Suppose $\theta$ obeys condition~(ii) from the hypothesis of the theorem,
i.e. $\theta(u) \ge 1$ on $[1, b]$, where $b>1$ is the same as in Lemma~\ref{lem:norms}, (b)
(the proof in the other case is the same).
Then by (\ref{lem-norm-3}) and Lemma~\ref{lem:norms} we have
for $0<\delta \le \min\{1, \frac{\diam M}{3}\}$
\begin{align*}
\|\theta (\delta \sqrt L) (x,.)\|_\infty
&\ge \theta (\delta \sqrt L)(x,x)
\geq \ONE_{[1,b]}(\delta \sqrt L) (x,x)\\
&=\ONE_{[\delta^{-1}, \delta^{-1}b]}(\sqrt L) (x,x)
\ge c_5|B(x, \delta )|^{-1}.
\end{align*}
On the other hand
\begin{align*}
\|\theta (\delta \sqrt L) (x,.)\|_2^2
&=\theta^2 (\delta \sqrt L) (x,x)
\geq c_5|B(x, \delta )|^{-1},
\end{align*}
where for the last estimate we proceeded as above.
Thus so far we have
\begin{equation}\label{thm-norm-2}
\begin{aligned}
& \|\theta (\delta \sqrt L) (x,.)\|_p \le c_2|B(x, \delta )|^{1/p-1},
\quad 0<p\le \infty,\\
&\|\theta (\delta \sqrt L) (x,.)\|_\infty \ge c_5|B(x, \delta )|^{-1}
\quad\hbox{and}\quad
\|\theta (\delta \sqrt L) (x,.)\|_2^2 \ge c_5|B(x, \delta )|^{-1}.
\end{aligned}
\end{equation}
Now, for $0<p<\infty$ the left-hand side estimate in (\ref{est-norm-1}) follows from
the estimates in (\ref{thm-norm-2}) in a~standard manner.
Indeed, set $f:=\theta (\delta \sqrt L) (x,.)$.
If~$0<p<2$, then using (\ref{thm-norm-2}) we get
$$
c_5|B(x, \delta )|^{-1} \le \|f\|_2^2 \le \|f\|_p^p\|f\|_\infty^{2-p}
\le c\|f\|_p^p|B(x, \delta )|^{-2+p},
$$
which implies $\|f\|_p \ge c'|B(x, \delta )|^{1/p-1}$.
If $2<p<\infty$, we use (\ref{thm-norm-2}) and H\"{o}lder's inequality to obtain
$$
c_5|B(x, \delta )|^{-1} \le \|f\|_2^2 \le \|f\|_p\|f\|_{p'}
\le c\|f\|_p^p|B(x, \delta )|^{1/p'-1}
\quad (1/p+1/p'=1).
$$
This leads again to $\|f\|_p \ge c'|B(x, \delta )|^{1/p-1}$.
$\qed$

\subsection{Finite dimensional spectral spaces}\label{finite-dim-spaces}

It is easy to see that in the case when $\mu(M)<\infty$ the spectrum of $L$ is discrete and
the respective eigenspaces are finitely dimensional. This and some other related simple
facts are collected in the following statement,
where we adhere to the notation from the previous subsections.


\begin{proposition}\label{prop:finite-dimension}
The following claims are equivalent:

$(a)$
$\diam M <\infty$.

$(b)$
$\mu (M) <\infty$.

$(c)$ There exists $\delta>0$ such that
$\int_M \mu(B(x,\delta))^{-1} d\mu(x) < \infty$
and hence we have
$\int_M \mu(B(x,r))^{-1} d\mu(x) < \infty$
for all $r>0$.

$(d)$ The spectrum of the operator $L$ is discrete and of the form
$0 \le \lam_1<\lam_2<\dots$,
$$
\bL^2 = \sum\bigoplus_{j} \cH_{\lambda_j},
\;\;\hbox{where}\;\; \cH_{\lambda_j} = {\rm Ker}\,(L-\lambda_j \Id),
\quad\hbox{and}\;\; \dim(\cH_{\lambda_j}) <\infty.
$$

$(e)$ There exists $t>0$ such that
$$
\|e^{-tL}\|_{HS}^2=\int_M\int_M |p_t(x, y)|^2 d\mu(x) d\mu(y) =\int_M p_{2t}(x, x) d\mu(x) <\infty,
$$
and hence this is true for all $t>0$.

$(f)$ There exists
$\lambda \ge 1$ $($and hence $\forall \lambda \ge 1)$
$\Sigma_\lambda^\infty = \Sigma_\lambda^1$  $(= \Sigma_\lambda^p$ for all $1\leq p \leq \infty)$.

\smallskip

Furthermore, if one of the above holds, then for $\lam \ge 1$
\begin{equation}\label{finite-dim-2}
\dim (\Sigma_\lam) \sim  \int_M \mu(B(x,\lam^{-1}))^{-1} d\mu(x)
\quad\hbox{and}\quad
\dim (\Sigma_{\sqrt \lam} )\sim\| e^{-\lam L}\|^2_{HS},
\end{equation}
where $\Sigma_\lam =  \sum\bigoplus_{ \sqrt{\lambda_j} \le \lam } \cH_{\lambda_j}$.
In addition,
\begin{equation}\label{finite-dim-3}
p_t(x,y)= \sum_{j\ge 1} e^{-\lambda_j}P_{\cH_j}(x,y) , \quad
P_{\cH_j}(x,y)= \sum_{l=1}^{\dim(\cH_j)} e_j^l(x)\overline{e_j^l(y)},
\end{equation}
where $\{e_j^l: l=1,\dots, \dim(\cH_j)\}$ is an orthonormal basis for
$\cH_j$, $Le_j^l = \lambda_j e_j^l$.
The~convergence is uniform and $p_t(x,y)$ is a positive definite kernel.
\end{proposition}

\noindent
{\bf Proof.}
As already shown in Proposition~\ref{prop:prop-doubl-space}, (a) and (b) are equivalent.
Note that, since in our setting closed balls are compact,
(a) or (b) is also equivalent to the compactness of $M$.

Clearly $(b)$ implies $(f)$ as
$ \Sigma_\lambda^1 \subset \Sigma_\lambda^\infty \subset \bL^\infty  \subset \bL^1 $ and
$\Sigma_\lambda^\infty \cap \bL^1 =\Sigma_\lambda^1$.

To show that $(f)$ implies $(b)$,
assume  $ \Sigma_\lambda^\infty = \Sigma_\lambda^1$.
Then if $\theta \in C_0^\infty(\bR_+)$, $\theta \equiv 1$ in the neighborhood of $0$
and $\supp \theta \subset [0,\lambda]$
we have
$\theta(\sqrt L)f \in  \Sigma_\lambda^\infty = \Sigma_\lambda^1$
$\forall f \in \bL^\infty$.
Hence $1= \theta(\sqrt L)(1) \in \bL^1$, which implies $\mu(M) <\infty.$

Assume that (a)-(b) hold and fix $x_0\in M$.
Then using (\ref{D1})-(\ref{D2}) we get
$$
|B(x_0, 1)| \le 2^d (1+\rho(x_0, x))^d|B(x, 1)| \le (4/\delta)^d (1+\rho(x_0, x))^d|B(x, \delta)|,
\; 0<\delta\le 1,
$$
which readily implies
$$
\int_M |B(x,\delta)|^{-1} d\mu(x)
\le (4/\delta)^d|B(x_0, 1)|(1+D)|M|
< \infty.
$$
Thus (a)-(b) imply (c).

For the other direction, assume that (c) holds
and let $\XX_\delta$ be a maximal $\delta$-net on $M$ with a companion disjoint partition
$\{A_\xi\}_{\xi\in \XX_\delta}$ of $M$ as in Proposition~\ref{prop:delta-net}.
Then we use (\ref{D1})-(\ref{D2}) again to obtain
\begin{align*}
\# X_\delta \le 2^d\sum_{\xi\in\XX_\delta}\frac{|A_\xi|}{|B(\xi, \delta)|}
\le 8^d \sum_{\xi\in\XX_\delta}\int_{A_\xi}\frac{1}{|B(x, \delta)|} d\mu(x)
= 8^d\int_M |B(x, \delta)|^{-1}d\mu(x).
\end{align*}
Hence $\# X_\delta <\infty$, which readily implies $ \diam (M) <\infty$.
So, (c) implies (a).

Since $\int_M p_t(x,y)^2 d\mu(y) = p_t(x,x)$,
the equivalence of (c) and (e) is immediate from (\ref{Gauss-local}).

It remains to show that (c) and (d) are equivalent.
Suppose (c) holds true.
Since $E_\lam^2=E_\lam$, we have
\begin{equation}\label{EE}
\int_M |E_\lambda(x, y)|^2 d\mu(y) = \int_M E_\lambda(x, y)E_\lambda(y, x) d\mu(y)
=E_\lam^2(x, x)= E_\lam(x, x) 
\end{equation}
and hence, using Lemma~\ref{lem:norms},
\begin{align*}
\int_M\int_M |E_\lambda(x, y)|^2 d\mu(x) d\mu(y)
&=\int_M E_\lambda(x,x) dx
=\int_M \ONE_{[0, \sqrt{\lambda}]} (\sqrt L)(x,x) d\mu(x)\\
&\le c\int_M |B(x, \lambda^{-1/2})|^{-1} d\mu(x) <\infty,
\quad \lam \ge 1.
\end{align*}
Therefore, $E_\lambda$ ($\lam\ge 1$) is a Hilbert-Schmidt operator on $\LL^2$
and hence its spectrum is discrete.
Suppose
$\{e_j\}_{j\in J}$ is an orthonormal family, verifying $E_\lambda e_j =e_j$, and put
$$
H(x,y) = \sum_{j \in J} e_j(x) \overline{e_j(y)}.
$$
Evidently, $H^2=H$ and as in (\ref{EE})
$
\int_M | H(x,y)|^2 d\mu(y) = H(x,x) = \sum_{j \in J} |e_j(x)|^2.
$
On the other hand $E_\lambda H = H E_\lambda =H$ and hence
\begin{align*}
H(x,x)
&= \int_M E_\lambda(x,y) H(y,x) d\mu(y)\\
&\le \Big(\int_M |E_\lambda(x,y)|^2 d\mu(y)\Big)^{1/2}
\Big(\int_M |H(y,x)|^2 d\mu(y)\Big)^{1/2}\\
&\le \sqrt{E_\lambda (x,x)}\sqrt{H (x,x)}.
\end{align*}
Consequently,
$H(x,x) \le E_\lambda (x,x)$.
Thus
\begin{align*}
\# J &= \int_M \sum_{j \in J} |e_j(x)|^2 d\mu(x)
= \int_M H(x,x) dx \le \int_M E_\lambda(x,x) dx\\
&=\int_M \ONE_{[0, \sqrt{\lambda}]} (\sqrt L)(x,x) dx
\le c\int_M |B(x,\lam^{-1/2})|d\mu(x)<\infty.
\end{align*}
Therefore,
$
\dim (\Sigma_{\sqrt{\lambda}}) \le c\int_M |B(x,\lam^{-1/2})|d\mu(x)<\infty,
$
which shows that (c) implies~(d).

Finally, assume that (d) holds true.
Let $\{e_j\}_{j\in J}$ be an orthonormal basis of $\Sigma_\lam$, $\lam\ge 1$.
Then
$E_\lam(x, y)=\sum_{j\in J} e_j(x)\overline{e_j(y)}$,
where $\# J= \dim(\Sigma_\lam)$.
Now, using Lemma~\ref{lem:norms} we infer
\begin{align*}
c_3\int_M |B(x,\lam^{-1/2})|d\mu(x)
&\le \int_M \ONE_{[0, \sqrt{\lambda}]} (\sqrt L)(x,x) dx\\
& = \int_M E_\lam(x, x)d\mu(x)
= \dim(\Sigma_{\sqrt\lam})<\infty.
\end{align*}
Thus (d) implies (c).

The estimates in (\ref{finite-dim-2}) follow from above.
The last assertion of the theorem is Mercer's theorem (see \cite{FeMe}).
$\qed$

\section{Sampling theorem and cubature formula}\label{sec:sampling-cubature}
\setcounter{equation}{0}

Basic tools for constructing decomposition systems (frames) for various spaces will be
a sampling theorem for $\Sigma_\lambda^p$ and a cubature formula for $\Sigma_\lambda^1$.
In turn these results will rely on the nearly exponential localization of operator kernels
induced by smooth cut-off functions $\varphi$ (Theorem~\ref{thm:main-local-kernels}):
If $\varphi \in C^\infty(\bR_+)$,
 $\supp\varphi \subset [0,b]$, $b>1$, $0\le \varphi \le 1$, and $\varphi = 1$ on $[0,1]$,
then there exists a constant $\alpha>0$ such that
for any $\delta >0$ and $x,y,x' \in M$
\begin{align}
|\varphi(\delta \sqrt{L})(x,y) | &\le \Css D_{\delta, \sig}(x,y) \label{local-Phi}
\quad\hbox{and}\\
|\varphi(\delta \sqrt{L})(x,y)- \varphi(\delta \sqrt{L})(x',y) |
&\leq \Css \Big(\frac{\rho(x,x')}{\delta}\Big)^\alpha D_{\delta, \sig}(x,y),
\; \rho(x,x')\le \delta. \label{Lip-Phi}
\end{align}
Here $\Css>1$ depends on $\varphi$, $\sig$ and the other parameters,
but is independent of $x, y, x'$ and $\delta$.

\smallskip

The main ingredient in our constructions will be the following
Marcinkiewicz-Zygmund inequality for $\Sigma_\lambda^1$,
where maximal $\dd-$nets (see \S\ref{sec:dd-nets}) will be utilized.


\begin{proposition}\label{prop:Marcinkiewicz}
Given $\lambda \geq 1$, let $\XX_\dd$ be a maximal $\dd-$net on $M$ with
$\dd:=\frac \gamma \lambda$, where $0<\gamma \le 1$.
Suppose $\{A_\xi\}_{\xi\in\XX_\dd}$ is a companion disjoint partition of $M$
consisting of measurable sets such that
$B(\xi, \dd/2) \subset  A_\xi \subset B(\xi,\dd)$,  $\xi \in \XX_\dd$.
Then for any $f \in \Sigma_\lambda^p$, $1\leq p <\infty$,
\begin{equation}\label{Marcink1}
\sum_{\xi \in \XX_\dd} \int _{A_\xi}    |f(x) - f(\xi)|^p dx
\leq  [\CPhi  \gamma^{\alpha } \Cdiam]^p \| f \|_p^p,
\end{equation}
and for any $f \in \Sigma_\lambda^\infty$
\begin{equation}\label{Marcink2}
\sup_{\xi \in \XX_\dd} \sup_{x\in A_\xi} |f(x)-f(\xi)|
\le  \CPhi  \gamma^{\alpha }  \Cdiam \|f\|_\infty,
\end{equation}
where $\CPhi$ is the constant from $(\ref{local-Phi})-(\ref{Lip-Phi})$
with $\sig_*:=2d +1$ and $\Cdiam=2^{2d+1}$. 
\end{proposition}

\noindent
{\bf Proof.}
Suppose $\varphi$ is a cut-off function as in (\ref{local-Phi})-(\ref{Lip-Phi}).
Then we have
$f = \int_M  \varphi (\lambda^{-1} \sqrt L) (\cdot,y) f(y)dy$
for $f\in \Sigma_\lambda^p$, $1\leq p \leq \infty$,
and using (\ref{Lip-Phi}) with $\delta=\lambda^{-1}$ we obtain for $1\leq p < \infty$

\begin{align*}
&\sum_{\xi \in \XX_\dd} \int _{A_\xi}|f(x) - f(\xi)|^p dx\\
&\qquad\quad =\sum_{\xi \in \XX_\dd} \int_{A_\xi}
\Big|\int_M  \Big[\varphi (\lambda^{-1} \sqrt L)(x,y)-\varphi (\lambda^{-1}\sqrt L)(\xi,y)\Big]f(y)dy \Big|^p dx\\
&\qquad\quad \leq   \CPhi^p \sum_{\xi \in \XX_dd} \int_{A_\xi}
\Big(\int_M ( \lambda \rho(x,\xi))^\alpha D_{\delta, \sig_*}(x,y )|f(y)|dy\Big)^p dx\\
&\qquad\quad \leq  \CPhi^p\gamma^{\alpha p}\int_M \Big(\int_M D_{\delta, \sig_*}(x,y )|f(y)|dy\Big)^pdx
\leq  [\CPhi  \gamma^{\alpha }\Cdiam]^p \| f \|_p^p,
\end{align*}
where for the last inequality we used Proposition~\ref{prop:young}.
The proof of (\ref{Marcink2}) is similar.
$\qed$

\subsection{The sampling theorem}

The following sampling theorem will play an important role in the sequel.


\begin{theorem}\label{thm:sampling}
Let  $0< \gamma <1$ and
\begin{equation}\label{gamma}
\CPhi  \gamma^{\alpha }\Cdiam \le \frac 12,
\end{equation}
where $\CPhi$ is the constant from $(\ref{Lip-Phi})$ with $\sig_*:=2d +1$
and $\Cdiam=2^{2d+1}$. 
For a~given $\lambda  \geq 1$ let $\XX_\dd$ be a~maximal $\dd-$net on $M$ 
with
$\dd:=\frac \gamma \lambda$ and
suppose $\{A_\xi\}_{\xi\in\XX_\dd}$ is a companion disjoint partition of $M$
consisting of measurable sets such that
$B(\xi, \dd/2) \subset  A_\xi \subset B(\xi,\dd)$,  $\xi \in \XX_\dd$.
Then for any $f\in \Sigma_\lambda^p$, $1\le p<\infty$,
\begin{equation}\label{samp1}
\frac 12 \| f \|_p \le \Big(\sum_{\xi \in \XX_\dd} |A_\xi|| f(\xi)|^p\Big)^{1/p} \le 2 \| f \|_p
\end{equation}
and for $f\in \Sigma_\lambda^\infty$
\begin{equation}\label{samp2}
\frac 12 \| f \|_\infty \le \sup_{\xi \in \XX_\dd} |f(\xi)|\le \| f \|_\infty.
\end{equation}
Furthermore, if $0<\gamma<1$ is selected so that
\begin{equation}\label{gamma-eps}
\CPhi  \gamma^{\alpha }\Cdiam \le \small\frac{\eps}{3},
\end{equation}
$($instead of $(\ref{gamma}))$
for a given $0 < \eps <1$,
then for any $f\in \Sigma_\lambda^p$, $1\le p\le 2$,
\begin{equation}\label{samp3}
(1-\eps)\| f \|_p^p \le \sum_{\xi \in \XX_\dd} |A_\xi|| f(\xi)|^p \le (1+\eps)\| f \|_p^p.
\end{equation}
\end{theorem}

\noindent
{\bf Proof.}
We first prove (\ref{samp3}).
It is easy to see that
\begin{equation}\label{simp-ineq2}
\frac 1{(1+\delta)^{p-1}} |a |^p \leq \frac 1{\delta^{p-1}} |a-b|^p + |b|^p
\quad \hbox{if $0<\delta<1$, $a,b \in \bC$ and $1\le p$.}
\end{equation}
which implies :
\begin{equation}\label{simp-ineq}
(1-\delta) |a |^p \leq \frac 1{\delta^{p-1}} |a-b|^p + |b|^p
\quad \hbox{if $0<\delta<1$, $a,b \in \bC$ and $1\le p\le 2$.}
\end{equation}
This inequality with $\delta:=\eps/3$ implies
\begin{align}
(1-\eps/3) \int_{A_\xi}|f(x)|^p dx
&\leq \frac 1{(\eps/3)^{p-1}}\int _{A_\xi}|f(x)-f(\xi)|^p dx + | A_\xi || f(\xi)|^p,\label{samp-est1}\\
(1-\eps/3)| A_\xi||f(\xi)|^p
&\leq \frac 1{(\eps/3)^{p-1}}\int_{A_\xi}|f(x) - f(\xi)|^p dx + \int_{A_\xi}|f(x)|^p dx.\label{samp-est2}
\end{align}
Summing up estimates (\ref{samp-est1}) over $\xi\in\XX_\dd$, we get
\begin{align*}
(1-\eps/3)\|f\|_p^p
&\le \frac 1{(\eps/3)^{p-1}}\sum_{\xi \in \XX_\dd} \int _{A_\xi}|f(x) - f(\xi)|^p dx
+  \sum_{\xi \in \XX_\dd} | A_\xi || f(\xi)|^p\\
&\le \frac 1{(\eps/3)^{p-1}}  [\CPhi\gamma^{\alpha }\Cdiam]^p \| f \|_p^p
+  \sum_{\xi \in \XX_\dd} | A_\xi || f(\xi)|^p\\
&\le \frac{\eps}{3}\| f \|_p^p + \sum_{\xi \in \XX_\dd} | A_\xi || f(\xi)|^p,
\end{align*}
which implies the left-hand side estimate in (\ref{samp3}).
Here for the second estimate we used (\ref{Marcink1}).

Similarly, we sum up estimates (\ref{samp-est2}) and use again (\ref{Marcink1}) to obtain
\begin{align*}
(1-\eps/3)\sum_{\xi \in \XX_\dd} | A_\xi || f(\xi)|^p
&\le \frac 1{(\eps/3)^{p-1}}\sum_{\xi \in \XX_\dd} \int _{A_\xi}|f(x) - f(\xi)|^p dx
+  \|f\|_p^p\\
&\le \frac 1{\eps^{p-1}}  [\CPhi\gamma^{\alpha }\Cdiam]^p \| f \|_p^p
+   \|f\|_p^p
\le (1+\eps/3)\| f \|_p^p,
\end{align*}
which readily yields the right-hand side estimate in (\ref{samp3}).

To establish (\ref{samp1}) note that (using (\ref{simp-ineq2}))
$\frac 1{2^{p-1}} |a|^p \leq |a-b|^p + |b|^p$ for $a, b\in \bC$ and $1\leq p <\infty$,
which leads to
\begin{align*}
\frac 1{2^{p-1}} \int_{A_\xi}  |f(x)|^p dx
&\leq  \int _{A_\xi}    |f(x) - f(\xi)|^p dx + | A_\xi |    | f(\xi)|^p,\\
\frac 1{2^{p-1}}  | A_\xi |    | f(\xi)|^p
&\leq  \int _{A_\xi}    |f(x) - f(\xi)|^p dx +\int_{A_\xi}  |f(x)|^p dx.
\end{align*}
Then one proceeds exactly as above and obtains (\ref{samp1}).
The proof of (\ref{samp2}) is simpler and will be omitted.
$\qed$


\begin{remark}
Observe that under the assumptions of Theorem~\ref{thm:sampling} one has,
using $(\ref{doubling})$ and $(\ref{D1})$,
$$ 
(4/\gamma)^{-d} |B(\xi, \lambda^{-1})|
\leq 2^{-d} |B(\xi, \gamma \lambda^{-1})|
\leq  |A_\xi| \leq |B(\xi, \gamma \lambda^{-1})|  \leq  |B(\xi, \lambda^{-1})|,
\;\; \xi\in\XX_\dd.
$$ 
Then estimates $(\ref{samp1})$ imply that
for $f\in \Sigma_\lambda^p$, $1\le p<\infty$,
\begin{equation}\label{samp4}
\small \frac 12(\gamma/4)^{d/p}
\Big(\sum_{\xi \in \XX_\dd}  |B(\xi, \lambda^{-1})|| f(\xi)|^p\Big)^{1/p}
\le \| f \|_p
\le 2\Big(\sum_{\xi \in \XX_\dd}|B(\xi, \lambda^{-1})|| f(\xi)|^p\Big)^{1/p}.
\end{equation}
Also, note that estimates $(\ref{samp3})$ are immediate for $p=\infty$ with the usual modification
and hold when $2<p<\infty$ with some modification of the constant in $(\ref{gamma-eps})$ ($ \gamma $ depends on $p$).
We do not elaborate on this since we shall only need $(\ref{samp3})$ for $p=2$.
\end{remark}

\subsection{Cubature formula for \boldmath $\Sigma_\lambda^1$}

In this subsection we utilize the Marcinkiewicz-Zygmund inequality from Proposition~\ref{prop:Marcinkiewicz}
for the construction of a cubature formula on $\Sigma_\lambda^1$. 


\begin{theorem}\label{thm:quadrature}
Let $0< \gamma < 1$ and
\begin{equation}\label{def-gamma}
\CPhi  \gamma^{\alpha }\Cdiam = \frac 14.
\end{equation}
Let $\lambda \ge 1$ and suppose $\XX_\ddelta$ is a maximal $\ddelta$-net on $M$
with $\ddelta :=\gamma \lambda^{-1}$.
Then there exist positive constants $($weights$)$
$\{\ww^\lambda_\xi\}_{\xi \in \XX_\ddelta}$ such that
\begin{equation}\label{quadrature}
\int_M f(x) d\mu(x) = \sum_{\xi \in \XX_\ddelta} \ww^\lambda_\xi f(\xi),
\quad f \in \Sigma^1_\lambda,
\end{equation}
and
\begin{equation}\label{quad-weight}
\small\frac 23|B(\xi, \ddelta/2)|
\le \ww^\lambda_\xi
\le 2|B(\xi, \ddelta)|, \quad \xi \in \XX_\ddelta.
\end{equation}
\end{theorem}

We shall derive this theorem 
from Proposition~\ref{prop:Marcinkiewicz} and a version of the~Hahn-Banach theorem for ordered
linear spaces. We next give a theorem of Bauer of this sort
(adapted to the case of linear normed spaces)
that best serves our purposes and refer the reader to \cite{AL} for its proof.


\begin{theorem}[Bauer]\label{thm:Bauer}
Suppose $E$ is a linear normed space, $F \subset E$ is a subspace of $E$,
and $C$ is a convex cone in $E$, which determines an order on $E$ $(f\le g$ if $g-f\in C)$.
Set $V:=\{f\in E: \|f\|\le 1\}$.
Let $\Lambda: F\to \bR$ be a linear functional on~$F$.
Then $\Lambda$ can be extended to a linear functional $\tilde\Lambda$ on $E$ which is
$(i)$ positive, i.e. $\tilde\Lambda(f)\ge 0$ if $f\in C$, and
$(ii)$ $|\tilde\Lambda(f)|\le \|f\|$ for $f\in E$,
if and only if
\begin{equation}\label{Bauer-con}
\Lambda(f) \ge -1 \quad \hbox{for all $f\in F\cap(V+C)$.}
\end{equation}
A simple rescaling shows that the theorem holds if the condition in
$(ii)$ above is replaced by $\tilde\Lambda(f)\le \Cstar\|f\|$ and the condition in $(\ref{Bauer-con})$
by $\Lambda(f) \ge -\Cstar$, where $\Cstar>0$ is a constant.
\end{theorem}

We next show how the Marcinkiewicz-Zygmund inequality implies the existence of a quadrature rule
in a general setting and then apply the result to our particular case.


\begin{proposition}\label{propo:Marc-qubature}
Suppose $(X, \mu)$ is a measure space and let $\cH$ be a space of $\mu-$integrable functions
defined everywhere on $X$.
Suppose $\{A_i\}_{i \in I} $ is a finite or countable disjoint partition of $X$,
i.e. $X=\cup_{i\in I}A_i$ and $A_i \cap A_j =\emptyset$ if $i\ne j$,
consisting of measurable subsets of $X$ of finite measure $(0<\mu(A_i) <\infty)$.
Let $\xi_i \in A_i$, $i\in I$.
Also, assume that there exists a constant $\alpha <\frac 12$ such that
\begin{equation}\label{cub}
\sum_{i \in I} \int_{A_i} |f(x)- f(\xi_i)| d\mu(x)  \leq \alpha \int_X |f(x)| d\mu(x),
\quad f \in \cH.
\end{equation}
Then there exist positive constants $\{\gamma_i\}_{i\in I}$ such that
\begin{equation}\label{cub2}
\int_X f(x) d\mu(x)  = \sum_{i \in I} \gamma_i f(\xi_i)
\quad \hbox{for}\;\; f \in \cH,
\end{equation}
and
\begin{equation}\label{cub3}
\frac{1-2\alpha}{1-\alpha} \mu(A_i) \le \gamma_i \le \frac{1+2\alpha}{1-\alpha}\mu(A_i),
\quad i\in I.
\end{equation}
\end{proposition}

\noindent
{\bf Proof.}
Consider the discrete positive measure
$d\nu := \sum_{i\in I} \mu(A_i) \delta_{\xi_i}$ on $X$,
supported on the set $\XX:=\{\xi_i: i\in I\}$,
and let $\LL^1(\nu)$ be the respective (weighted discrete) $L^1$-space.
By (\ref{cub}) we obtain for $f\in \cH$
$$
\Big|\int_X  f d\mu - \int_X f d\nu \Big| \leq \alpha \int_X |f| d\mu,
$$
and
\begin{equation}\label{IN}
(1- \alpha) \| f \|_{\LL^1(\mu)} \leq  \| f \|_{\LL^1(\nu)} \leq (1 + \alpha) \| f \|_{\LL^1(\mu)},
\end{equation}
Hence
$$
\int_X  f d\mu - \int_X f d\nu
\geq -  \alpha \int_X |f| d\mu
\geq -\frac \alpha{1-\alpha} \int_X |f| d\nu,
$$
which readily implies
\begin{equation}\label{IN2}
\int_X  f d\mu - \frac{1-2\alpha}{1-\alpha} \int_X f d\nu
\geq  -\frac \alpha{1-\alpha} \int_X (|f|-f) d\nu.
\end{equation}
On the other hand, (\ref{IN}) yields that the operator
$J :  f \in \cH \mapsto \{f(\xi_i)\}_{i\in I} \in  \LL^1(\nu)$
is continuous and, moreover, if $J(\cH) = \tilde{\cH} \subset \LL^1(\nu)$,
then the operator
$$
J^{-1}:  g \in  \tilde{\cH} \mapsto \int_X J^{-1}(g)d\mu
$$
is well-defined and continuous, and by (\ref{IN2})
\begin{equation}\label{est-J-1}
\int_X J^{-1}(g) d\mu - \frac{1-2\alpha}{1-\alpha}  \int_X g \,d\nu
\geq - \frac{2\alpha}{1-\alpha} \int_X (|g|-g) d\nu.
\end{equation}
Let the linear functional $\Lambda: \tilde\cH \mapsto \R$ be defined by
\begin{equation}\label{def-Lam}
\Lambda: g \in  \tilde{\cH} \mapsto \Lambda (g)
:=  \int_X J^{-1}(g) d\mu -   \frac{1-2\alpha}{1-\alpha}  \int_X g d\nu.
\end{equation}
We next apply Theorem~\ref{thm:Bauer} with
$E=\LL^1(\nu)$, $F=\tilde\cH$,
$$\hbox{
$C=\{f\in \LL^1(\nu): f(\xi)\ge 0, \xi\in\XX\}$,
\quad
$V=\{f\in\LL^1(\nu):\|f\|_{\LL^1(\nu)}\le 1\}$,
}
$$
and the linear functional $\Lambda$ from (\ref{def-Lam}).
Evidently, in this case, $f\in F\cap(V+C)$ if and only if
$f \in   \tilde{\cH}$ and $f$ can be represented in the form
$f = g+h$, where $\|g\|_{\LL^1(\nu)} \leq 1$ and $h\ge 0$.
Then by (\ref{est-J-1}) it follows that
\begin{align*}
\Lambda(f) &=\Lambda(g+h)
\ge - \frac{2\alpha}{1-\alpha}\int_X (|g+h|-g-h)d\nu\\
&\ge - \frac{2\alpha}{1-\alpha}\int_X (|g|-g)d\nu
\ge - \frac{4\alpha}{1-\alpha}\int_X |g|d\nu
\geq - \frac{ 4\alpha}{1-\alpha}=:-\Cstar.
\end{align*}
Applying now Theorem~\ref{thm:Bauer} we conclude
that there exists a positive continuous extension $\tilde\Lambda$ of $\Lambda$ to $\LL^1(\nu)$
such that $\|\tilde\Lambda\|\le \Cstar=\frac{4\alpha}{1-\alpha}$.
However, as is well-known  (see e.g. \cite{DS}) $(\LL^1(\nu))^*=\LL^\infty(\nu)$.
Therefore, there exists a sequence
$\beta\in \LL^\infty(\nu)$, $\beta=\{\beta_i\}_{i\in I}$, such that
$\|\beta \|_{\infty} =\sup_{i\in I} \beta_i \le \frac{4\alpha}{1-\alpha}$ and
$$
\tilde\Lambda(f)=\sum_{i\in I} f(\xi_i)\beta_i\mu(A_i),
\quad f\in\LL^1(\nu).
$$
Since $\tilde\Lambda$ is positive, we have $\beta_i \ge 0$, $i\in I$.
Consequently, for any $f\in \cH$
$$
\Lambda(f)=\int_X f d\mu - \frac{1-2\alpha}{1-\alpha} \sum_i \mu(A_i) f(\xi_i)
= \sum_i \beta_i \mu(A_i) f(\xi_i),
$$
where $0 \le \beta_i \le \frac{2 \alpha}{1-\alpha}$, which leads to
$\int_X f(x) d\mu(x) = \sum_{i \in I} \gamma_i f(\xi_i)$ for $f \in \cH$,
where
$$
\frac{1-2\alpha}{1-\alpha} \mu(A_i)
\le \gamma_i
\le \frac{1-2\alpha}{1-\alpha} \mu(A_i ) + \frac{4\alpha}{1-\alpha}\mu(A_i)
=  \frac{1+2\alpha}{1-\alpha} \mu(A_i).
$$
The proof is complete.
$\qed$

\subsubsection*{\bf Proof of theorem \ref{thm:quadrature}.}

Let $\XX_\ddelta$ be a maximal $\ddelta-$net on $M$ with $ \ddelta=\frac \gamma \lambda$.
Then by 
Proposition~\ref{prop:Marcinkiewicz} we have
$$
\sum_{\xi \in \XX_\ddelta} \int _{A_\xi}|f(x) - f(\xi)| d\mu(x)
\le  \CPhi\gamma^{\alpha }\Cdiam\|f\|_{\LL^1}.
$$
If $\gamma>0$ and $\CPhi\gamma^{\alpha }\Cdiam \le \frac 14$,
then Theorem~\ref{thm:quadrature} follows at once from Proposition~\ref{propo:Marc-qubature}.
$\qed$

\section{Construction of frames}\label{sec:frames}
\setcounter{equation}{0}

An important part of our development in this article is the construction of well-localized decomposition systems 
for spaces of functions or distributions in the general setting of this article.
The goal will be to construct a pair of dual frames, where the elements of both frames
are band limited and have nearly exponential space localization.

\subsection{A natural (Littlewood-Paley type) frame for \boldmath $\LL^2$}\label{natural-frame}

We begin with the construction of a well-localized frame based on the kernels of spectral operators
considered in \S\ref{local-kernels}.

Let $\Phi\in C^\infty (\bR_+)$,
$\Phi(u)=1$ for $u\in[0, 1]$,
$0\le \Phi \le 1$, and
$\supp \Phi \subset [0, b]$, where $b>1$ is the constant from Theorem~\ref{thm:norms}.
Set $\Psi(u):=\Phi(u)-\Phi(bu)$ and note that $0\le \Psi \le 1$ and
$\supp \Psi\subset [b^{-1}, b]$.
We shall also assume that $\Phi$ is selected so that
$\Psi(u) \ge c>0$ for $u\in [b^{-3/4}, b^{3/4}]$.
We set
\begin{equation}\label{def-Psi-j}
\Psi_0(u):=\Phi(u) \quad\hbox{ and }\quad \Psi_j(u):=\Psi(b^{-j}u),\;\; j\ge 1.
\end{equation}
Clearly, $\Psi_j\in C^\infty (\bR_+)$, $0\le \Psi_j \le 1$, $\supp \Psi_0 \subset [0, b]$,
$\supp \Psi_j \subset [b^{j-1}, b^{j+1}]$, $j\ge 1$, and
$\sum_{j\ge 0}\Psi_j(u) = 1$ for $u\in\bR_+$.
By Corollary~\ref{cor:Littlewood-Paley} we have the following Littlewood-Paley decomposition
\begin{equation}\label{repres-Psi-j}
f= \sum_{j\ge 0} \Psi_j(\sqrt L)f
\quad\hbox{for $f\in \LL^p$, $\;1\le p\le\infty$. \;$(\LL^\infty:=\UCB)$}
\end{equation}
From above it follows that
\begin{equation}\label{sum-Psi-j}
\frac 12 \le \sum_{j\ge 0} \Psi_j^2(u) \leq 1, \quad u\in \R_+,
\end{equation}
and since
$
\|\Psi_j(\sqrt L) f\|_2^2
= \langle \Psi_j(\sqrt L) f, \Psi_j(\sqrt L) f \rangle
= \langle \Psi_j^2(\sqrt L) f,  f \rangle,
$
we get
$$
\sum_{j \ge 0} \|\Psi_j(\sqrt L) f\|_2^2
= \int_0^\infty  \sum_{j\ge 0} \Psi_j^2(u) d \langle F_u f,f\rangle,
$$
and using (\ref{sum-Psi-j}) we arrive at
\begin{equation}\label{frame1}
\frac 12  \| f\|_2^2 \le \sum_{j \ge 0} \|\Psi_j(\sqrt L) f \|_2^2 \le \| f\|_2^2,
\quad f\in \LL^2.
\end{equation}

Here we introduce a constant $0 <\eps <1$ that is sufficiently small and
will be specified later on in 
(\ref{def-gamma-eps}).
Choose $0< \gamma <1$ so that
\begin{equation}\label{pick-gamma}
\CPhi \gamma^{\alpha }\Cdiam = \eps/3,
\end{equation}
where $\CPhi$ is the constant from (\ref{local-Phi})-(\ref{Lip-Phi}) with $\sig_*:=2d+1$,
and $\Cdiam:=2^{2d+1}$.
For any $j\ge 0$ let $\XX_j \subset M$ be a maximal $\ddj-$net on $M$ $($see Proposition~\ref{prop:delta-net}$)$
with
$\ddj:=\gamma b^{-j-2}$ and
suppose $\{A_\xi^j\}_{\xi\in\XX_j}$ is a companion disjoint partition of $M$
consisting of measurable sets such that
$B(\xi, \ddj/2) \subset  A_\xi^j \subset B(\xi,\ddj)$, $\xi \in \XX_j$,
as in Proposition~\ref{prop:delta-net}.
By Theorem~\ref{thm:sampling} we have
\begin{equation}\label{sampling-L2}
(1-\eps)\| f \|_2^2 \le \sum_{\xi \in \XX_j} |A_\xi^j|| f(\xi)|^2 \le (1+\eps)\| f \|_2^2
\quad\hbox{for}\quad f\in \Sigma_{b^{j+2}}^2.
\end{equation}
By the definition of $\Psi_j$ it follows that
$\Psi_j(\sqrt L) f \in \Sigma^2_{b^{j+1}}$ for $f\in \LL^2$,
and hence (\ref{frame1}) and (\ref{sampling-L2}) imply
\begin{equation}\label{frame2}
\frac 14 \|f\|_2^2 \le
\sum_{j \geq 0} \sum_{\xi \in \XX_j} |A^j_\xi || \Psi_j(\sqrt L) f (\xi) |^2
\leq  2\| f\|_2^2,
\quad f\in L^2.
\end{equation}
Note that
\begin{align*}
\Psi_j(\sqrt L) f (\xi)
&= \int_M f(u) \Psi_j(\sqrt L)(\xi,u)d\mu(u)\\
&= \int_M f(u) \overline{ \Psi_j(\sqrt L)(u,\xi )} d\mu(u)
= \langle f,  \Psi_j(\sqrt L)( . ,\xi  ) \rangle.
\end{align*}
Consider the system $\{\psi_{j \xi}\}$ defined by
\begin{equation}\label{def-frame}
\psi_{j \xi} (x):= |A^{j}_\xi |^{1/2} \Psi_j(\sqrt L)(x,\xi),
\quad \xi \in \XX_j, j\ge 0.
\end{equation}
From the above observation and (\ref{frame2}) it follows that
$\{\psi_{j \xi}: \xi\in \XX_j, j\ge 0\}$
is a frame for $\LL^2$.

We next record the main properties of this system.


\begin{proposition}\label{prop:frame-prop}
$(a)$ {\rm Localization:} For any $\ssigma>0$ there exist a constant $c_\ssigma >0$
such that for any $\xi\in\XX_j$, $j\ge 0$, we have
\begin{equation}\label{prop-psi-1}
|\psi_{j \xi} (x)|
\le c_\ssigma|B(\xi, b^{-j})|^{-1/2}(1+b^j\rho(x, \xi))^{-\ssigma}
\end{equation}
and if $\rho(x, y) \le b^{-j}$
\begin{equation}\label{prop-psi-Lip}
|\psi_{j \xi} (x)- \psi_{j \xi} (y)|
\le c_\ssigma|B(\xi, b^{-j})|^{-1/2}(b^j\rho(x, y))^\alpha (1+b^j\rho(x, \xi))^{-\ssigma},
\quad \alpha>0.
\end{equation}

$(b)$ {\rm Norms:}
\begin{equation}\label{prop-psi-2}
\|\psi_{j \xi}\|_p \sim |B(\xi, b^{-j})|^{\frac 1p-\frac 12},
\quad 0< p \leq \infty.
\end{equation}
The constants involved in the previous equivalence depend of $p.$

$(c)$ {\rm Spectral localization:}
$\psi_{0\xi}\in \Sigma_b^p$ if $\xi\in \XX_0$ and
$\psi_{j\xi}\in \Sigma_{[b^{j-1}, b^{j+1}]}^p$
if $\xi\in \XX_j$, $j\ge 1$, $0<p\le\infty$.

$(d)$
The system $\{\psi_{j \xi}\}$ is a frame for $\LL^2$, namely, 
\begin{equation}\label{frame3}
 4^{-1}\|f\|_2^2 \le
\sum_{j \geq 0} \sum_{\xi \in \XX_j} |\langle f, \psi_{j \xi}\rangle|^2
\leq  2\| f\|_2^2,
\quad \forall f\in \LL^2.
\end{equation}

\end{proposition}

\noindent
{\bf Proof.}
Estimates (\ref{prop-psi-1}) and (\ref{prop-psi-Lip}) follow by Theorem~\ref{thm:main-local-kernels};
(\ref{prop-psi-2}) follows by Theorem~\ref{thm:norms}.
The spectral localization is obvious by the definition.
Estimates (\ref{frame3}) follow by (\ref{frame2}).
$\qed$

\subsection{Dual frame}\label{dual-frame}

Our next (nontrivial) step is to construct a dual frame $\{\tilde\psi_{j \xi}\}$ to $\{\psi_{j \xi}\}$
with elements of similar space and spectral localization.
We begin this construction by introducing two new cut-off functions
by ``stretching"
$\Psi_0$ and $\Psi_1$ from \S\ref{natural-frame}:
\begin{equation}\label{def-LLam}
\LLam_0(u) := \Phi (b^{-1}u)\quad \hbox{and} \quad
\LLam_1(u) := \Phi (b^{-2}u)- \Phi(bu)
= \LLam_0(b^{-1}u) - \LLam_0(b^2u).
\end{equation}
Note that
$\supp\LLam_0 \subset [0, b^2]$, $\LLam_0(u)=1$ for $u\in [0, b]$, 
$\supp\LLam_1 \subset [b^{-1}, b^3]$, $\LLam_1(u)=1$ for $u\in [1, b^2]$,
and $0\le \LLam_0, \LLam_1 \le 1$.
Therefore,
\begin{equation}\label{prop-LLam}
\LLam_0(u)\Psi_0(u) = \Psi_0(u), \quad
\LLam_1(u)\Psi_1(u) = \Psi_1(u).
\end{equation}
We shall also use the 
cut-off function
$\TTheta(u):=\Phi(b^{-3}u)$.
Note that
$\supp \TTheta \subset [0, b^4]$,
$\TTheta(u)=1$ for $u\in[0, b^3]$, and $\TTheta\ge 0$.
Hence,
$\TTheta(u) \LLam_j(u) = \LLam_j(u)$, $j=0,1$.

\smallskip

\noindent
{\bf Parameter \boldmath $\sig$:} The dual frame under construction will depend on a parameter
$\sig > 2d+1$
that can be selected as large as we wish. It will govern the localization
properties of the dual frame elements.

With $\sig > 2d+1$ already selected we next record the localization properties of
the operators generated by the above selected functions.
Let $f=\Gamma_0$ or $f=\Gamma_1$ or $f=\Theta$. Then by Corollary~\ref{cor:LS-local-kernels}
there exists a constant $c_\sig>1$ such that for $\delta >0$ and $0\le m\le \sigma$
we have
\begin{equation}\label{gen-local}
|L^m f(\delta \sqrt{L})(x, y)| \le c_\sig \delta^{-2m}D_{\delta, 2\sig}(x, y).
\end{equation}

\smallskip

We now select the constant $0< \eps <1$ so that
\begin{equation}\label{def-gamma-eps}
\frac 1{2\eps} = c_\sig^3 2^{8\sig+9d+10}.
\end{equation}
Recall that the constant $\gamma$, which depends on $\eps$, was defined in (\ref{pick-gamma})
so that
$
\CPhi  \gamma^{\alpha }\Cdiam = \eps/3.
$

The next lemma will be instrumental in the construction of the dual frame.


\begin{lemma}\label{lem:instrument}

Given $\lambda\ge 1$, let $\XX_\dd$ be a~maximal $\dd-$net on $M$ 
with
$\dd:=\gamma\lambda^{-1}b^{-3}$ and
suppose $\{A_\xi\}_{\xi\in\XX_\dd}$ is a companion disjoint partition of $M$
consisting of measurable sets such that
$B(\xi, \dd/2) \subset  A_\xi \subset B(\xi,\dd)$,  $\xi \in \XX_\dd$,
just as in Proposition~\ref{prop:delta-net}.
Set $\kappa_\xi := \frac{1}{1+\varepsilon}|A_\xi |\sim  |B(\xi, \delta)|.$
Let $\LLam=\LLam_0$ or $\LLam=\Gamma_1$.
Then there exists an operator $\TT_\lambda: \LL^2\to\LL^2$
of the form $\TT_\lambda = \Id + \SSS_\lambda$ such that

$(a)$
$$
\|f \|_2\le \|\TT_\lambda f\|_2 \le \frac 1{1-2\varepsilon}\|f \|_2
\quad \forall f \in \LL^2.
$$

$(b)$
$L^m\SSS_\lambda$ with $0\le m\le \sigma$ is an integral operator with kernel
$L^m\SSS_\lambda(x, y)$ verifying
$$
|L^m\SSS_\lambda(x,y)| \le c\lambda^{2m} D_{\lambda^{-1}, \sss}(x,y),\quad x,y\in M.
$$

$(c)$
$\SSS_\lambda (\LL^2)\subset \Sigma_{\lambda b^2}^2$ if $\LLam = \LLam_0$ and
$\SSS_\lambda (\LL^2)\subset \Sigma_{[\lambda b^{-1}, \lambda b^3]}^2$ if $\LLam = \LLam_1$.

$(d)$ For any $f\in \LL^2$ such that $ \Gamma(\lambda^{-1} \sqrt L)f =f$ we have
\begin{equation}\label{instr-1}
f(x) =  \sum_{\xi \in \XX_\dd}  \kappa_\xi  f(\xi)  \TT_\lambda [\LLam_\lambda(\cdot, \xi)](x),
\quad x\in M,
\end{equation}
where $\LLam_\lambda(\cdot, \cdot)$ is the kernel of the operator
$\LLam_\lambda:=\LLam(\lambda^{-1} \sqrt L)$.
\end{lemma}


\noindent
{\bf Proof.}
By Theorem~\ref{thm:sampling} we have 
$$
(1-\eps)\|f\|^2_2
\le \sum_{\xi \in \XX_\dd} |A_\xi||f(\xi)|^2
\le (1+\eps)\| f \|^2_2
\quad\hbox{for}\; f\in \Sigma_{\lam b^3}^2,
$$
and setting
$\kappa_\xi := \frac{1}{1+\varepsilon}|A_\xi|$
we get
\begin{equation}\label{instr-2}
(1- 2\eps)\|f\|^2_2\le \sum_{\xi \in \XX_\dd}  \kappa_\xi|f(\xi)|^2\le \| f \|^2_2
\quad\hbox{for}\; f\in \Sigma_{\lam b^3}^2.
\end{equation}
Denote briefly $\Theta_\lambda:=\Theta(\lambda^{-1} \sqrt L)$ and
let $\Theta_\lambda(\cdot, \cdot)$ be the kernel of this operator.
Consider now the positive self-adjoint operator $U_\lambda$ with kernel
$$
U_\lambda(x,y)
= \sum_{\xi \in \XX_\dd} \kappa_\xi  \Theta_\lambda(x,\xi)\Theta_\lambda(\xi,y).
$$
By (\ref{gen-local})
$|\Theta_\lambda (x,y)|  \le c_\sig D_{\lambda^{-1}, 2\sss}(x,y)$ for $x, y\in M$.
Therefore, taking into account that $\delta=\gamma\lambda^{-1}b^{-3} < \lambda^{-1}$
and $2\sss > 2d+1$ we can apply (\ref{discr-comp}) to obtain
\begin{equation}\label{local-U}
|U_\lambda (x,y)|  \le c_\sig c_\sharp D_{\lambda^{-1}, 2\sss}(x,y),
\quad c_\sharp:=2^{2\sig+3d+3}
\end{equation}
Also, if $f\in\Sigma_{\lambda b^3}^2$, then
$
\langle U_\lambda f, f  \rangle =   \sum_{\xi \in \XX_\dd}  \kappa_\xi \  | f(\xi)|^2
$
and hence, using (\ref{instr-2}),
\begin{equation}\label{R1}
(1-2\varepsilon) \| f \|_2^2 \leq\langle U_\lambda f, f  \rangle  \leq \|f \|_2^2
\quad \hbox{for } f\in\Sigma_{\lambda b^3}^2.
\end{equation}
Denote briefly $\LLam_\lambda:= \LLam(\lambda^{-1} \sqrt L)$
and let $\LLam_\lambda(\cdot, \cdot)$ be the kernel of this operator
(recall that $\LLam = \LLam_0$ or $\LLam = \LLam_1$).
We define yet another self-adjoint kernel operator by
$$
R_\lambda := \LLam_\lambda(\Id - U_\lambda) \LLam_\lambda
= \LLam_\lambda^2 - \LLam_\lambda U_\lambda \LLam_\lambda.
$$
Set $V_\lambda:= \LLam_\lambda U_\lambda \LLam_\lambda$
and denote by $V_\lambda (\cdot, \cdot)$ its kernel.
Since
$ \Theta(u)\LLam(u) = \LLam(u)$, we have
\begin{align*}
V_\lambda (x, y)
&= \sum_{\xi \in \XX_\dd}  \kappa_\xi
\int_M \int_M \LLam_\lambda(x,u)\Theta_\lambda(u,\xi)\Theta_\lambda(\xi,v)\LLam_\lambda(v,y) du dv\\
&= \sum_{\xi \in \XX_\dd}  \kappa_\xi \LLam_\lambda(x,\xi)\LLam_\lambda(\xi,y).
\end{align*}
Now, by (\ref{gen-local}), (\ref{Comp}), and (\ref{local-U})
it follows that for $0\le m\le \sigma$
\begin{align*}
|L^mR_\lambda(x,y)|
&\le c_\sig^2 c_\star \lambda^{2m}D_{\lambda^{-1}, 2\sss}(x,y)
+
c_\sig^3 c_\sharp c_\star^2 \lambda^{2m}D_{\lambda^{-1}, 2\sss}(x,y)\\
&\le 2c_\sig^3 c_\sharp c_\star^2 \lambda^{2m}D_{\lambda^{-1}, 2\sss}(x,y).
\end{align*}
Here $c_\star:=2^{2\sig+2d+2}$ and as above $c_\sharp:=2^{2\sig+3d+3}$.
To simplify our notation we set
$C_\sig:=2c_\sig^3 c_\sharp c_\star^2= c_\sigma^32^{6\sig+7d+8}$.
Thus we have
\begin{equation}\label{est-Rxy}
|L^mR_\lambda(x,y)|
\le C_\sig \lambda^{2m}D_{\lambda^{-1}, 2\sss}(x,y), \quad 0\le m\le \sigma.
\end{equation}

By the definition of $R_\lambda$ we have
$$
\langle R_\lambda f, f  \rangle
= \|\LLam_\lambda f\|_2^2 -
\langle U_\lambda \LLam_\lambda f , \LLam_\lambda f \rangle
\quad \hbox{for $f\in\LL^2$.}
 $$
Since $ \LLam_\lambda (\LL^2) \subset \Sigma^2_{\lambda b^3}$, then
$\Theta_\lambda \LLam_\lambda f = \LLam_\lambda f$, and by (\ref{R1})
$$
(1-2\varepsilon) \| \LLam_\lambda f \|^2_2
\le \langle U_\lambda \LLam_\lambda f , \LLam_\lambda f \rangle
\le  \| \LLam_\lambda f \|^2_2,
\quad f\in\LL^2.
$$
Therefore,
$$
0 \leq  \langle R_\lambda f, f  \rangle
\le 2\eps \| \LLam_\lambda f \|^2_2
\le 2\eps \| f \|_2^2,
\quad f\in\LL^2,
$$
where for the last inequality we used that $ \|\LLam\|_\infty \le 1$.
Consequently,
$$
\| R_\lambda\|_{2\rightarrow 2} \leq 2\eps <1
\quad \hbox{and}\quad
(1-2\eps)\|f\|_2 \le  \|( \Id-R_\lambda) f\|_2 \le \| f \|_2,
\quad f\in\LL^2.
$$

We now define
$
\TT_\lambda :=(\Id- R_\lambda)^{-1} = \Id + \sum_{k\ge 1} R_\lambda^k =: \Id + \SSS_\lambda.
$
Clearly
\begin{equation}\label{B0}
\|f\|_2 \le  \|\TT_\lambda f \|_2 \le  \frac 1{1-2\eps}\|f \|_2
\quad  \forall f \in \LL^2.
\end{equation}
If  $\LLam_\lambda f = f$, then
$$
f= \TT_\lambda (f -R_\lambda f)
= \TT_\lambda \big( f- \LLam_\lambda f + V_\lambda f \big)
=  \TT_\lambda V_\lambda f.
$$
On the other hand, if $\LLam_\lambda f =f$, then
$
(V_\lambda f)(x) = \sum_{\xi \in \XX_\dd}  \kappa_\xi  f(\xi)   \LLam_\lambda(x,\xi)
$
and hence 
\begin{equation}\label{B1}
f(x) =  \sum_{\xi \in \XX_\dd}  \kappa_\xi  f(\xi)  \TT_\lambda[\LLam_\lambda(\cdot,\xi)](x).
\end{equation}
Note that by construction 
\begin{equation}\label{B2}
\SSS_\lambda: \LL^2 \mapsto \Sigma^2_{\lambda b^3}
\quad\hbox{if $\LLam=\LLam_0$ and} \quad
\SSS_\lambda: \LL^2 \mapsto \Sigma^2_{[\lambda b^{-1}, \lambda b^3]}
\quad \hbox{if $\LLam=\LLam_1$.}
\end{equation}


It remains to establish the space localization of the kernel $L^m\SSS_\lambda (x, y)$
of the operator~$L^m\SSS_\lambda$.
Our method borrows from \cite{L}.
Consider first the case $m=0$.
Denoting by $R_\lambda^k(x,y)$ the kernel of $R_\lambda^k$,
we have
$$
|\SSS_\lambda (x,y) |  \le  \sum_{k\ge 1}  |R_\lambda^k(x,y)|.
$$
But since $R_\lambda^k = \Theta_\lambda R_\lambda^k \Theta_\lambda$
we get by (\ref{gen-local}) with $f=\Theta$ and the fact that
$\|R_\lambda\|_{2\rightarrow 2} \le 2\eps$,
applying Proposition~\ref{prop:prod-oper},
\begin{equation}\label{kernel-Rk-1}
|R_\lambda^k(x,y)|
\le  \frac{c_dc_\sig^2 \|R_\lambda\|^k_{2\rightarrow 2}}
{\sqrt{|B(x, \lambda^{-1})| |B(y, \lambda^{-1})|}}
\le  \frac{(2\eps)^k c_dc_\sig^2}{\sqrt{|B(x, \lambda^{-1})|  |B(y, \lambda^{-1})|}},
\end{equation}
where $c_d:=2^{4d+4}$.
On the other hand, applying repeatedly estimate (\ref{Comp}) $k$ times using (\ref{est-Rxy})
with $m=0$ we obtain
\begin{equation}\label{kernel-Rk-2}
|R_\lambda^k(x,y)|\le  C_\sig^k c_\star^{k-1} D_{\lambda^{-1}, 2\sss}(x,y),
\quad c_\star=2^{2\sig+2d+2}.
\end{equation}
Therefore, for any $K\in \bN$
\begin{align*}
|\SSS_\lambda (x,y)|
&\le  \sum_{k=1}^K C_\sig^k c_\star^{k-1} D_{\lambda^{-1}, 2\sss}(x,y)
+ \sum_{k>K} \frac{(2\eps)^k c_dc_\sig^2}{\sqrt{|B(x, \lambda^{-1})|  |B(y, \lambda^{-1})|}}\\
&\le \frac{C_\sig}{\sqrt{|B(x, \lambda^{-1})||B(y, \lambda^{-1})|}}
\Big\{ \frac 1{(1+ \lambda \rho(x,y))^{2\sss}}
\sum_{k=1}^K   (c_\star C_\sig)^{k-1} + \frac{(2\eps)^{K+1}}{1-2\eps} \Big\}\\
&\le \frac{C_\sig}{\sqrt{|B(x, \lambda^{-1})|  |B(y, \lambda^{-1})|}}
\Big\{ \frac 1{(1+ \lambda \rho(x,y))^{2\sss}}
\frac{ (c_\star C_\sig)^{K}}{c_\star C_\sig-1}
+ \frac{(2\eps)^{K+1}}{1-2\eps} \Big\}\\
&\le  \frac{2C_\sig}{\sqrt{|B(x, \lambda^{-1})|
|B(y, \lambda^{-1})|}}
\Big\{ \frac{(c_\star C_\sig)^{K-1}}{(1+ \lambda \rho(x,y))^{2\sss}} +(2\eps)^{K+1}\Big\}.
\end{align*}
Choose $K\ge 1$ so that
$(\frac{1}{2\eps})^{K-1} \le (1+ \lambda \rho(x,y))^\sss <(\frac{1}{2\eps})^{K}$
and note that $\frac 1{2\eps} = c_\star C_\sig$
by (\ref{def-gamma-eps}).
Then from above we get
\begin{equation}\label{kernel-S}
|\SSS_\lambda (x,y)|
\le  \frac{4C_\sig}{\sqrt{|B(x, \lambda^{-1})|  |B(y, \lambda^{-1})|}}
\frac 1{ (1+ \lambda \rho(x,y))^\sss}
= 4C_\sig D_{\lambda^{-1}, \sss} (x,y).
\end{equation}


Let $1\le m\le \sigma$.
Since $L^m R_\lambda^k = L^m\Theta_\lambda R_\lambda^k \Theta_\lambda$, with slight modification
of the argument above, (\ref{gen-local}) implies that (\ref{kernel-Rk-1})
holds for the kernel $L^mR^k(\cdot, \cdot)$ with an~additional factor $\lambda^{2m}$ to the right.
On the other hand (\ref{est-Rxy}) implies that estimate (\ref{kernel-Rk-2}) also holds for
$L^mR^k(\cdot, \cdot)$ with an additional factor $\lambda^{2m}$ to the right.
Then proceeding exactly as above it follows that estimate (\ref{kernel-S}) holds for
$L^mS(\cdot, \cdot)$ with an additional factor $\lambda^{2m}$ to the right.
This completes the proof of the lemma.
$\qed$

\smallskip

Armed with this lemma we can now complete the construction of the dual frame.
We shall utilize the functions and operators introduced in \S\ref{natural-frame} and above.

Denote briefly
$\LLam_{\lam_0}:= \LLam_0(\sqrt L)$ and
$\LLam_{\lam_j}:= \LLam_1(b^{-j+1}\sqrt L)$ for $j\ge 1$,
$\lam_j:=b^{-j+1}$.
Observe that since $\LLam_0(u)=1$ for $u \in [0, b]$ and $\LLam_1(u)=1$ for $u\in [1, b^2]$,
then $\LLam_{\lam_0}(\Sigma_b^2)=\Sigma_b^2$ and
$\LLam_{\lam_j}(\Sigma_{[b^{j-1}, b^{j+1}]}^2) = \Sigma_{[b^{j-1}, b^{j+1}]}^2$, $j\ge 1$.
On the other hand, it is readily seen that
$\Psi_0(\cdot, y) \in \Sigma_{b}^2$ and
$\Psi_j(\cdot, y) \in \Sigma_{[b^{j-1}, b^{j+1}]}^2$ if $j\ge 1$.
Therefore, we can apply Lemma~\ref{lem:instrument} with
$\XX_j$ and $\{A_\xi^j\}_{\xi\in\XX_j}$ from \S\ref{natural-frame}, and
$\lambda=\lambda_j=b^{j-1}$
to obtain
\begin{equation}\label{rep-Psi-j}
\Psi_j(\sqrt L)(x, y)=
\sum_{\xi \in \XX_j}  \kappa_\xi  \Psi_j(\xi, y) \TT_{\lambda_j} [\LLam_{\lambda_j}(\cdot, \xi)](x),
\quad \kappa_\xi=(1+\eps)^{-1}|A_\xi|.
\end{equation}
By (\ref{def-frame}) we have
$\psi_{j\xi}(x)= |A_\xi^j|^{1/2}\Psi_j(\xi, x)$ and we now set
\begin{equation}\label{rep-psi-tpsi}
\tilde\psi_{j\xi}(x):=c_\eps |A_\xi^j|^{1/2}\TT_{\lambda_j} [\LLam_{\lambda_j}(\cdot, \xi)](x),
\quad \xi\in\XX_j, \quad c_\eps:=(1+\eps)^{-1}.
\end{equation}
Thus $\{\tilde\psi_{j\xi}: \xi\in\XX_j, j\ge 0\}$ is the desired dual frame.
Observe immediately that (\ref{rep-Psi-j}) takes the form
\begin{equation}\label{rep2-Psi-j}
\Psi_j(\sqrt L)(x, y)=
\sum_{\xi \in \XX_j} \psi_{j\xi}(y)\tilde\psi_{j\xi}(x).
\end{equation}

We next record the main properties of the dual frame $\{\tilde\psi_{j\xi}\}$.


\begin{theorem}\label{thm:dual-frame}
$(a)$ {\rm Representation:}
For any $f\in\LL^p$, $1\le p\le \infty$, 
we have
\begin{equation}\label{rep-L2}
f = \sum_{j\ge 0}\sum_{\xi \in \XX_j} \langle f, \tilde\psi_{j\xi}\rangle \psi_{j\xi}
= \sum_{j\ge 0}\sum_{\xi \in \XX_j} \langle f, \psi_{j\xi}\rangle \tilde\psi_{j\xi}
\quad\hbox{in}\;\; \LL^p.
\end{equation}

$(b)$ {\rm Frame:}
The system $\{\tilde\psi_{j \xi}\}$ as well as $\{\psi_{j \xi}\}$
is a frame for $\LL^2$, namely, there exists a constant $c>0$ such that
\begin{equation}\label{frame-tpsi}
 c^{-1}\|f\|_2^2 \le
\sum_{j \geq 0} \sum_{\xi \in \XX_j} |\langle f, \tilde\psi_{j \xi}\rangle|^2
\leq  c\| f\|_2^2,
\quad \forall f\in \LL^2.
\end{equation}

$(c)$ {\rm Space localization:} For any $\xi\in\XX_j$, $j\ge 0$, and $0\le m\le \sigma$
\begin{equation}\label{prop-tpsi-1}
|L^m\tilde\psi_{j \xi} (x)|
\le c_\ssigma b^{2jm}|B(\xi, b^{-j})|^{-1/2}(1+b^j\rho(x, \xi))^{-\ssigma},
\end{equation}
and if $\rho(x, y) \le b^{-j}$
\begin{equation}\label{prop-tpsi-Lip}
|\tilde\psi_{j \xi} (x)- \tilde\psi_{j \xi} (y)|
\le c_\ssigma|B(\xi, b^{-j})|^{-1/2}(b^j\rho(x, y))^\alpha (1+b^j\rho(x, \xi))^{-\ssigma}.
\end{equation}
Here $\ssigma>2d+1$ is the parameter of the dual frame selected 
in the beginning of \S\ref{dual-frame}.

$(d)$ {\rm Spectral localization:}
$\tilde\psi_{0\xi}\in \Sigma_b^p$ if $\xi\in \XX_0$ and
$\tilde\psi_{j\xi}\in \Sigma_{[b^{j-2}, b^{j+2}]}^p$ if $\xi\in \XX_j$, $j\ge 1$,
$d/\sig<p\le\infty$.

$(e)$ {\rm Norms:}
\begin{equation}\label{prop-tpsi-2}
\|\tilde\psi_{j \xi}\|_p \sim |B(\xi, b^{-j})|^{\frac 1p-\frac 12}
\quad\hbox{for} \;\; d/\sig< p \le \infty.
\end{equation}
\end{theorem}


\noindent
{\bf Proof.}
By the definition of $\tilde\psi_{j\xi}$ in (\ref{rep-psi-tpsi}) and Lemma~\ref{lem:instrument}
we have
$$ 
\tilde\psi_{j\xi}(x)
:=c_\eps |A_\xi^j|^{1/2}\TT_{\lambda_j} [\LLam_{\lambda_j}(\cdot, \xi)](x)
= c_\eps |A_\xi^j|^{1/2}\big[\LLam_{\lambda_j}(x, \xi)
+\SSS_{\lambda_j} [\LLam_{\lambda_j}(\cdot, \xi)](x)\big].
$$ 
Then estimate (\ref{prop-tpsi-1}) follows from
the localization of $L^m\LLam_{\lam_j}(\cdot, \cdot)$ given by (\ref{gen-local}),
Lemma~\ref{lem:instrument} (b), and (\ref{Comp}).
Estimate (\ref{prop-tpsi-Lip}) follows by the fact $\LLam_{\lam_j}(\cdot, \cdot)$ is Lip $\alpha$,
given by Theorem~\ref{thm:main-local-kernels}, and
the localization of $S_{\lam_j}(\cdot, \cdot)$, given in Lemma~\ref{lem:instrument}~(b),
exactly as in the proof of Theorem~\ref{thm:local-kernels}.

To establish representation (\ref{rep-L2}) we note that
(\ref{rep2-Psi-j}), (\ref{prop-tpsi-1}), and (\ref{discr-comp}) readily imply
$
\sum_{\xi \in \XX_j} |\psi_{j\xi}(y)||\tilde\psi_{j\xi}(x)|
\le c D_{b^{-j}, \sig-d}(x, y).
$
Then (\ref{rep-L2}) follows by (\ref{repres-Psi-j}) and (\ref{rep2-Psi-j}).


The estimate
\begin{equation}\label{norm-tilde-psi-1}
\|\tilde{\psi}_{j\xi}\|_p \le c|B(\xi, b^{-j})|^{\frac 1p - \frac 12}
\quad\hbox{ for } \;\; d/\sig<p\le \infty
\end{equation}
follows by (\ref{prop-tpsi-1}) and (\ref{tech1}).
On the other hand, Lemma~\ref{lem:instrument} (a) and Theorem~\ref{thm:norms} yield
$$
\|\tilde\psi_{j\xi}\|_2
\ge c |B(\xi, b^{-j})| \|\LLam_{\lambda_j}(\cdot, \xi)\|_2
\ge c'>0.
$$
From this and (\ref{norm-tilde-psi-1}) one easily derives
$\|\tilde{ \psi}_{j\xi}\|_p
\ge  c|B(\xi, b^{-j})|^{\frac 1p - \frac 12}$
for $0<p\le \infty$
(see the proof of Theorem~\ref{thm:norms}).
%


For the proof of (\ref{frame-tpsi}) we shall employ the following lemma which will be instrumental
in the development of Besov spaces later on as well.


\begin{lemma}\label{lem:technical}
$(a)$ For any $f\in\LL^p$, $1\le p\le\infty$,
\begin{equation}\label{technical-1}
\Big(\sum_{\xi\in\XX_j}\|\langle f, \tilde\psi_{j\xi} \rangle\psi_{j\xi}\|_p^p\Big)^{1/p} \le c\|f\|_p,
\quad\forall\, j\ge 0.
\end{equation}
$(b)$ For any sequence of complex numbers $\{a_\xi\}_{\xi\in\XX_j}$, $j\ge 0$, and $1\le p\le\infty$,
\begin{equation}\label{technical-2}
\big\|\sum_{\xi\in\XX_j} a_\xi\psi_{j\xi}\big\|_p
\le c\Big(\sum_{\xi\in\XX_j} \|a_\xi\psi_{j\xi}\|_p^p\Big)^{1/p}.
\end{equation}
%
Above each of the $\ell^p$-norms is replaced by the $\sup$-norm when $p=\infty$.
Also $(a)$ and $(b)$ hold with the roles of $\{\psi_{j\xi}\}$ and $\{\tilde\psi_{j\xi}\}$ interchanged.
The constant $c>0$ is independent of $f$, $\{a_\xi\}$, and $j$.
\end{lemma}

\noindent
{\bf Proof.} We shall need the following simple inequalities
\begin{equation}\label{technical-3}
\sum_{\xi\in\XX_j}|\tilde\psi_{j\xi}(x)|\|\psi_{j\xi}\|_1  \le c
\quad\hbox{and}\quad
\sum_{\xi\in\XX_j} |\psi_{j\xi}(x)|\|\psi_{j\xi}\|_1 \le c,
\quad x\in M,
\end{equation}
where the roles of $\{\psi_{j\xi}\}$ and $\{\tilde\psi_{j\xi}\}$ can be switched.
Using (\ref{prop-tpsi-1}) with $m=0$ and (\ref{prop-psi-2}) we obtain
$$
\sum_{\xi\in\XX_j}|\tilde\psi_{j\xi}(x)|\|\psi_{j\xi}\|_1
\le c\sum_{\xi\in\XX_j}(1+b^j\rho(x, \xi))^{-\sigma} \le c<\infty,
$$
where for the last inequality we used (\ref{discr-tech11}) and the fact that $\sigma \ge 2d+1$.
This gives the left-hand side inequality in (\ref{technical-3}).
The proof of the other inequality is the same.


Estimate (\ref{technical-1}) is immediate from (\ref{technical-3}) when $p=1$.
In the case $p=\infty$ (\ref{technical-1}) follows readily by the inequality
$\|\tilde\psi_{j\xi}\|_1\|\psi_{j\xi}\|_\infty \le c<\infty$ which is a consequence of
(\ref{prop-psi-2}) and (\ref{prop-tpsi-2}).

To prove (\ref{technical-1}) in the case $1<p<\infty$ we just apply
H\"{o}lder's inequality ($1/p+1/p'=1$) and obtain
\begin{align*}
\|\langle f, \tilde\psi_{j\xi} \rangle\psi_{j\xi}\|_p^p
&\le \Big(\int_M|f(x)|\tilde\psi_{j\xi}(x)|d\mu(x)\Big)^p\|\psi_{j\xi}\|_p^p\\
&= \Big(\int_M|f(x)|\tilde\psi_{j\xi}(x)|^{1/p}|\tilde\psi_{j\xi}(x)|^{1/p'}d\mu(x)\Big)^p\|\psi_{j\xi}\|_p^p\\
& \le \int_M|f(x)|^p|\tilde\psi_{j\xi}(x)|d\mu(x) \|\tilde\psi_{j\xi}\|_1^{p-1}\|\psi_{j\xi}\|_p^p.
\end{align*}
This coupled with the obvious inequality
$$
\|\tilde\psi_{j\xi}\|_1^{p-1}\|\psi_{j\xi}\|_p^p
\le (\|\tilde\psi_{j\xi}\|_1 \|\psi_{j\xi}\|_\infty)^{p-1}\|\psi_{j\xi}\|_1
\le c\|\psi_{j\xi}\|_1,
$$
using
$\|\tilde\psi_{j\xi}\|_1\|\psi_{j\xi}\|_\infty \le c<\infty$ as above, 
leads to
$$
\sum_{\xi\in\XX_j}\|\langle f, \tilde\psi_{j\xi} \rangle\psi_{j\xi}\|_p^p
\le c\int_M|f(x)|^p \sum_{\xi\in\XX_j}|\tilde\psi_{j\xi}(x)|\|\psi_{j\xi}\|_1 d\mu(x)
\le c\|f\|_p^p.
$$
Here we used (\ref{technical-3}). This confirms the validity of (\ref{technical-1}).


We now turn to the proof of (\ref{technical-2}).
This inequality is obvious when $p=1$.
In the case $p=\infty$ inequality (\ref{technical-2}) follow easily from the right-hand side
inequality in (\ref{technical-3}) and the fact that
$\|\psi_{j\xi}\|_1\|\psi_{j\xi}\|_\infty \le c<\infty$, see (\ref{prop-psi-2}).

To prove (\ref{technical-2}) in the case $1<p<\infty$ we apply the discrete H\"{o}lder inequality
and the right-hand side inequality in (\ref{technical-3}) to obtain
\begin{align*}
\big|\sum_{\xi\in\XX_j} a_\xi\psi_{j\xi}(x)\big|^p
&\le \Big[\sum_{\xi\in\XX_j} |a_\xi|\|\psi_{j\xi}\|_1^{-1}
\big(|\psi_{j\xi}(x)|\|\psi_{j\xi}\|_1\big)^{1/p} \big(|\psi_{j\xi}(x)|\|\psi_{j\xi}\|_1\big)^{1/p'}\Big]^p\\
&\le \sum_{\xi\in\XX_j} |a_\xi|^p\|\psi_{j\xi}\|_1^{-p} |\psi_{j\xi}(x)|\|\psi_{j\xi}\|_1
\Big(\sum_{\xi\in\XX_j}|\psi_{j\xi}(x)|\|\psi_{j\xi}\|_1\Big)^{p-1}\\
&\le c\sum_{\xi\in\XX_j} |a_\xi|^p\|\psi_\xi\|_1^{1-p} |\psi_{j\xi}(x)|.
\end{align*}
Integrating both sides we get
$$
\big\|\sum_{\xi\in\XX_j} a_\xi\psi_{j\xi}\big\|_p^p
\le c\sum_{\xi\in\XX_j} |a_\xi|^p\|\psi_{j\xi}\|_1^{2-p}
\le c\sum_{\xi\in\XX_j} |a_\xi|^p\|\psi_{j\xi}\|_p^p.
$$
Here we used that
$\|\psi_{j\xi}\|_1^{2-p} \sim \|\psi_{j\xi}\|_p^p$, which follows by (\ref{prop-psi-2}).
The proof of Lemma~\ref{lem:technical} is complete.
$\qed$

\smallskip


We are now in a position to complete the proof of Theorem~\ref{thm:dual-frame}.
From (\ref{frame1}) applying (\ref{technical-2}) we get
\begin{align*}
\|f\|_2^2 \le 2\sum_{j\ge 0} \|\Psi_j(\sqrt L)f\|_2^2
\le 2\sum_{j\ge 0}\Big\|\sum_{\xi\in\XX_j}\langle f, \tilde\psi_{j\xi} \rangle \psi_{j\xi}\Big\|_2^2
\le c \sum_{j\ge 0}\sum_{\xi\in\XX_j}|\langle f, \tilde\psi_{j\xi} \rangle|^2,
\end{align*}
which confirms the left-hand side inequality in (\ref{frame-tpsi}).
For the other direction, we first note that since $\supp \Psi_j \subset [b^{j-1}, b^{j+1}]$ and
$\tilde\psi_{j\xi} \in \Sigma_{[b^{j-2}, b^{j+2}]}$ we have by (\ref{repres-Psi-j})
$
\langle f, \tilde\psi_{j\xi}\rangle
= \sum_{\nu=j-2}^{j+2} \big\langle \Psi_\nu(\sqrt{L})f, \tilde\psi_{j\xi}\big\rangle
$
(here $\Psi_\nu:=0$ if $\nu<0$)
and hence
\begin{align*}
\sum_{\xi\in\XX_j}|\langle f, \tilde\psi_{j\xi}\rangle|^2
\le 5\sum_{\nu=j-2}^{j+2} \sum_{\xi\in\XX_j} |\langle \Psi_\nu(\sqrt{L})f, \tilde\psi_{j\xi}\rangle|^2
\le c\sum_{\nu=j-2}^{j+2}\|\Psi_\nu(\sqrt{L})f\|_2^2.
\end{align*}
Here we used (\ref{technical-1}).
Summing up the above inequalities and using (\ref{frame1}) we obtain the right-hand side inequality
in (\ref{frame-tpsi}).
This completes the proof of Theorem~\ref{thm:dual-frame}.
$\qed$

\subsection{Frames in the case when \boldmath $\{\Sigma_\lam^2\}$  possess the polynomial property}
\label{sec:polyn-prop}

The~construction of frames with the desired excellent space and spectral localization is
simple and elegant in the case when the spectral spaces $\Sigma_\lambda^2$
have the polynomial property in the sense of the following


\begin{definition}\label{def:product-prop}

Let $\{F_\lambda, \lambda\ge 0\}$ be the spectral resolution associated with the operator $\sqrt{L}$;
then
$\sqrt L = \int_0^\infty \lambda dF_\lambda$.
We say that the associated spectral spaces
$$
\Sigma_\lambda^2= \{ f \in \bL^2: F_\lambda f =f\}
$$
have the polynomial property if there exists a constant $\kappa>1$ such that
\begin{equation}\label{polyp}
\Sigma_\lambda^2\cdot \Sigma_\lambda^2 \subset \Sigma^1_{\kappa\lambda},
\quad \hbox{i.e.}\quad f, g\in \Sigma_\lambda^2 \Longrightarrow fg\in \Sigma_{\kappa\lambda}^1.
\end{equation}
\end{definition}

The construction begins with two pairs of cut-off functions
$\Psi_0, \Psi, \tilde\Psi_0, \tilde\Psi \in C^\infty(\bR_+)$ with the following properties:
\begin{align*}
&\supp \Psi_0, \tilde\Psi_0\subset [0, b],
\quad \supp \Psi, \tilde\Psi\subset [b^{-1}, b],
\quad
0\le \Psi_0, \Psi, \tilde\Psi_0, \tilde\Psi\le 1,\\
&\Psi_0(u), \tilde\Psi_0(u) \ge c>0,\;\; u\in [0, b^{3/4}], \quad
\Psi(u), \tilde\Psi(u) \ge c>0,\;\; u\in [b^{-3/4}, b^{3/4}],\\
& \Psi_0(u)=1 \;\; \hbox{and} \;\; \tilde\Psi_0(u)=1, \quad u\in [0,1], \quad\hbox{and}\\
&\Psi_0(u)\tilde\Psi_0(u)+\sum_{j\ge 1} \Psi(b^{-j}u)\tilde\Psi(b^{-j}u)=1, \quad u\in \bR_+.
\end{align*}
As in \S\ref{natural-frame}, here $b>1$ is the constant from Theorem~\ref{thm:norms}.
The construction of functions with these properties is quite simple and well-known
and will be omitted.
It is worth pointing out that given $\Psi_0, \Psi$,
then $\tilde\Psi_0, \tilde\Psi$ can be easily constructed with
the above properties (see e.g. \cite{F-J-W}, Lemma (6.9)).

Denote $\Psi_j(u):=\Psi(b^{-j}u)$ and $\tilde\Psi_j(u):=\tilde\Psi(b^{-j}u)$.
Then from above we have
\begin{equation}\label{decomp-unity-2}
\sum_{j\ge 0} \Psi_j(u)\tilde\Psi_j(u)=1, \quad u\in \bR_+.
\end{equation}
This and Proposition~\ref{prop:app-identity} imply the following Calder\'{o}n type decomposition
\begin{equation}\label{calderon}
f=\sum_{j\ge 0} \Psi_j(\sqrt{L})\tilde\Psi_j(\sqrt{L})f, \quad f\in \LL^p, \; 1\le p \le\infty.
\end{equation}

The key idea is that the polynomial property (\ref{polyp}) of the spectral spaces can be used to
discretize the above expansion and as a result to obtain the desired frames.
Indeed, observe first that
$\supp \Psi_0, \tilde\Psi_0\subset [0, b]$ and
$\supp \Psi_j, \tilde\Psi_j\subset [b^{j-1}, b^{j+1}]$, $j\ge 1$.
From this and above it follows that
$\Psi_j(\sqrt{L})$, $\tilde\Psi_j(\sqrt{L})$ are kernel operator whose kernels have
nearly exponential localization and
$\Psi_j(\sqrt{L})(x, \cdot)\in \Sigma_{b^{j+1}}$ and
$\tilde\Psi_j(\sqrt{L})(\cdot, y)\in \Sigma_{b^{j+1}}$.
We now invoke the cubature formula from Theorem~\ref{thm:quadrature}.
With $0<\gamma<1$ the constant from (\ref{def-gamma})
and $\kappa>1$ from (\ref{polyp}), we select a maximal $\delta$-net, say $\XX_j$, on $M$
with $\delta:=\gamma\kappa^{-1}b^{-j-1}\sim b^{-j}$.
Theorem~\ref{thm:quadrature} provides a~cubature formula of the form
$$
\int_M f(x) d\mu(x) = \sum_{\xi\in\XX_j} \ww_{j\xi} f(\xi)
\quad\hbox{for}\quad f\in \Sigma_{\kappa b^{j+1}}^1,
$$
where $\frac 23|B(\xi, \delta/2)| \le \ww_{j\xi} \le 2|B(\xi, \delta)|$.
Since $\Psi_j(\sqrt{L})(x,\cdot)\tilde\Psi_j(\sqrt{L})(\cdot,y) \in \Sigma_{\kappa b^{j+1}}^1$
due to (\ref{polyp}), we get
\begin{align}\label{Psi-Psi}
\Psi_j(\sqrt{L})\tilde\Psi_j(\sqrt{L})(x,y)
&= \int_M \Psi_j(\sqrt{L})(x,u)\tilde\Psi_j(\sqrt{L})(u,y) d\mu(u)\\
&=\sum_{\xi\in\XX_j} \ww_{j\xi} \Psi_j(\sqrt{L})(x,\xi)\tilde\Psi_j(\sqrt{L})(\xi,y).\notag
\end{align}
We now define the frame elements by
\begin{equation}\label{def-psi-xi}
\psi_{j\xi}(x):=\sqrt{\ww_{j\xi}}\Psi_j(\sqrt{L})(x,\xi),
\;\;
\tilde\psi_{j\xi}(x):=\sqrt{\ww_{j\xi}}\tilde\Psi_j(\sqrt{L})(x,\xi),
\;\; \xi\in\XX_j, \; j\ge 0.
\end{equation}

We next present the main properties of the system
$\{\psi_{j\xi}\}$, $\{\tilde\psi_{j\xi}\}$.


\begin{proposition}\label{prop:frame-2}
$(a)$ {\rm Frame property:}
For any $f\in\LL^p$, $1\le p\le\infty$, $(\LL^\infty:=\UCB)$ we have
\begin{equation}\label{frame2-rep}
f = \sum_{j\ge 0}\sum_{\xi \in \XX_j} \langle f, \tilde\psi_{j\xi}\rangle \psi_{j\xi}
= \sum_{j\ge 0}\sum_{\xi \in \XX_j} \langle f, \psi_{j\xi}\rangle \tilde\psi_{j\xi}
\quad\hbox{in}\quad \LL^p
\end{equation}
and
\begin{equation}\label{frame2-norm}
\|f\|_2^2 =
\sum_{j \geq 0} \sum_{\xi \in \XX_j}
\overline{\langle f, \tilde\psi_{j \xi}\rangle} \langle f, \psi_{j \xi}\rangle
\quad \forall f\in \LL^2.
\end{equation}

$(b)$ {\rm Space localization:} For any $\sig>0$ there exists a constant $c_\sig>0$ such that
for any $\xi\in\XX_j$, $j\ge 0$,
\begin{equation}\label{local-frame2}
|\psi_{j \xi} (x)|, |\tilde\psi_{j \xi} (x)|
\le c_\ssigma|B(\xi, b^{-j})|^{-1/2}(1+b^j\rho(x, \xi))^{-\ssigma},
\end{equation}
and if $\rho(x, y) \le b^{-j}$
\begin{equation}\label{Lip-frame2}
|\psi_{j \xi} (x)- \psi_{j \xi} (y)|
\le c_\ssigma|B(\xi, b^{-j})|^{-1/2}(b^j\rho(x, y))^\alpha (1+b^j\rho(x, \xi))^{-\ssigma}.
\end{equation}
Here $\alpha>0$ is the global parameter from $(\ref{lip})$
and the same inequality hold for $\tilde\psi_{j \xi}$ in place of $\psi_{j \xi}$.

$(c)$ {\rm Spectral localization:}
$\psi_{0\xi}, \tilde\psi_{0\xi}\in \Sigma_b^p$ if $\xi\in \XX_0$ and
$\psi_{j\xi},\tilde\psi_{j\xi}\in \Sigma_{[b^{j-1}, b^{j+1}]}^p$
if $\xi\in \XX_j$, $j\ge 1$, $0<p\le\infty$.

$(d)$ {\rm Norms:}
\begin{equation}\label{norm-frame2}
\|\psi_{j \xi}\|_p \sim \|\tilde\psi_{j \xi}\|_p \sim |B(\xi, b^{-j})|^{\frac 1p-\frac 12},
\quad 0< p \le \infty.
\end{equation}
\end{proposition}

\noindent
{\bf Proof.}
Identities (\ref{frame2-rep}) follow immediately from (\ref{calderon}) and (\ref{Psi-Psi}).
For the proof of (\ref{frame2-norm}), denote
$
S_Nf = \sum_{j=0}^N\sum_{\xi \in \XX_j} \langle f, \tilde\psi_{j\xi}\rangle \psi_{j\xi}
$
and observe that
$$
\|f\|_2^2 =\lim_{N\to \infty}\langle f, S_Nf \rangle
= \lim_{N\to \infty}\sum_{j=0}^N\sum_{\xi \in \XX_j}
\overline{\langle f, \tilde\psi_{j \xi}\rangle} \langle f, \psi_{j \xi}\rangle
= \sum_{j\ge 0}\sum_{\xi \in \XX_j}
\overline{\langle f, \tilde\psi_{j \xi}\rangle} \langle f, \psi_{j \xi}\rangle.
$$
The localization and Lipschitz property of the frame elements given in
(\ref{local-frame2}) and (\ref{Lip-frame2}) follow by Theorem~\ref{thm:main-local-kernels}.
The claimed spectral localization is obvious.
The~norm bounds in (\ref{norm-frame2}) follow by Theorem~\ref{thm:norms}.
$\qed$

\smallskip

An interesting special case of the above construction occurs when we choose
$\Psi_0 = \tilde\Psi_0$ and $\Psi=\tilde\Psi$.
Then $\psi_{j\xi}=\tilde\psi_{j\xi}$ and $\{\psi_{j\xi}\}$ is a tight frame for $\LL^2$,
i.e.
$$
\|f\|_2^2 = \sum_{j\ge 0}\sum_{\xi \in \XX_j}
|\langle f, \psi_{j\xi}\rangle|^2,
\quad \forall f\in \LL^2.
$$


\noindent
{\bf Remark.}
The  polynomial property (\ref{polyp}) of the spectral spaces apparently is valid when the spectral
functions are polynomials. This simple fact has been utilized for construction of frames
on the sphere \cite{NPW1}, on the interval with Jacobi weights \cite{PX1}, on the ball \cite{PX2},
and in the context of Hermite \cite{PX3} and Laguerre \cite{KPPX} expansions.

\section{Besov spaces}\label{besov-spaces}
\setcounter{equation}{0}

We shall follow the general idea of using spectral decompositions, e.g. \cite{Peetre, Triebel1, Triebel2},
to introduce (inhomogeneous) Besov spaces in the general set-up of this paper.
As explained in the introduction, we shall only consider Besov spaces $B^s_{pq}$
with $s>0$ and $1\le p \le \infty$.
The Besov spaces $B^s_{pq}$ with full range of indices are treated in the follow-up paper \cite{KP}. For another approach to Besov spaces under heat kernel estimates, but a polynomial upper bound on the volume instead of the volume doubling condition, see \cite{BDY}.

To introduce Besov spaces we assume that there are given two (Littlewood-Paley) functions
$\varphi_0, \varphi\in C^\infty(\bR_+)$ such that
\begin{align}
&\supp \varphi_0 \subset   [0, 2] ,\;\;
\varphi_0^{(\nu)}(0) = 0  \hbox{ for } \nu\ge 1,\;\;
|\varphi_0(\lam)| \ge c>0 \;\hbox{ for } \lam\in [0, 2^{3/4}], \label{cond_phi}\\
&\supp \varphi \subset   [1/2, 2], \;\;
|\varphi(\lam)| \ge c>0 \;\hbox{ for } \lam\in [2^{-3/4}, 2^{3/4}]. \label{cond_psi}
\end{align}
Then
$|\varphi_0(\lam)| +\sum_{j\ge 1} |\varphi(2^{-j}\lam)| \ge c >0$
for $\lam \in  [0, \infty)$.
Set $\varphi_j(\lambda):= \varphi(2^{-j}\lambda)$ for $j\ge 1$.


\begin{definition}\label{def-B-spaces}
Let $s>0$, $1\le p\le \infty$, and $0<q \le \infty$.
The Besov space  $B_{pq}^{s}=B_{pq}^{s}(L)$
is defined as the set of all $f \in \LL^p$ such that
\begin{equation}\label{def-Besov-space1}
\|f\|_{B_{pq}^{s}} :=
\Big(\sum_{j\ge 0} \Big(2^{s j}
\|\varphi_j(\sqrt{L}) f(\cdot)\|_{\Lp}
\Big)^q\Big)^{1/q} <\infty.
\end{equation}
Here the $\ell^q$-norm is replaced by the sup-norm if $q=\infty$.
\end{definition}

Note that by Proposition~\ref{prop:character-Besov} below it follows that the definition of
the Besov spaces $B_{pq}^{s}$ is independent of the specific selection of $\varphi_0$, $\varphi$
satisfying (\ref{cond_phi})-(\ref{cond_psi}).
Also, $B_{pq}^{s}$ are (quasi-)Banach spaces,
which are continuously embedded in $\LL^p$ as will be seen below.

\subsection{Characterization of Besov spaces via linear approximation from \boldmath $\{\Sigma_t^p\}$}
\label{sec:char-Besov}

Here we show that the Besov spaces $B_{pq}^s$ with $s>0$ and $p\ge 1$
are in fact the approximation spaces of linear approximation from $\Sigma_t^p$, $t\ge 1$.
As in \S\ref{linear-app}, we let $\cE_t(f)_p$ denote the best approximation of
$f \in \LL^p$ from $\Sigma_t^p$
and $A_{pq}^s$ will denote the associated approximation spaces,
defined in (\ref{def:app-space-p})-(\ref{def:app-space-infty}).


\begin{proposition}\label{prop:character-Besov}
Let $s > 0$, $1 \le p \le \infty$, and $0 < q \le \infty$.
Then $f \in B_{pq}^s$ if and only if  $f\in A_{pq}^s$.
Moreover,
\begin{equation}\label{character-Besov2}
\|f\|_{B_{pq}^s} \sim \|f\|_{A_{pq}^s}
:= \|f\|_p +
\Big(\sum_{j\ge 0}\big(2^{s j}\cE_{2^j}(f)_p \big)^q\Big)^{1/q}.
\end{equation}
\end{proposition}

\noindent
{\bf Proof.}
Let $\varphi_j$ be as in the definition of the Besov spaces with the additional property:
$
\sum_{j\ge 0} \varphi_j(\lambda)=1
$
for $\lambda\in [0, \infty)$ (see \S\ref{app-identity}).
Suppose $f\in \LL^p$.
Then by Corollary~\ref{cor:Littlewood-Paley} we have
$f=\sum_{j\ge 0}\varphi_j(\sqrt L)f$
and since $\varphi_j(\sqrt L)f\in \Sigma_{[2^{j-1}, 2^{j+1}]}^p$ we obtain
$$
\cE_{2^m}(f)_p\le \sum_{j\ge m} \|\varphi_j(\sqrt L)f\|_p
$$
and the standard Hardy inequality
\begin{equation}\label{hardy}
\sum_{m\ge 0}\Big(2^{sm}\sum_{j\ge m}b_j\Big)^q
\le c \sum_{m\ge 0}\Big(2^{sm}b_m\Big)^q,
\quad b_j \ge 0,\; s>0,\; 0<q\le \infty,
\end{equation}
leads to the estimate
$\|f\|_{A_{pq}^s} \le c\|f\|_{B_{pq}^s}$.

For the estimate in the other direction,
we note that for any $g\in \Sigma_{2^{j-1}}^p$ we have
$\varphi_j(\sqrt L) f = \varphi_j(\sqrt L) (f-g)$
and hence
$$
\|\varphi_j(\sqrt L) f\|_p = \|\varphi_j(\sqrt L) (f-g)\|_p \le c\|f-g\|_p,
$$
where we used the boundedness of the operator $\varphi_j(\sqrt L)$
on $\LL^p$.
This implies $\|\varphi_j(\sqrt L) f\|_p \le c\cE_{2^{j-1}}(f)_p$, $j\ge 1$,
and obviously $\|\varphi_0(\sqrt L) f\|_p \le c\|f\|_p$.
We use these estimates in the definition of $B_{pq}^s$ to obtain
$\|f\|_{B_{pq}^s} \le c\|f\|_{A_{pq}^s}$.
$\qed$

We next record the heat kernel characterization of Besov spaces.
Denote
\begin{equation}\label{heat-Besov}
\|f\|_{B^s_{pq}(H)} := \| f \|_p +
\Big(\int_0^1 \big(t^{-s/2} \|(tL)^m e^{-tL}f\|_p \big)^q \frac{dt}t\Big)^{1/q}
\end{equation}
with the usual modification for $q=\infty$.


\begin{corollary}\label{cor:Heat-Besov}
For admissible indices $s, p, q$ a function $f\in B^s_{pq}$ if and only if
$\|f\|_{B^s_{pq}(H)}<\infty$ and if $f\in B^s_{pq}$, then
$\|f\|_{B^s_{pq}} \sim \|f\|_{B^s_{pq}(H)}$.
\end{corollary}

This corollary follows readily by Proposition~\ref{prop:character-Besov} taking into account
Remark~\ref{semig}.

\subsection{Comparison of Lipschitz spaces and \boldmath $B^s_{\infty \infty}$}
The Lipschitz space $\lip \gamma$, $\gamma>0$, is defined as the set of all $f\in\LL^\infty$
such that
\begin{equation}\label{lip-norm}
\|f\|_{\lip \gamma} := \| f \|_\infty + \sup_{x\neq y} \frac{| f(x)-f(y) |}{\rho^\gamma(x,y)} <\infty.
\end{equation}

We would like to record next the fact that in the setting of this article the spaces
$\lip s$ and $B^s_{\infty \infty}$ coincide provided $0<s<\alpha$,
where $\alpha$ is the structural constant from (\ref{lip}).


\begin{proposition}\label{prop:lip}
The following continuous embeddings hold:
$(a)$ For any $s>0$
$$\lip s \subset B^s_{\infty \infty}.$$
$(b)$ For any $0<s< \alpha$
$$B^s_{\infty \infty} \subset  \lip s.$$
\end{proposition}

\noindent
{\bf Proof.}
(a) Let $f\in \lip s$ and choose
$\theta \in C^\infty[0, \infty)$ so that $\theta \ge 0$, $\theta \equiv 1$ on $[0,1]$
$\supp \theta \subset[0,2]$.
Then using Theorem~\ref{thm:main-local-kernels} and (\ref{INT}) we obtain for $t\geq 1$
and $k>s+3d/2$
\begin{align*}
|\theta (t^{-1} \sqrt L)f(x)- f(x)|
&= \Big|\int_M \theta(t^{-1} \sqrt L))(x,y) [f(y)-f(x)] d\mu(y)\Big|\\
&\le c\|f\|_{\lip s}\int_M D_{t^{-1},k} (x,y) \rho^s(x,y) d\mu(y)\\
&\le ct^{-s}\|f\|_{\lip s}\int_M D_{t^{-1},k-s} (x,y) d\mu(y)
\le ct^{-s}\| f\|_{\lip s}.
\end{align*}
On the other hand $\theta (t^{-1} \sqrt L)f \in \Sigma_{2t}^\infty$ and hence
$\cE_{2t}(f)_\infty \le \|\theta (t^{-1} \sqrt L)f-f\|_\infty$.
From this and above we infer
$\cE_{2t}(f)_\infty \le ct^{-s}\|f\|_{\lip s}$, which implies (a).

(b)
Let $\varphi_0:=\theta$ with $\theta$ the function from above.
Set
$\varphi(\lambda):=\theta(\lambda) - \theta(2\lambda)$ and
$\varphi_j(\lambda):=\varphi(2^{-j}\lambda)$.
Then
$\sum_{j\ge 0}\varphi_j(\lambda) = 1$ for $\lambda\ge 0$,
$\supp \varphi_0\subset [0, 2]$ and $\supp\varphi_j\subset [2^{j-1}, 2^{j+1}]$, $j\ge 1$.
Now, assuming that $f \in B^s_{\infty \infty }$ we apparently have
$\varphi_0(\sqrt L)f  \in \Sigma^\infty_2$,
$\varphi_j(\sqrt L)f  \in \Sigma^\infty_{2^{j+1}}$, and
by the Littlewood-Paley decomposition (Corollary~\ref{cor:Littlewood-Paley})
$f = \sum_{j \ge 0}\varphi_j(\sqrt L)f$.
Evidently, $B^s_{\infty \infty}$ can be defined using the above constructed functions $\{\varphi_j\}$
and hence
$\|\varphi_j(\sqrt L)f \|_\infty \le c2^{-js}\|f\|_{B^s_{\infty \infty}}$, $j\ge 0$.
Therefore, using (\ref{lip2}) we have for $0<s<\alpha$ and any $J\ge 1$
\begin{align*}
|f(x)-f(y)|
& \le \sum_{j \ge 0}|\varphi_j(\sqrt L)f(x)- \varphi_j(\sqrt L)f(y)|\\
&\le c\sum_{j=0}^J \|\varphi_j(\sqrt L)f\|_\infty \big(2^{j}\rho(x,z)\big)^\alpha
+2\sum_{j>J} \|\varphi_j(\sqrt L)f  \|_\infty\\
&\le c\|f\|_{B^s_{\infty \infty}}
\Big(\sum_{j=0}^J 2^{-js}\big(2^{j}\rho(x,z)\big)^\alpha + \sum_{j>J} 2^{-js}\Big)\\
&\le c\|f\|_{B^s_{\infty \infty}}\big(2^{J(\alpha-s)}\rho(x,z)^\alpha + 2^{-Js}\big).
\end{align*}
Assuming that $0<\rho(x, y) \le 1$ we choose $J\ge 1$ so that $2^{-J}\sim \rho(x, y)$ and the above yields
$|f(x)-f(y)|\le c\|f\|_{B^s_{\infty \infty}}\rho(x, y)^s$.
If $\rho(x, y)>1$ this estimate is immediate from
$\|f\|_\infty\le c\|f\|_{B^s_{\infty \infty}}$,
which follows trivially using the decomposition of $f$ from above.
This completes the proof of (b).
$\qed$

\subsection{Frame decomposition of Besov spaces}

Our aim here is to show that the Besov spaces introduced by Definition~\ref{def-B-spaces}
can be characterized in terms of respective sequence norms of the frame coefficients
of functions, using the frames constructed in \S\ref{sec:frames}.
We shall utilize the pair of dual frames $\{\psi_{j\xi}\}$, $\{\tilde\psi_{j\xi}\}$
constructed in \S\S\ref{natural-frame}-\ref{dual-frame} or in \S\ref{sec:polyn-prop}.
To make the idea of frame decomposition of $B^s_{pq}$ more transparent we first introduce
the sequence B-spaces $b^s_{pq}$.


\begin{definition}\label{def:b-spaces}
For $s>0$, $1\le p\le\infty$, and $0<q \le \infty$
the sequence space $b_{pq}^s$
is defined as the space of all complex-valued sequences
$a:=\{a_{j\xi}: j\ge 0, \xi\in \XX\}$ such that
\begin{equation}\label{def-tilde-berpq}
\|a\|_{b_{pq}^s}
:=\Bigl(\sum_{j\ge 0}b^{jsq}
\Bigl[\sum_{\xi\in \XX_j}\Big(|B(\xi, \bb^{-j})|^{1/p-1/2}|a_{j\xi}|\Big)^p
\Bigr]^{q/p}\Bigr)^{1/q} <\infty.
\end{equation}
Here $b>1$ is the constant from \S\ref{sec:frames},
and the $\ell^p$ or $\ell^q$ norm is replaces by the $\sup$-norm if $p=\infty$ or $q=\infty$.
\end{definition}

In our further analysis we shall use the ``analysis" and ``synthesis" operators
defined by
\begin{equation}\label{anal_synth_oprts}
S_{\til\psi}: f\rightarrow \{\langle f, \til\psi_{j\xi}\rangle\}
\quad\text{and}\quad
T_{\psi}: \{a_{j\xi}\}\rightarrow \sum_{j\ge 0}\sum_{\xi\in \XX_j}a_{j\xi}\psi_{j\xi}.
\end{equation}


\begin{theorem}\label{thm:B-character}
Let $s>0$, $1\le p\le\infty$, and $0<q \le \infty$.
Then the operators
$S_{\til\psi}: B_{pq}^s \rightarrow  b_{pq}^s$ and
$T_{\psi}: b_{pq}^s \rightarrow B_{pq}^s$
are bounded and $T_{\psi} S_{\til\psi}=\Id$ on $B_{pq}^s$.
Consequently, $f\in B_{pq}^s$ if and
only if $\{\langle f, \til\psi_{j\xi}\rangle\}\in b_{pq}^s$.
Moreover, if $f\in B_{pq}^s$, then
\begin{equation}\label{Bnorm-equivalence-1}
\|f\|_{B_{pq}^s}
\sim  \|\{\langle f,\tilde\psi_{j\xi}\rangle\}\|_{b_{pq}^s}
\sim \Big(\sum_{j\ge 0} b^{jsq}\Bigl[\sum_{\xi\in \XX_j}
\|\langle f,\til\psi_{j\xi}\rangle\psi_{j\xi}\|_p^p\Bigr]^{q/p}\Bigr)^{1/q}
\end{equation}
with the usual modification when $p=\infty$ or $q=\infty$.
Above the roles of $\{\psi_{j\xi}\}$ and $\{\tilde\psi_{j\xi}\}$ can be interchanged.
\end{theorem}

\noindent
{\bf Proof.}
Let $\Psi_j\in C^\infty_0$, $j\ge 0$,
be the functions from the definition of  the frames in \S\ref{natural-frame}.
Recall that $\supp \Psi_0\subset [0, b]$ and $\supp \Psi_j\subset [b^{j-1}, b^{j+1}]$, $j\ge 1$.
Also, $\sum_{j\ge 0}\Psi_j(u)=1$, $u\in\R_+$, and hence
$f=\sum_{j\ge 0}\Psi_j(\sqrt{L})f$ for $f\in \LL^p$.
It is easy to see that
Proposition~\ref{prop:character-Besov} implies (with the obvious modification when $q=\infty$)
\begin{equation}\label{equiv-norms}
\|f\|_{B^s_{pq}} \sim \|f\|_{A^s_{pq}}
\sim \|f\|_p + \Big(\sum_{j\ge 0}\big(b^{js}\cE_{b^j}(f)_p\big)^q\Big)^{1/q}
\sim \Big(\sum_{j\ge 0}\big(b^{js}\|\Psi_j(\sqrt{L})f\|_p\big)^q\Big)^{1/q}.
\end{equation}
Here the second equivalence follows by the monotonicity of $\cE_t(f)_p$ and
the last equivalence follows exactly as in the proof of Proposition~\ref{prop:character-Besov}.

Let $f\in B^s_{pq}$ and assume $q<\infty$ (the case $q=\infty$ is easier).
By (\ref{prop-psi-2}) and (\ref{def-tilde-berpq}) it follows that
\begin{equation}\label{est-S}
\|S_{\tilde \Psi}f\|_{b_{pq}^s}
=  \|\{\langle f,\tilde\psi_{j\xi}\rangle\}\|_{b_{pq}^s}
\sim \Big(\sum_{j\ge 0} b^{jsq}\Bigl[\sum_{\xi\in \XX_j}
\|\langle f,\til\psi_{j\xi}\rangle\psi_{j\xi}\|_p^p\Bigr]^{q/p}\Bigr)^{1/q}.
\end{equation}
Using that $f=\sum_{j\ge 0}\Psi_j(\sqrt{L})f$,
$\Psi_j(\sqrt{L})(\cdot, y) \in \Sigma^2_{[b^{j-1}, b^{j+1}]}$, and
$\tilde\psi_{j\xi} \in \Sigma^2_{[b^{j-2}, b^{j+2}]}$ we obtain
$$
\langle f,\til\psi_{j\xi}\rangle\psi_{j\xi}
= \sum_{\nu=j-2}^{j+2}\langle \Psi_\nu(\sqrt{L})f,\til\psi_{j\xi}\rangle\psi_{j\xi},
\quad \xi\in \cX_j,
$$
where
$\Psi_\nu(\sqrt{L}):= 0$ if $\nu <0$.
This readily implies
$$
\sum_{\xi\in \XX_j}\|\langle f,\til\psi_{j\xi}\rangle\psi_{j\xi}\|_p^p
\le c\sum_{\nu=j-2}^{j+2}\|\langle \Psi_\nu(\sqrt{L})f,\til\psi_{j\xi}\rangle\psi_{j\xi}\|_p^p
\le c\sum_{\nu=j-2}^{j+2}\|\Psi_\nu(\sqrt{L})f\|_p^p.
$$
Here for the last inequality we used Lemma~\ref{lem:technical}, (a).
We insert the above in (\ref{est-S}) and use (\ref{equiv-norms}) to obtain
$
\|S_{\tilde\Psi}f\|_{b^s_{pq}}
\le c \|f\|_{B^s_{pq}}.
$
Hence the operator $S_{\til\psi}: B_{pq}^s \rightarrow  b_{pq}^s$ is bounded.


To prove the boundedness of $T_{\psi}: b_{pq}^s \rightarrow B_{pq}^s$,
we assume that $a=\{a_{j\xi}\}\in b^s_{pq}$ and denote briefly
$
f=T_{\psi}a=\sum_{j\ge 0}\sum_{\xi\in\XX_j} a_{j\xi}\psi_{j\xi}.
$
Assume $q<\infty$ (the case $q=\infty$ is easier).
Using (\ref{technical-2}), H\"{o}lder's inequality if $q>1$, and (\ref{prop-psi-2}) we obtain
\begin{equation}\label{estimate-f}
\|f\|_p \le c\sum_{j\ge 0}\Big(\sum_{\xi\in\XX_j} \|a_{j\xi}\psi_{j\xi}\|_p^p\Big)^{1/p}
\le c\Big(\sum_{j\ge 0}\Big[b^{sj}\sum_{\xi\in\XX_j}
\|a_{j\xi}\psi_{j\xi}\|_p^p\Big]^{q/p}\Big)^{1/q}
\le c\|a\|_{b^s_{pq}}.
\end{equation}
Therefore, $T_{\psi}a$ is well-defined.
Further, since $\psi_{j\xi}\in \Sigma_{b^{j+1}}^p$ and applying again (\ref{technical-2})
we get
$$
\cE_{b^j}(f)_p \le \big\|\sum_{m\ge j}\sum_{\xi\in\XX_m} a_{m\xi}\psi_{m\xi}\big\|_p
\le c\sum_{m\ge j}\Big(\sum_{\xi\in\XX_m} \|a_{m\xi}\psi_{m\xi}\|_p^p\Big)^{1/p}.
$$
This and the Hardy inequality (\ref{hardy}) give
$
\Big(\sum_{j\ge 0}\big(b^{js}\cE_{b^j}(f)_p\big)^q\Big)^{1/q} \le c\|a\|_{b^s_{pq}}.
$
In turn, this and (\ref{equiv-norms}) yield
$\|f\|_{B^s_{pq}} \le c\|a\|_{b^s_{pq}}$.
Thus the operator $T_{\psi}: b_{pq}^s \rightarrow B_{pq}^s$ is also bounded.

The identity $T_{\psi} S_{\til\psi}=\Id$ on $B_{pq}^s$ follows by (\ref{rep-L2}).
$\qed$

\subsection{Embedding of Besov spaces}

Finally we show that the Besov spaces $B^s_{pq}$ embed ``correctly".


\begin{proposition}\label{B-embedding}
Let $1\le p\le p_1<\infty$, $0<q\le q_1\le \infty$, $0<s_1\le s<\infty$.
Then we have the continuous embedding
\begin{equation}\label{B-embed}
B_{pq}^{s} \subset B_{p_1q_1}^{s_1}
\quad\mbox{if}\quad
s/d-1/p=s_1/d-1/p_1.
\end{equation}
\end{proposition}

\noindent
{\bf Proof.} This assertion follows easily by Proposition \ref{prop:Nikolski}.
Indeed, let $\{\varphi_j\}_{j\ge 0}$ be the functions from the definition of Besov spaces
(Definition~\ref{def-B-spaces}).
Given $f\in B_{pq}^{s}$
we evidently have $\varphi_j(\sqrt{L}) f\in \Sigma_{2^{j+1}}^p$ and
using  (\ref{norm-relation1})
\begin{align*}
\|\varphi_j(\sqrt{L}) f(\cdot)\|_{p_1}
\le c2^{jd(1/p_1-1/p)}\|\varphi_j(\sqrt{L}) f(\cdot)\|_{p},
\end{align*}
which readily implies
$\|f\|_{B_{p_1q_1}^{s_1}} \le c\|f\|_{B_{pq}^{s}}$.
$\qed$

\smallskip

Compare the above result with \cite{CSa}, where embeddings between Besov spaces defined via the heat semigroup
are proved under an assumption of polynomial decay  of the heat kernel.

\section{Heat kernel on $[-1,1]$ induced by the Jacobi operator}\label{sec:Jacobi}
\setcounter{equation}{0}

We consider the case when
$M =[-1, 1]$, $d\mu(x) = \W(x) dx$, where
$$
\W(x)=w(x)= (1-x)^\alpha (1+x)^\beta, \quad \alpha , \beta > -1,
$$
and
$$
Lf(x) =- \frac{[w(x) a(x)f'(x)]'}{w(x)}, \quad a(x) =(1-x^2),
\quad D(L) = C^2[-1,1].
$$
Integrating by parts we get
$
\cE(f,g) =\langle Lf, g \rangle = \int_{-1}^1 a(x) f'(x) g'(x) w(x) dx.
$
Clearly, the domain $D(\overline \cE)$ of the closure $\overline{\cE}$ of $\cE$ is given by the set of
weakly differentiable functions $f$ on $]-1,1[$ such that
$$
\|f\|^2_{\cE}
= \int_{-1}^1 |f(x)|^2 w(x) dx +   \int_{-1}^1 a(x) |f'(x)|^2 w(x) dx <\infty.
$$
Note that $D(L)\supset \cP$ the space of all polynomials,
and $L(\cP_k) \subset \cP_k$, $k\ge 0$, with $\cP_k$ being the space of all polynomials of degree $k$.
As is well known \cite{Sz} the (normalized) Jacobi polynomials $P_k$, $k=0, 1, \dots$, are eigenfunctions
of $L$, i.e. $LP_k=\lambda_kP_k$ with $\lambda_k=k(k+\alpha+\beta+1)$.
By the density of polynomials in $L^2([-1, 1], \mu)$ it follows that
$$
e^{-t\bar L} (f) = \sum_{k\ge 0} e^{-\lambda_k t} \langle f, P_k \rangle P_k, \quad t>0.
$$

We next show that  $e^{-t\bar L}$ is submarkovian.
Let $\Phi_\varepsilon \in C^\infty(\R)$ and $0 \leq \Phi_\varepsilon' \leq 1$.
Then for any $f \in C^2[-1,1]$ we have
$(\Phi_\varepsilon (f))' =\Phi'_\varepsilon (f)f' \in C^2[-1,1]$
and
\begin{align*}
\cE(\Phi_\varepsilon (f),\Phi_\varepsilon (f))
=\int_{-1}^1 a(x) |(\Phi_\varepsilon (f))'|^2w(x) dx
\le \int_{-1}^1 a(x)|f'(x)|^2w(x) dx= \cE(f,f).
\end{align*}
Hence $e^{-t\bar L}$ is submarkovian (see \S\ref{sec:Dirichle-spaces}).

Moreover, this Dirichlet space is evidently strongly local and regular
and also $\Gamma(f,g)(u)= a(u) f'(u)g'(u)$.

We now compute the  intrinsic  metric. We have for $x, y\in [-1, 1]$, $x<y$,
\begin{align*}
\rho(x,y) &= \sup\{u(x)-u(y): u\in C^2[-1,1], a(x)|u'(x)| \le 1\}\\
&= \sup\{\int_x^y u'(t) dt: u \in C^2[-1,1], a(t) |u'(t)|^2 \leq 1\}\\
&= \int_x^y \frac{dt}{\sqrt{a(t)}} dt = |\arccos x-\arccos y|.
\end{align*}
Evidently, the topology generated by this metric is the usual topology on $[-1, 1]$,
and $[-1, 1]$ is complete.

\smallskip

It remains to verify the doubling property of the measure and the scale-invariant Poincar\'e inequality.

\subsection{Doubling property of the measure}

The doubling property of the measure $d\mu(x)=\W(x)dx$ follows readily by
the following estimates on $|B(x, r)|$:
For any $x\in [-1, 1]$ and $0<r\le \pi$
\begin{equation}\label{measure-ball}
c_1|B(x, r)| \le r(1-x+r^2)^{\alpha+1/2}(1+x+r^2)^{\beta+1/2} \le c_2|B(x, r)|,
\end{equation}
where $c_1, c_2>0$ are constants depending only on $\alpha$ and $\beta$.

To prove these estimates, assume that $x=\cos\theta$, $0\le \theta\le \pi$.
Then evidently
$|B(x, r)|=\int_{\cos[\pi\wedge(\theta+r)]}^{\cos[0\vee (\theta-r)]}\W(u)du$,
where  $a\vee b:= \max\{a, b\}$ and $a\wedge b:= \min\{a, b\}$ as usual.
Assume $0\le x\le 1$ and $0<r\le \pi/4$.
The following chain of similarities with constants depending only on $\alpha, \beta$ is quite obvious:
\begin{align*}
|B(x, r)|
&\sim\int_{\cos(\theta+r)}^{\cos[0\vee (\theta-r)]}(1-u)^\alpha du
\sim \int^{\theta+r}_{0\vee (\theta-r)}(1-\cos \varphi)^\alpha\sin\varphi d\varphi\\
& \sim \int^{\theta+r}_{0\vee (\theta-r)}\varphi^{2\alpha+1} d\varphi
\sim r(\theta+r)^{2\alpha+1} \sim r(\sin\theta+r)^{2\alpha+1}\\
&\sim r(\sqrt{1-x^2}+r)^{2\alpha+1} \sim r(1-x+r^2)^{\alpha+1/2},
\end{align*}
which implies (\ref{measure-ball}).
The case when $-1\le x <0$ and $0<r\le \pi/4$ is similar
and in the case $\pi/4<r\le \pi$ we obviously have $|B(x, r)|\sim 1$,
which again leads to (\ref{measure-ball}).

\subsection{Poincar\'e inequality}
As was explained in \S\ref{sec:Dirichle-spaces} a critical ingredient in establishing Gaussian bounds
for the heat kernel is the scale-invariant Poincar\'e inequality,
which we establish next.


\begin{theorem}\label{thm:poincare}
For any $f\in D(\overline \cE)$ and an interval $I=[a,b] \subset [-1,1]$
\begin{equation}\label{poincare}
\int_I  |f(x)- f_I|^2 w(x) dx
\le c (\diam_\rho (I))^2\int_I |f'(x)|^2 (1-x^2) w(x) dx
\end{equation}
where
$\diam_\rho (I)= \arccos a - \arccos b$, 
$f_I = \frac 1{w(I)}\int_I f(x) w(x) dx$  with $w(I)= \int_Iw(x)dx$,
and $c>0$ is a constant depending only on $\alpha, \beta$.

\end{theorem}

\noindent
{\bf Proof.}
Denote briefly
$w[c, d]:= \int_c^dw(u)du$.
We have for $I=[a,b] \subset [-1,1]$ and $x\in I$
\begin{align*}
f(x) - f_I
&= \frac 1{w(I)} \int_I  (f(x)-f(y)) w(y) dy
= \frac 1{w(I)} \int_I  \int_y^x f'(u)du w(y) dy\\
&= \int_I f'(u) K(x,u) du,
\end{align*}
where
$
K(x,u) :=  \frac 1{w(I)} \{w[a,u ] \ONE_{[a,x]} -  w[u,b]\ONE_{[x, b]}\}.
$
It is easy to see that
$$
\int_I |K(x,u)| du =\frac 1{w(I)} \int_I |x-y| w(y) dy \le \frac 12(b-a)
\quad\hbox{and}
$$
$$
\int_I |K(x,u)|w(x) dx =\frac{2w[a,u ]w[u,b ]}{w(I)}.
$$
Using the above we obtain
\begin{align*}
\int_I |f(x) - f_I |^2 w(x) dx
&= \int_I \Big|\int_I f'(u) K(x,u)du\Big|^2 y(x) dx\\
&\le \int_I\Big( \int_I  |f'(u)|^2 | K(x,u)| du \int_I | K(x,v)| dv \Big) w(x) dx\\
&\le \frac 12(b-a)\int_I |f'(u)|^2 \Big(\int_I K(x,u)|w(x)dx\Big))du\\
&=(b-a)  \int_I  |f'(u)|^2   \frac{w[a,u ]w[u,b ]}{w(I)}du.
\end{align*}
Therefore, the theorem will be proved if we show that
\begin{equation}\label{Poincare2}
(b-a)  \frac{w[a,u ] w[u,b ]}{w(I)}
\le  cw(u)(1-u^2)\Big(\int_a^b \frac {dz}{\sqrt{1-z^2}}\Big)^2
\end{equation}
for some constant $c>0$ depending only on $\alpha, \beta$.

Suppose
$[a,b]\subset [-1/2,1]$.
Then it is readily seen that
\begin{equation}\label{poincare3}
w(u)(1-u^2)\Big(\int_a^b \frac {dz}{\sqrt{1-z^2}}\Big)^2
\ge 2^{-\beta}(1-u)^{\alpha +1}(\sqrt{1-a}- \sqrt{1-b})^2.
\end{equation}
On the other hand, since
$w(x)\le 2^{|\beta|}(1-x)^\alpha$, we have
$$
\frac{w[a,u ]w[u,b ]}{w(I)}
\le \frac{2^{3|\beta|}}{\alpha +1}\frac{[ (1-a)^{\alpha +1} - (1-u)^{\alpha +1} ]
[ (1-u)^{\alpha +1} - (1-b)^{\alpha +1} ] }
{(1-a)^{\alpha +1} - (1-b)^{\alpha +1}}.
$$
We need the following inequality whose proof is straightforward:
If $\gamma >0$ and $0 \leq A \leq X \leq B$, then
\begin{equation}\label{inequality}
\frac{ (X^\gamma-A^\gamma)(B^\gamma  -X^\gamma)}{B^\gamma-A^\gamma}
\le (\gamma \vee 1) X^\gamma  \frac{\sqrt B - \sqrt A}{\sqrt B + \sqrt A}.
\end{equation}
Applying this inequality we get
\begin{align*}
(b-a) \frac{w[a,u ]w[u,b ]}{w(I)}
&\le  2^{3|\beta|}\big(\frac 1{\alpha+1} \vee 1\big)
(1-u)^{\alpha +1}(b-a) \frac{\sqrt{1-a}- \sqrt{1-b}}{\sqrt{1-a}+ \sqrt{1-b}}\\
&= 2^{3|\beta|}\big(\frac 1{\alpha+1} + 1\big)(1-u)^{\alpha +1}  (\sqrt{1-a}- \sqrt{1-b})^2.
\end{align*}
This coupled with (\ref{poincare3}) gives (\ref{Poincare2}).
The proof of (\ref{Poincare2}) in the case when $I= [a,b]\subset [-1,1/2]$ is the same.

Let now $-1 \le a < -1/2 <1/2 <b \le 1$.
Suppose $u\in [0, b]$ (the case when $u\in [a, 0)$ is similar).
Then evidently $w[a, u]\sim 1$, $w(I)\sim 1$,
$\int_a^b \frac {dz}{\sqrt{1-z^2}} \sim 1$
and (\ref{Poincare2}) follows by
$$
w[u, b]\le 2^{|\beta|}\int_u^1(1-y)^\alpha dy \le \frac{2^{|\beta|}}{\alpha+1}(1-u)^{\alpha+1}
\;\, \hbox{and}\;\, w(u)(1-u^2) \sim (1-u)^{\alpha+1}.
$$
The proof of the theorem is complete.
$\qed$

\subsection{Gaussian bounds on the heat kernel associated with the Jacobi operator}

As a consequence of the Poincar\'e inequality and the doubling property of the measure,
established above, we obtain (\S\ref{sec:Dirichle-spaces})
Gaussian bounds for the heat kernel $p_t(x, y)$ associated
with the Jacobi operator:


\begin{theorem}\label{thm:HK-Jacobi}
For any $x, y\in [-1, 1]$ and $0<t\le 1$,
\begin{equation}\label{HK-Jacobi-1}
\frac{c_1'\exp\{-\frac{c_1\rho^2(x, y)}{t}\}}{\sqrt{|B(x,\sqrt t)| |B(y,\sqrt t)}}
\le  p_t(x,y)
\le \frac{c_2'\exp\{-\frac{c_2\rho^2(x, y)}{t}\}}{\sqrt{|B(x,\sqrt t)| |B(y,\sqrt t)}}.
\end{equation}
Here
$|B(x, \sqrt t)| \sim \sqrt t(1-x+t)^{\alpha+1/2}(1+x+t)^{\beta+1/2}$,
$\rho(x, y)=|\arccos x - \arccos y|$ or
$\rho(x, y)=|\theta - \phi|$ if
$x= \cos \theta$ and $y= \cos \phi$,  $0  \le \theta, \phi\le \pi$,
and $c_1, c_2, c_1', c_2'>0$ are constants depending only on $\alpha$ and $\beta$.

Furthermore,
\begin{equation}\label{HK-Jacobi-2}
p_t(x,y) = \sum_{k\ge 0} e^{-\lambda_k t}  P_k(x)  P_k(y), \quad \lambda_k = k(k+\alpha + \beta +1),
\end{equation}
where the series converges uniformly.
\end{theorem}

The above results and Theorem~\ref{thm:main-local-kernels} yield the nearly exponential localization
of kernels as in the following


\begin{corollary}\label{cor:HK-Jacobi}
Let $f\in C^\infty_0(\R_+)$ and $f^{(2\nu+1)}(0)=0$, $\nu\ge 0$ and consider the kernel
$\Lambda_\delta(x, y) = \sum_{k\ge 0} f(\delta\sqrt{\lambda_k}) P_k(x)P_k(y)$,
$0<\delta\le 1$.
Then for any $\sigma>0$ there exists a constant $c_\sigma>0$ such that
\begin{equation}\label{HK-Jacobi-3}
|\Lambda_\delta(x, y)| \le c_\sigma(|B(x, \delta)||B(x, \delta)|\big)^{1/2}
\Big(1+\frac{\rho(x, y)}{\delta}\Big)^{-\sigma},
\end{equation}
where $|B(\cdot, \delta)|$ and $\rho(x, y)$ are as above.
\end{corollary}

This result is more complete than the similar estimate (2.14) in \cite{PX1} (see also \cite{IPX1, IPX2})
which is proved 
under the restriction $\alpha, \beta > -1/2$.


\end{document}